\documentclass[12pt]{article}
\usepackage{bm}
\usepackage{mathptmx}
\usepackage{amsfonts}
\usepackage{latexsym}
\usepackage{graphicx}
\usepackage{color}
\usepackage{ifpdf}
\newtheorem{theorem}{Theorem}[section]

\newtheorem{lemma}[theorem]{Lemma}

\newtheorem{proposition}{Proposition}[section]
\newtheorem{definition}{Definition}[section]
\newtheorem{remark}{Remark}[section]
\usepackage{amsfonts,amssymb,eucal,amsmath}
\title{\bf On Periodic orbits of  the Planar $N$-body Problem   \thanks{Supported by NSFC(No.11701464)}}

\author{{ Xiang Yu\footnote{Email:yuxiang@swufe.edu.cn, xiang.zhiy@gmail.com}} \\
\small \it School of Economic and Mathematics, Southwestern
University of Finance and Economics, \\
\small \it Chengdu 611130, China}

\date{}

\begin{document}
\maketitle

\begin{abstract}
By introducing a new coordinate system, we prove that there are abundant new periodic orbits near relative equilibrium solutions of the $N$-body problem. We consider only  Lagrange relative equilibrium of the three-body problem and Euler-Moulton relative equilibrium of the $N$-body problem, although we believe that there are similar results for general relative equilibrium solutions.  All of these periodic orbits lie on a $2d$-dimensional central manifold of  the planar $N$-body problem.  Besides $d$ one parameter family of periodic orbits  which are well known as Lyapunov's orbits or Weinstein's orbits, we further prove that periodic
orbits are unexpectedly abundant: generically the relative measure of the closure of the set of periodic orbits near relative equilibrium solutions on the $2d$-dimensional central manifold is close to 1. These abundant  periodic
orbits are named Conley-Zender's orbits, since  to find them is based on an extended result of Conley and Zender on the local existence result for periodic orbits near an elliptic equilibrium point of a Hamiltonian. In particular, the results provide some evidences to support the well known claim of Poincar\'{e} on the conjecture of  periodic orbits of the $N$-body problem. \\


{\bf 2010AMS Subject Classification} {70F10 \and 34C25 \and 70K42 \and 70K43 \and 70K45}.

\end{abstract}

\section{Introduction}
\ \ \ \
We consider $N $ particles with positive mass moving in an Euclidean space ${\mathbb{R}}^2$ interacting under
the laws of universal gravitation. Let the $k$-th particle have mass $m_k$ and position
$\mathbf{r}_k \in {\mathbb{R}}^2$ ($k=1, 2, \cdots, N$), then the equation of motion of the $N$-body problem is written
\begin{equation}\label{eq:Newton's equation1}
m_k \ddot{\mathbf{r}}_k =\sum_{1 \leq j \leq N, j \neq k} \frac{m_km_j(\mathbf{r}_j-\mathbf{r}_k)}{|\mathbf{r}_j-\mathbf{r}_k|^3}, ~~~~~~~~~~~~~~~k=1, 2, \cdots, N.
\end{equation}
where $|\cdot|$  denotes the Euclidean norm in ${\mathbb{R}}^2$. Since these equations
are invariant by translation, we can assume that the center of mass stays at the origin.

The importance of periodic orbits (or solutions) are especially emphasized by Poincar\'{e}  in his celebrated work \emph{Les m{\'e}thodes nouvelles de la m{\'e}canique c{\'e}leste} \cite{poincare1892methodes}. As a matter of fact, Poincar\'{e} wrote extensively on periodic orbits in  \cite{poincare1892methodes}.

In particular, on periodic orbits of the
three-body problem, Poincar\'{e} wrote
: ``\ldots In fact, there is a zero probability for the initial conditions of the motion to
be precisely those corresponding to a periodic solution. However, it can happen
that they differ very little from them, and this takes place precisely in the case
where the old methods are no longer applicable. We can then advantageously take
the periodic solution as first approximation, as intermediate orbit, to use Gyld\'{e}n's language.

There is even more: here is a fact which I have not been able to demonstrate
rigorously, but which seems very probable to me, nevertheless.

Given equations of the form defined in art. 13 and any particular solution of
these equations, we can always find a periodic solution (whose period, it is true,
is very long), such that the difference between the two solutions is as small as
we wish, during as long a time as we wish. In addition, these periodic solutions
are so valuable for us because they are, so to say, the only breach by which we
may attempt to enter an area heretofore deemed inaccessible."  (\cite{poincare1892methodes},  ch. 3, a. 36).

This Poincar\'{e}'s conjecture, that is, periodic orbits of the $N$-body problem are dense,  was often quoted by Birkhoff as a main motivation for his works on fixed point theorems and related topics (see \cite{birkhoff1927dynamical,Birkhoff1915The}). In particular, Birkhoff also pointed out that ``The question of the existence of periodic orbits is of great importance" (\cite{Birkhoff1915The},  p. 267). However, the Poincar\'{e}'s conjecture is  far from settled. According to what I know, only a weakened version of Poincar\'{e}'s  conjecture for the restricted three-body problem has been established by G\'{o}mez and  Llibre \cite{gomez1981a}.

Poincar\'{e} suggested two approaches to establish the existence of periodic
orbits in the $N$-body problem. One approach is global and the other is local.

The first approach is purely variational:  since a solution of  the  $N$-body problem is a critical
point of the corresponding Lagrangian action, an action
minimizer should be a classical solution. However, due to the potential of the $N$-body problem is singular at collision configurations, the main
problem involved in variational minimizations is that collisions could occur for an action
minimizer and this may prevent an action
minimizer from being a true
solution.  As a result,  the successful time of obtaining periodic
orbits by variational methods in the $N$-body problem was  much later than that in general Hamiltonian systems. Indeed, as recently as 2000, Chenciner
and Montgomery got the first and well known result on periodic
orbits of the $N$-body problem  by variational methods \cite{Chenciner2000a}, then  variational approaches  are extensively exploited to study the $N$-body problem
with numerous work appearing in recent years,   please see \cite[etc]{chenciner2002action,
ferrario2004existence,zhang2004nonplanar,Chen2008}
and the references therein.

The second approach is based
on the principle of analytic continuation of periodic orbits. The method is effective and feasible all the time.
Since the work  of Poincar\'{e}, by the continuation method,  there is a good deal of literature  on the existence and nature of periodic
orbits of the $N$-body problem, especially the restricted three-body problem. Please see \cite[etc]{poincare1892methodes,siegel1971lectures,
Meyer2009Introduction}
and the references therein.

Besides, the fixed-point method developed by  Birkhoff (see \cite{Birkhoff1915The}) also goes
back to Poincar\'{e}. And many different methods have been used to establish the existence of periodic
orbits in the $N$-body problem, such as:  the Lagrangian manifold intersection
theory (see Weinstein  \cite{weinstein1973normal}), the method of majorants (see Lyapunov \cite{liapounoff1907probleme} and Siegel \cite{siegel1950periodische}),
and so on.

For the fixed-point method, let's  take note of the remarkable work of Moser \cite{moser1977proof} which proves and extends the well known the Birkhoff-Lewis fixed point theorem. This work can be used to establish the existence of infinitely
many periodic
orbits  near an equilibrium point of a Hamiltonian and is made available for this paper to establish the existence of infinitely
many periodic
orbits in the $N$-body problem. Nevertheless, in this paper we highlight the  method of Conley and Zender on the local existence result for periodic orbits near an elliptic equilibrium point of a Hamiltonian (see \cite{conley1983index}), and we will establish the existence of abundant periodic
orbits in the $N$-body problem by virtue of the  method of Conley and Zender. This method is  based upon KAM theory and the Birkhoff-Lewis fixed point theorem. Indeed Conley and Zender observed that the invariant tori constructed by the KAM Theorem lie in the closure of the set of periodic orbits whose existence is ensured by Birkhoff-Lewis fixed point theorem.

However, we remark that  the orbits Poincar\'{e} investigates are
not periodic in the standard sense \cite{poincare1892methodes,siegel1971lectures,
Meyer2009Introduction,montgomery1998n}. Indeed, the concept of periodic
orbits is related to the selected reference frame physically. For example, the motion of an object is periodic in an inertial frame, but probably not  in another inertial frame.
Thus periodic orbits in the $N$-body problem are \emph{relative periodic orbits} traditionally:  a  motion $\mathbf{r}(t)$
will be called \emph{relative periodic} with period $T$ if there is an orthogonal transformation $A$  such that $\mathbf{r}(t+T)=A \mathbf{r}(t)$. Or in Meyer's words, relative periodic orbits are not necessarily periodic in fixed space, but are
periodic when the rotational symmetry is eliminated, i.e.,
in the `reduced space' (see \cite{Meyer2009Introduction}).

In this paper, we follow  Poincar\'{e} and Siegel and others  to find (relative) periodic orbits of  the $N$-body problem. We will show that there are abundant periodic
orbits near relative equilibrium solutions of  the $N$-body problem. Here we consider only  Lagrange relative equilibria of the three-body problem and Euler-Moulton relative equilibria of the $N$-body problem, although we believe that there are similar results for general relative equilibrium solutions.

It's well known that generically there are three one parameter family of periodic orbits near Lagrange relative equilibrium of the three-body problem for appropriate masses by the Lyapunov's Center Theorem (see the following Theorem \ref{Lyapunov}), more specifically, there are one trivial family of periodic
orbits which are just Keplerian elliptic orbits generated by Lagrange configuration and  two one parameter family of Lyapunov's periodic orbits near Lagrange relative equilibrium (see \cite{siegel1971lectures}). And there are $N-1$ one parameter family of periodic orbits near Euler-Moulton relative equilibrium of the $N$-body problem by the Weinstein's Theorem (see the following Theorem \ref{weinstein}), more specifically,  there are one trivial family of periodic
orbits which are just Keplerian elliptic orbits generated by Euler-Moulton configuration and  $N-2$ one parameter family of Weinstein's periodic orbits near Euler-Moulton relative equilibrium (see \cite{Moeckel1994}).

Based upon a new coordinate system,  we explore periodic orbits near relative equilibrium solutions. In the coordinate system, the degeneracy of the equations of motion according to intrinsic symmetrical characteristic of the $N$-body problem can easily be reduced. As a matter of fact,  we obtain a system of practical equations of motion which are suitable for describing  the motion of  particles in a neighbourhood of relative equilibrium solutions by virtue of the coordinate system. At first, the coordinate system is originated in the work \cite{yu2019problem} on  the problem of the infinite spin, since we found that the coordinate system is suitable for collision
orbits. It is an exciting fact to see that the coordinate system is also useful to find periodic orbits of the $N$-body problem.

First, we rediscover the above Lyapunov's orbits or Weinstein's orbits.
All of these periodic orbits lie on  central manifolds of  the planar $N$-body problem. If adding one trivial family of periodic orbits, we say that there are generically $d$ one parameter family of periodic orbits on a $2d$-dimensional central manifold of the planar $N$-body problem. More specifically, for Lagrange relative equilibria, when the masses $m_1,m_2,m_3$ of particles satisfy the condition
\begin{equation}
m_1 m_2+ m_3 m_2+ m_1 m_3 < \frac{(m_1+m_2+m_3)^2}{27},\nonumber
\end{equation}
$d=3$; for Euler-Moulton relative equilibria of the $N$-body problem,
$d=N-1$.

By virtue of Conley-Zender's Theorems (see the following Theorems \ref{Conley and Zender} and \ref{Conley and Zendergeneral}), we further prove that periodic
orbits are abundant:generically the relative measure of the closure of the set of periodic orbits near relative equilibrium solutions on the $2d$-dimensional central manifold is close to 1. The fact echoes the above observation of Poincar\'{e} that periodic orbits of the $N$-body problem are dense, so,  as it were, we obtain some evidences to support Poincar\'{e}'s conjecture.

The paper is structured as follows. In \textbf{Section 2}, we recall some classical aspects of Hamiltonian system, in particular, we give a general statement of Conley-Zender's Theorems on the center manifold. In \textbf{Section 3}, we recall some notations,  and some preliminary results including a moving frame and equations of motion in \cite{yu2019problem,yu2019stability}. In \textbf{Section 4}, we give the Hamiltonian form of equations of motion near relative equilibrium solutions. In \textbf{Section 5}, we  discuss the problem of periodic orbits near five relative equilibrium solutions of the planar three-body problem in detail.  In \textbf{Section 6}, we  discuss the problem of periodic orbits near Euler-Moulton relative equilibria of the planar $N$-body problem by the perturbed method and the inductive method. It's noteworthy that the troublesome task in the discussion of this paper is obtaining the  Birkhoff normal of the problem.

\section{Classical Results of Hamiltonian System}
\indent\par

In this section, let's recall some aspects of Hamiltonian system that will be needed later.

We consider a $\mathrm{C}^{l}$-smooth ($2\leq l \leq \omega$ and $\mathrm{C}^{\omega}$ means  analytic)  Hamiltonian system, with $n$ degrees of freedom,
having the origin as an equilibrium point:
\begin{equation}\label{Hamiltonian system}
H(p,q) =  H_2 (p,q)+\cdots + H_m (p,q)+ O(\|(p,q)\|^{m+1}),
\end{equation}
where $H_k$ is a homogeneous polynomial of degree $k$ in $(p,q)$ for every $2\leq k \leq m$; hereafter $O(\|(p,q)\|^{m+1})$  stands for the terms which vanish at the origin  together with all its partial derivatives of
the first $m$-th order.

We will not specify the order of differentiability explicitly in each
point. From the context it will be clear what regularity is needed in a certain step.


To establish the existence of periodic
orbits in a neighbourhood of an equilibrium point of a Hamiltonian, there are the well known the Lyapunov's Center Theorem and the Weinstein's Theorem:
\begin{theorem}\label{Lyapunov}\emph{(\textbf{Lyapunov}\cite{liapounoff1907probleme,siegel1971lectures})}
Consider a $\mathrm{C}^{\omega}$-smooth  Hamiltonian $H$, suppose that the $2n$ eigenvalues $\lambda_1,\cdots,\lambda_n,-\lambda_1,\cdots,-\lambda_n$ of the
quadratic part $H_2$ of the Hamiltonian are all distinct. Let $\lambda_1$ be purely imaginary, and assume that none of $n-1$ quotients $\frac{\lambda_2}{\lambda_1},\cdots,\frac{\lambda_n}{\lambda_1}$  is an integer. Then there exists a family of real periodic
orbits to the Hamiltonian system which depend analytically on one
real parameter $\epsilon$, with $\epsilon=0$ corresponding to the equilibrium solution.
The period $T(\epsilon)$, likewise, is analytic in $\epsilon$ and, moreover, $T(0)=\frac{2\pi}{|\lambda_1|}$.
\end{theorem}
\begin{theorem}\label{weinstein}\emph{(\textbf{Weinstein}\cite{weinstein1973lagrangian,weinstein1973normal})}
Consider a $\mathrm{C}^{2}$-smooth  Hamiltonian $H$, if the
quadratic part $H_2$ of the Hamiltonian is positive definite, then for sufficiently small $\epsilon>0$,  there are at least $n$ geometrically distinct periodic
orbits on energy  surface $H=\epsilon$ whose periods are close to those
of the  linear system corresponding to $H_2$.
\end{theorem}

Given $m\geq 4$, assume that $\lambda_1,\cdots,\lambda_n$ are nonresonant up to order $m$:
\begin{equation}\label{resonant0}
  k_1 \lambda_1 + \cdots +k_n \lambda_n\neq0, ~~~~~~~~~for~(k_1,   \cdots,  k_n) \in \mathbb{Z}^n  ~such ~that~ 1\leq |k_1| + \cdots +|k_n| \leq m.
\end{equation}
The well-known Birkhoff theorem \cite{birkhoff1927dynamical,siegel1971lectures,Meyer2009Introduction} states that, in some
neighbourhood of the origin, there exists a symplectic change of
variables $(p,q)\mapsto (x,y)$,
near to the identity map, such that  in
the new variables the Hamiltonian function is reduced to a Birkhoff normal
form $\mathcal{H}_{Bm}(x,y)$ of degree $m$ up to terms of degree higher than $m$:
\begin{equation}
H(p,q)=\mathcal{H}(x,y)=\mathcal{H}_{Bm}(x,y)+ O(\|(x,y)\|^{m+1}), \nonumber
\end{equation}
here
\begin{equation}
\mathcal{H}_{Bm}(x,y) = \sum_{k=1}^{n}\lambda_k x_k y_k +\cdots  \nonumber
\end{equation} is
a polynomial  in symplectic  variables $x,y$ that is actually
a polynomial of degree $[m/2]$ in the variables $x_k y_k$.

First we confine ourself to  the eigenvalues of the
quadratic part $H_2$ of the Hamiltonian are all distinct and purely imaginary.  Then
in suitable symplectic
coordinates, the
quadratic part $H_2$ takes the form
\begin{equation}
H_2 = \sum_{j= 1}^{n} \frac{\omega_j (p_j^2+q_j^2)}{2}.
\end{equation}
Here every $\omega_j$ is
called  a characteristic frequency, and $\varpi=(\omega_1,\cdots,\omega_n)$ is
called  the frequency vector.

Now when the frequency vector $\varpi$ is nonresonant up to order $m$, there exists a symplectic change of
variables $(p,q)\mapsto (x,y)$  in some
neighbourhood of the origin,
near to the identity map, such that  in
the new variables the Hamiltonian function is reduced to a Birkhoff normal
form $\mathcal{H}_{Bm}(\varrho)$ of degree $m$ up to terms of degree higher than $m$:
\begin{equation}
H(p,q)=\mathcal{H}(x,y)=\mathcal{H}_{Bm}(\varrho)+ O(|x|+|y|)^{m+1}. \nonumber
\end{equation}
Here the Birkhoff normal
form $\mathcal{H}_{Bm}(\varrho)$ in symplectic  variables $x,y$ is actually
a polynomial of degree $[m/2]$ in the variables $\varrho_j=\frac{x_j^2+y_j^2}{2}$.

Therefore it's natural  to consider a nearly-integrable Hamiltonian written in action-angle
variables $\varrho,\varphi$ defined by $x_j=\sqrt{2\varrho_j}\cos\varphi_j,y_j=\sqrt{2\varrho_j}\sin\varphi_j$
\begin{equation}\label{Hamiltonian system1}
\mathcal{H}(\varrho,\varphi)=\mathcal{H}_{Bm}(\varrho)+ \mathcal{R}_m(\varrho,\varphi),
\end{equation}
where $\mathcal{R}_m(\varrho,\varphi)=O(\|\varrho\|)^{[m/2]+1}$, here $\|\varrho\| = \max_{1\leq j\leq n}|\varrho_j|$.

Let's recall the concept of $(c,\upsilon)$-Diophantine and the important concept of nondegenerate (see \cite{arnold2006dynamical}):
\begin{definition}
A frequency vector $\varpi$  is said to be $(c,\upsilon)$-Diophantine for some $c, \upsilon>0$ if  we have
\begin{equation}
| (k, \varpi)|\geq \frac{c}{|k|^\upsilon}, ~~~~~~~~~\forall k \in \mathbb{Z}^n ~such ~that ~~|k|\neq 0. \nonumber
\end{equation}
A $(c,\upsilon)$-Diophantine frequency vector $\varpi$  is also said to be strongly
incommensurable.
\end{definition}
\begin{definition}
The Hamiltonian (\ref{Hamiltonian system}) or (\ref{Hamiltonian system1}) is called to be nondegenerate
in a neighbourhood of the equilibrium point if
\begin{equation}
det\left( \frac{\partial^2 \mathcal{H}_{l}}{\partial \varrho^2}|_{\varrho=0}\right)\neq 0.\nonumber
\end{equation}
\end{definition}

Then it is well known that:
\begin{theorem}\label{KAM} \emph{(\textbf{KAM} \cite{arnol1963small,moser1968lectures})}
Consider a sufficiently smooth  Hamiltonian with a nonresonant frequency vector  up to order $4$ in a neighbourhood of an equilibrium point, if the Hamiltonian is nondegenerate, then the Hamiltonian has invariant tori close to the tori of the
linearized system. These tori form a set whose relative measure in the polydisc $\|\varrho\|<\epsilon$ tends to 1 as $\epsilon\rightarrow 0$.
\end{theorem}
\begin{remark}\label{KAM1}
As a matter of fact, on the relative measure of the set of invariant tori in the polydisc $\|\varrho\|<\epsilon$, it's well known that: if the frequency vector of  the Hamiltonian in KAM Theorem \ref{KAM} is  nonresonant up to order $l\geq 4$, then the relative measure of the set of invariant tori in the polydisc $\|\varrho\|<\epsilon$ is
at least $1-O(\epsilon^{\frac{l-3}{4}})$;  if the frequency vector $\varpi$  satisfies the strong
incommensurability condition, i.e., $(c,\upsilon)$-Diophantine condition, then this measure is $1-O(\exp (-\tilde{c}\epsilon^{\frac{-1}{\upsilon+1}}))$ for a positive number $\tilde{c} = const$, now it'necessary that the Hamiltonian is analytic. Please see \cite{arnold2006dynamical}
and the references therein.
\end{remark}

It turns out that the invariant tori constructed by the KAM Theorem lie in the closure of the set of periodic orbits. This is the following well known result of Conley and Zender.
\begin{theorem}\label{Conley and Zender}\emph{(\cite{conley1983index})}
Consider a sufficiently smooth  Hamiltonian with a nonresonant frequency vector  up to order $m\geq 4$ in a neighbourhood of an equilibrium point, if the Hamiltonian is nondegenerate, then the relative measure of the closure of the set of periodic orbits in  the polydisc $\|\varrho\|<\epsilon$ is
at least $1-O(\epsilon^{\frac{m-3}{4}})$.
\end{theorem}
\begin{remark}
As Conley and Zender pointed out, the
periods of these periodic orbits could be very large. Furthermore,
according to the proof of the above theorem, for any given $T>0$, it's easy to see that there must exist a  periodic orbit such that whose  period is an integer multiple of $T$.
\end{remark}

Next, let's consider more general  Hamiltonian $H$ whose quadratic part $H_2$ has   eigenvalues being not purely imaginary. Suppose  the $2n$ eigenvalues $\lambda_1,\cdots,\lambda_n,-\lambda_1,\cdots,-\lambda_n$ of the
quadratic part $H_2$ of the Hamiltonian are of the following form
\begin{equation}
\begin{array}{c}
  Re \lambda_k=0 ~~~~~~~~~~~~~~~~~~1\leq k \leq d\\
  Re \lambda_k \neq 0 ~~~~~~~~~~~~~~~~~~d< k \leq n.
\end{array}\nonumber
\end{equation}

We introduce the notation $\lambda_k = \mathbf{i} \omega_k, 1\leq k \leq d$ for the purely imaginary
eigenvalues, and assume that  $\omega_1,\cdots,\omega_{d}$ are nonresonant up to order $4$:
\begin{equation}\label{resonant1}
  k_1 \omega_1 + \cdots +k_{d} \omega_{d} \neq0, ~~~~~~~~~for~(k_1,   \cdots,  k_{d}) \in \mathbb{Z}^{d}  ~such ~that~ 1\leq |k_1| + \cdots +|k_{d}| \leq 4;
\end{equation}
Then the well-known Birkhoff theorem \cite{birkhoff1927dynamical,siegel1971lectures,Meyer2009Introduction} states that, in some
neighbourhood of the origin, there exists a symplectic change of
variables
\begin{equation}
(p_1,\cdots,p_{d},p_{d+1},\cdots,p_{n},q_1,\cdots,q_{d},q_{d+1},\cdots,q_{n})\mapsto (x,X,y,Y) \nonumber
\end{equation}
near to the identity map, such that  in
the new variables the Hamiltonian function is reduced to a Birkhoff normal
form:
\begin{equation}\label{reducednormalform}
\begin{array}{c}
  H(p,q)=\mathcal{H}(x,X,y,Y)= \\
  \sum_{k=1}^{d}\omega_k \varrho_k + X^\top \Omega Y +\frac{1}{2}\sum_{j,k=1}^{d}\omega_{jk} \varrho_j \varrho_k + O(\|(x,y)\|^{5}+\|(X,Y)\|^{3}+\|(x,y)\|\|(X,Y)\|^{2}),
\end{array}
\end{equation}
where $\Omega$ is an $(n-d)\times(n-d)$ matrix, $\varrho_k=\frac{x_k^2+y_k^2}{2}$, and $(\omega_{jk})$ is a real symmetric
$d\times d$ matrix.

For the Hamiltonian system
\begin{equation}\label{reducedHamiltoniansystem}
\left\{
             \begin{array}{lr}
             \dot{x}=-\frac{\partial\mathcal{H}(x,X,y,Y)}{\partial y} &  \\
             \dot{y}=\frac{\partial\mathcal{H}(x,X,y,Y)}{\partial x} &\\
             \dot{X}=-\frac{\partial\mathcal{H}(x,X,y,Y)}{\partial Y} &  \\
             \dot{Y}=\frac{\partial\mathcal{H}(x,X,y,Y)}{\partial X} &
             \end{array}
\right.
\end{equation}
it's obvious that there exists an $2d$-dimensional invariant center manifold. The question is that whether the flow on the center manifold
is again described by a reduced Hamiltonian system? As Mielke pointed out in \cite{mielke2006hamiltonian}, the answer of this question is yes, and ``although this fact seems to be
well-known in the realm of Hamiltonian theory, one hardly finds references to this".
As a general result, Mielke \cite{mielke2006hamiltonian}  showed that there are always symplectic
coordinates $(x,X,y,Y)$ such that the center manifold is given by $(X,Y)\equiv0$. In these
coordinates the reduction is trivial: The reduced Hamiltonian is $\mathcal{H}(x,y)=\mathcal{H}(x,0,y,0)$
and $(x,y)$ are symplectic coordinates on the center manifold. This result goes
back to Moser \cite{moser1977proof}.
\begin{lemma}\label{Mielke}\emph{(\cite{mielke2006hamiltonian})}
Let $\mathcal{H}(x,X,y,Y)\in \mathrm{C}^{l}(\mathbb{R}^{2n},\mathbb{R})$ be a Hamiltonian
having an equilibrium in the origin with a (non-trivial) center manifold. Then for any $m<l$  there exists an analytical symplectic transformation $(x,X,y,Y) = \psi (u,\textsc{U},v,\textsc{V})$ such that the
center manifold is given by $(\textsc{U},\textsc{V})=F^{c}(u,v)=O(\|(u,v)\|^{m})$. Taking $(u,v)$ as symplectic
coordinates on the center manifold and using $\widetilde{\mathcal{H}}(u,v)=\mathcal{H}(\psi (u,0,v,0))$ as Hamiltonian
gives the correct terms up to order $2m$ of the true reduced Hamiltonian system.
\end{lemma}
As a result, we have
\begin{theorem}\label{reduced Hamiltonian}\emph{(\cite{mielke2006hamiltonian,Shilnikov})}
Consider a $\mathrm{C}^{l}$-smooth ($4\leq l \leq \omega$) Hamiltonian $\mathcal{H}(x,X,y,Y)$ (\ref{reducednormalform}) in a neighbourhood of the origin. Assume $d<n$, then  there exists a $2d$-dimensional invariant center manifold $\mathcal{W}^{c}_{loc}: (X,Y)=F^{c}(x,y)=O(\|(x,y)\|^{3})$ of class $C^{l-1}$, and we can take $(x,y)$ as symplectic
coordinates on the center manifold and using $\widetilde{\mathcal{H}}(x,y)=\mathcal{H}(x,0,y,0)$ as Hamiltonian
gives the correct terms up to order $6$ of the true reduced Hamiltonian system.
\end{theorem}
\begin{remark}
In general, we cannot claim
that the center manifold is analytical for analytical Hamiltonian systems. However, there is some hope that center manifolds for analytic Hamiltonian
systems are again analytic. For example, if the center manifold is two-dimensional and is filled with
periodic solutions of bounded period the the results is true, please see \cite{moser1958generalization,aulbach1985classical}. Indeed, one conjectured \cite{mielke2006hamiltonian}: ``If the center manifold of an analytic Hamiltonian system is completely filled
with bounded solutions for more restrictive: the reduced Hamiltonian has
a positive (negative) definite second derivative at the fixed point), then the
center manifold is analytic."
\end{remark}

Let's finish this section with the following general statement of Conley-Zender's Theorems on the center manifold which is an obvious corollary of Theorem \ref{Conley and Zender} and Theorem \ref{reduced Hamiltonian}.
\begin{theorem}\label{Conley and Zendergeneral}
Consider the  Hamiltonian $\mathcal{H}(x,X,y,Y)$ (\ref{reducednormalform}) in a neighbourhood of the origin, if $det(\omega_{jk})\neq0$, then the Hamiltonian system (\ref{reducedHamiltoniansystem}) has abundant periodic orbits on the $2d$-dimensional center manifold, indeed, the relative measure of the closure of the set of periodic orbits in  the polydisc $\|\varrho\|<\epsilon$ on the $2d$-dimensional center manifold is
at least $1-O(\epsilon^{\frac{1}{4}})$; if $\omega_1,\cdots,\omega_{d}$ are nonresonant up to order $m\geq 4$, then the relative measure  is
at least $1-O(\epsilon^{\frac{m-3}{4}})$.
\end{theorem}

\section{Preliminaries}
\label{Preliminaries}
\indent\par
In this section we recall some notations and definitions that have essentially been introduced in \cite{yu2019problem,yu2019stability}, but here some of them will be restated in a slightly more general form.     

Let $\mathbb{R}^2$: $ (\mathbb{R}^2)^N = \{ \mathbf{r} = (\mathbf{r}_1,  \cdots, \mathbf{r}_N):\mathbf{r}_j \in \mathbb{R}^2, j = 1,  \cdots, N \}$ denote the space of configurations for $N$ point particles in Euclidean space. Here $(\mathbb{R}^2)^N $ is considered as a column space. Then $\mathbf{r}\in (\mathbb{R}^2)^N $ can be written as
\begin{displaymath}
\mathbf{r} = (\xi_1, \xi_2, \cdots, \xi_{2N})^\top,
\end{displaymath}
here ``$\top$" denotes transposition of matrix. It's also true that $\mathbf{r}_j = (\xi_{2j-1},\xi_{2j})^\top$ for $j = 1, 2, \cdots, N$.

It's well known that the equations (\ref{eq:Newton's equation1}) of motion  are the Euler-Lagrange equations of the the action functional $\mathcal{{A}}$
defined by
\begin{equation}\label{action functional}
\mathcal{A}(\mathbf{r}(t)) = \int{ \mathcal{L}(\mathbf{r}(t),\dot{\mathbf{r}}(t)) dt},\nonumber
\end{equation}
where\begin{displaymath}
\mathcal{L}(\mathbf{r},\dot{\mathbf{r}}) = \mathcal{L} = \mathcal{K} + \mathcal{U} = \sum_j \frac{1}{2} m_j |\dot{\mathbf{r}}_j|^2  + \sum_{k<j}{\frac{m_k m_j}{r_{jk}}}
\end{displaymath}
is the Lagrangian function;  $\mathcal{K} (\dot{\mathbf{r}})= \sum_{j=1}^{N} {\frac{1}{2}{m_j |\dot{\mathbf{r}}_j|^2}}$ and $\mathcal{U}(\mathbf{r}) = \sum_{1\le\xi k<j\le\xi N} {\frac{m_k m_j }{r_{jk}}}$ are respectively the kinetic energy function and the opposite of the potential energy (force function).

As usual, consider the total energy, the angular momentum and the moment
of inertia, respectively, defined by
\begin{displaymath}
H(\mathbf{r},\dot{\mathbf{r}})= \mathcal{K}(\dot{\mathbf{r}})- \mathcal{U}(\mathbf{r}),
\end{displaymath}
\begin{displaymath}
\mathcal{J}(\mathbf{r}) = \sum_{j=1}^{N} {m_j {\mathbf{r}}_j \times {\dot{\mathbf{r}}}_j},
\end{displaymath}
\begin{displaymath}
I(\mathbf{r}) = \sum_{j=1}^{N} {{m_j |{\mathbf{r}}_j-{\mathbf{r}}_c|^2}},
\end{displaymath}
where $|{\mathbf{r}}_j| = \sqrt{\xi^2_{2j-1} + \xi^2_{2j}}$, $r_{jk}=|\mathbf{r}_k - \mathbf{r}_j|$, ${\mathbf{r}}_j \times {\dot{\mathbf{r}}}_j = \xi_{2j-1}\dot{\xi}_{2j} - \dot{\xi}_{2j-1}{\xi}_{2j}$, $\mathbf{r}_c = \frac{\sum_{k = 1}^{N} {m_k \vec{r}_k}}{\mathfrak{m}}$ is the center of mass and $\mathfrak{m}=\sum_{k = 1}^{N} m_k$ is the total mass.

Let $\mathfrak{M}$ be the matrix
\begin{displaymath}
diag(m_1,m_1,m_2,m_2, \cdots, m_N, m_N),
\end{displaymath} here ``diag" means diagonalmatrix. Let's introduce a scalar product and a metric  on the space $(\mathbb{R}^2)^N $:
\begin{displaymath}
\langle\mathbf{r},\mathbf{r}\rangle = \sum_{j=1}^{N} {{m_j |{\mathbf{r}}_j|^2}} = (\xi_1, \xi_2, \cdots, \xi_{2N})\mathfrak{M}(\xi_1, \xi_2, \cdots, \xi_{2N})^\top,
\end{displaymath}
\begin{displaymath}
\|\mathbf{r}\|=\sqrt{\langle\mathbf{r},\mathbf{r}\rangle},
\end{displaymath}
then the cartesian space $(\mathbb{R}^2)^N $ is a new Euclidean space.

Observe that the equations (\ref{eq:Newton's equation1}) of motion
are invariant by translation, there is usually an assumption that the center of mass $\mathbf{r}_c$ is at the origin.
We use $\mathcal{X} = \{ \mathbf{r} = (\mathbf{r}_1,\cdots, \mathbf{r}_N)\in (\mathbb{R}^2)^N: \sum_{k = 1}^{N} {m_k \mathbf{r}_k} = 0  \}$ to denote the space of configurations whose center of mass is at the origin; that is,
\begin{displaymath}
\mathcal{X} = \{ \mathbf{r} \in (\mathbb{R}^2)^N: \langle\mathbf{r},\mathcal{E}_1\rangle = 0, \langle\mathbf{r},\mathcal{E}_2\rangle = 0   \},
\end{displaymath}
where
\begin{displaymath}
\mathcal{E}_1 = (1,0,\cdots,1,0)^\top, \mathcal{E}_2 = (0,1,\cdots,0,1)^\top.
\end{displaymath}
Let $\Delta $ be the collision set in $(\mathbb{R}^2)^N $. Then
the set $\mathcal{X} \backslash \Delta$ is the space of collision-free configurations.\\

\begin{definition}
A configuration $\mathbf{r} \in \mathcal{X} \backslash \Delta$ is called a central configuration if there exists a constant $\lambda\in {\mathbb{R}}$ such that
\begin{equation}\label{centralconfiguration}
\sum_{j=1,j\neq k}^N \frac{m_jm_k}{|\mathbf{r}_j-\mathbf{r}_k|^3}(\mathbf{r}_j-\mathbf{r}_k)=-\lambda m_k\mathbf{r}_k,1\leq k\leq N,
\end{equation}
or
\begin{equation}\label{centralconfiguration1}
 \nabla \mathcal{U} = -\lambda\mathbf{r},
\end{equation}
where $\nabla \mathcal{U}$ is the gradient of $\mathcal{U}$ with respect to scalar product $\langle,\rangle$.
\end{definition}
The value of $\lambda$ in (\ref{centralconfiguration})(or (\ref{centralconfiguration1})) is uniquely determined by
\begin{equation}\label{lambda}
\lambda=\frac{\mathcal{U}(\mathbf{r})}{I(\mathbf{r})}.
\end{equation}
It's well known that a central configuration  is just a critical point of the function $I^{\frac{1}{2}}\mathcal{U}$.

\indent\par

Let $\mathbb{O}(2)$ and $S\mathbb{O}(2)$ be the orthogonal group and special orthogonal group of the plane respectively. Set
\begin{displaymath}
A(\theta)=\left(
                 \begin{array}{cc}
                   \cos\theta & -\sin\theta \\
                   \sin\theta & \cos\theta \\
                 \end{array}
               \right)
\in S\mathbb{O}(2).
\end{displaymath}

Given a configuration $\mathbf{r}$, let $\hat{\mathbf{r}}:= \frac{\mathbf{r}}{\|\mathbf{r}\|}$  be the unit vector corresponding to $\mathbf{r}$ henceforth.

For a configuration $\mathbf{r} = (\mathbf{r}_1,\cdots, \mathbf{r}_N)$, let
\begin{displaymath}
\mathbf{r}^\bot = (\mathbf{r}^\bot_1,\cdots, \mathbf{r}^\bot_N)
\end{displaymath}
 denote
\begin{displaymath}
A(\frac{\pi}{2})\mathbf{r} = (A(\frac{\pi}{2})\mathbf{r}_1,\cdots, A(\frac{\pi}{2})\mathbf{r}_N),
\end{displaymath}
 as an illustration, we have $\mathcal{E}_2 = \mathcal{E}^\bot_1$. Similarly, set
\begin{displaymath}
A^\bot(\theta)=A(\frac{\pi}{2})A(\theta)=\frac{d A(\theta)}{d\theta}.
\end{displaymath}

A central configuration $\mathcal{E}_3$ will be called \textbf{nondegenerate}, if  the kernel of the Hessian of $I^{\frac{1}{2}}\mathcal{U}$ evaluated at $\mathcal{E}_3$ is exactly $span\{\mathcal{E}_3, \mathcal{E}_4\}$, where $\mathcal{E}_4 = \mathcal{E}^\bot_3$ is another central configuration.

Given a central configuration
\begin{displaymath}
\mathcal{E}_3 = \mathbf{r} = (\mathbf{r}_1,\cdots, \mathbf{r}_N) = (\xi_1, \xi_2, \cdots, \xi_{2N})^\top,
\end{displaymath}
a straight forward computation shows that the Hessian of $I^{\frac{1}{2}}\mathcal{U}$ evaluated at $\mathcal{E}_3$ is
\begin{displaymath}
I^{\frac{1}{2}} (\lambda \mathfrak{M} + \mathfrak{B}) - 3 I^{-\frac{1}{2}} \lambda \mathfrak{M} \mathcal{E}_3 \mathcal{E}_3^\top \mathfrak{M},
\end{displaymath}
where $\mathfrak{B}$ is the Hessian of $\mathcal{U}$ evaluated at $\mathcal{E}_3$ and can be viewed as an $N\times N$ array of $2 \times 2$ blocks:
\begin{center}
$\mathfrak{B} = \left(
       \begin{array}{ccc}
         B_{11} & \cdots & B_{1N} \\
         \vdots & \ddots & \vdots \\
         B_{N1} & \cdots & B_{NN}\\
       \end{array}
     \right)
$
\end{center}
The off-diagonal blocks are given by:
\begin{center}
$B_{jk} = \frac{m_j m_k}{r^3_{jk}}[\mathbb{I}-\frac{3(\mathbf{r}_k - \mathbf{r}_j)(\mathbf{r}_k - \mathbf{r}_j)^\top}{r^2_{jk}}],
$
\end{center}
where $\mathbb{I}$ is the identity matrix of order 2. However, as a matter of notational convenience, the \textbf{identity matrix} of any order will always be denoted by $\mathbb{I}$, and the order of $\mathbb{I}$ can be determined according context.
The diagonal blocks are given by:
\begin{displaymath}
B_{kk} = -\sum_{1\leq j\leq N, j\neq k} B_{jk}.
\end{displaymath}

Let us investigate the matrix
\begin{displaymath}
\mathfrak{D}:= I^{\frac{1}{2}} (\lambda \mathbb{I} + \mathfrak{M}^{-1} \mathfrak{B}) - 3 I^{-\frac{1}{2}} \lambda \mathcal{E}_3 \mathcal{E}_3^\top \mathfrak{M}
\end{displaymath}
which can be viewed as the linearization of the gradient $\nabla \mathcal{U}$  at
the central configuration $\mathcal{E}_3$. Since the matrix $\mathfrak{D}$ is symmetric linear mapping with respect to the scalar product $\langle,\rangle$, there are $2N$ orthogonal eigenvectors of $\mathfrak{D}$  with respect to the scalar product $\langle,\rangle$.
It's easy to see that:
\begin{displaymath}
\mathfrak{D} \mathcal{E}_1 = I^{\frac{1}{2}} \lambda \mathcal{E}_1,
\end{displaymath}
\begin{displaymath}
\mathfrak{D} \mathcal{E}_2 = I^{\frac{1}{2}} \lambda \mathcal{E}_2,
\end{displaymath}
\begin{displaymath}
\mathfrak{D} \mathcal{E}_3 = 0,
\end{displaymath}
\begin{displaymath}
\mathfrak{D} \mathcal{E}_4 = 0.
\end{displaymath}
Therefore an orthogonal basis $\{\mathcal{E}_1, \mathcal{E}_2, \mathcal{E}_3, \mathcal{E}_4, \cdots, \mathcal{E}_{2N}\}$ can be chosen as  $2N$ orthogonal eigenvectors of $\mathfrak{D}$, that is,
\begin{displaymath}
 \mathfrak{D} (\mathcal{E}_1,  \cdots, \mathcal{E}_{2N}) = (\mathcal{E}_1,  \cdots, \mathcal{E}_{2N}) diag(\lambda_1,  \cdots, \lambda_{2N}),
\end{displaymath}
where $\lambda_j \in \mathbb{R}$ is the eigenvalue of $\mathfrak{D}$ corresponding to $\mathcal{E}_k$ ($k=1, 2, \cdots, 2N$), in addition, $\lambda_1=\lambda_2=I^{\frac{1}{2}} \lambda$, $\lambda_3=\lambda_4=0$.

Suppose
\begin{displaymath}
(\mathcal{E}_1, \mathcal{E}_2,\mathcal{E}_3, \mathcal{E}_4, \cdots, \mathcal{E}_{2N})^\top \mathfrak{M} (\mathcal{E}_1, \mathcal{E}_2, \mathcal{E}_3, \mathcal{E}_4, \cdots, \mathcal{E}_{2N}) = diag(g_{1},g_{2},g_{3},g_{4},  \cdots, g_{2N}),
\end{displaymath}
then
\begin{displaymath}
\begin{array}{c}
  g_{1}=g_{2}=  \sum_{k = 1}^{N} m_k = \mathfrak{m}, \\
  g_{3}=g_{4}=\|\mathcal{E}_3\|^2 =I(\mathcal{E}_3)=I,
\end{array}
\end{displaymath}
and $\mathcal{E}_k$ is a unit vector if and only if $g_{k}=I(\mathcal{E}_k)=\|\mathcal{E}_k\|^2=1$.

It follows that
\begin{equation}\label{Hessian00}
(\mathcal{E}_1,\cdots, \mathcal{E}_{2N})^{\top} (I^{\frac{1}{2}} (\lambda \mathfrak{M} + \mathfrak{B}) - 3 I^{-\frac{1}{2}} \lambda \mathfrak{M} \mathcal{E}_3 \mathcal{E}_3^\top \mathfrak{M}) (\mathcal{E}_1,\cdots, \mathcal{E}_{2N}) = diag(g_{1}\lambda_1, \cdots, g_{2N}\lambda_{2N}),
\end{equation}
\begin{equation}\label{Hessian11}
(\mathcal{E}_1,\cdots, \mathcal{E}_{2N})^{\top}  (\lambda \mathfrak{M} + \mathfrak{B})  (\mathcal{E}_1,\cdots, \mathcal{E}_{2N})
= diag( g_{1}\lambda, g_{1}\lambda,3 g_{3} \lambda, 0,  \frac{g_{5}\lambda_5}{\sqrt{g_{3}}}, \cdots,  \frac{g_{2N}\lambda_{2N}}{\sqrt{g_{3}}}).
\end{equation}

It is noteworthy that the subspaces $span\{\mathcal{E}_1, \mathcal{E}_2\}$, $span\{\mathcal{E}_3,\mathcal{E}_4\}$ and $span\{\mathcal{E}_5, \cdots, \mathcal{E}_{2N}\}$ of the space $(\mathbb{R}^2)^N $ are  invariant under the action of the transformation
\begin{displaymath}
\rho A\mathbf{r}=(\rho A\mathbf{r}_1, \rho A\mathbf{r}_2, \cdots, \rho A\mathbf{r}_N),\end{displaymath}
 where $A \in \mathbb{O}(2)$ (or $S\mathbb{O}(2)$) and $\rho>0$.

\subsection{Moving Frame}
\indent\par
Then, as in \cite{yu2019problem}, let's give a moving frame to describe the motion of the particles near a  relative equilibrium solution of the Newtonian $N$-body problem effectively.

For any configuration $\mathbf{r}\in \mathcal{X} \backslash  span\{\mathcal{E}_5, \cdots, \mathcal{E}_{2N}\}$, it's easy to see that  there exists a unique point $A(\theta(\mathbf{r}))\mathcal{E}_3$ on $\mathbf{S}$ such that
 \begin{equation}
\|A(\theta(\mathbf{r}))\mathcal{E}_3-\mathbf{r}\|= min_{\theta \in \mathbb{R}} \|A(\theta)\mathcal{E}_3 - \mathbf{r}\|,\nonumber
\end{equation}
where $\mathbf{S}= \{A(\theta)\mathcal{E}_3: \theta \in \mathbb{R}\}$ is a  circle in the space $(\mathbb{R}^2)^N $ with the origin as the center.

$\theta$ in the point $A(\theta(\mathbf{r}))$  can be continuously determined as a continuous function of the independent variable $\mathbf{r}$.

Set $\Xi_3 = A(\theta)\mathcal{E}_3, \Xi_4 = A(\theta)\mathcal{E}_4, \cdots, \Xi_{2N} = A(\theta)\mathcal{E}_{2N}$, then $\{\Xi_3, \Xi_{4}, \cdots, \Xi_{2N}\}$ is an orthogonal basis of $\mathcal{X}$, and
\begin{displaymath}
span\{\Xi_3, \Xi_4\}=span\{\mathcal{E}_3, \mathcal{E}_4\},~~~~~~span\{\Xi_5, \cdots, \Xi_{2N}\}=span\{\mathcal{E}_5, \cdots, \mathcal{E}_{2N}\}.
\end{displaymath}
$\{\Xi_3, \Xi_{4}, \cdots, \Xi_{2N}\}$  is a class of \textbf{moving frame} which is suitable for describing  the motion of  particles in a neighbourhood of relative equilibrium solutions.

Set $r = \|\mathbf{r}\|$, then $\mathbf{r}=r \hat{\mathbf{r}}$. In the moving frame, $\hat{\mathbf{r}}$ can be written as $\hat{\mathbf{r}} =  \sum_{k = 3}^{2N} {z_k \Xi_k}$. It's easy to see that $z_4=0$ and
\begin{equation}\label{z3}
z_3 = \sqrt{\frac{1 - \sum_{k = 5}^{2N} g_k z^2_k}{g_{3}}}.
\end{equation}
Then the total set of the variables $r, \theta, z_5, \cdots, z_{2N}$ can be thought as the  coordinates of $\mathbf{r}\in \mathcal{X} \backslash  span\{\mathcal{E}_5, \cdots, \mathcal{E}_{2N}\}$ in the moving frame.

Note that
\begin{displaymath}
min_{\theta \in \mathbb{R}} \|A(\theta)\mathcal{E}_3 - \mathbf{r}\| \geq \sqrt{g_{3}}
\end{displaymath}
if  $\mathbf{r}\in   span\{\mathcal{E}_5, \cdots, \mathcal{E}_{2N}\}$. Consequently, if
\begin{displaymath}
min_{\theta \in \mathbb{R}} \|A(\theta)\mathcal{E}_3 - \mathbf{r}\| < \sqrt{g_{3}},
\end{displaymath}
then  $\mathbf{r}\in  \mathcal{X} \backslash span\{\mathcal{E}_5, \cdots, \mathcal{E}_{2N}\}$. Therefore we have legitimate rights to use the coordinates $ (r, \theta, z_5, \cdots, z_{2N})$ in a neighbourhood of a relative equilibrium.

\subsection{Equations of Motion in the Moving Frame}
\indent\par

As in \cite{yu2019problem,yu2019stability}, we can write the equations of motion in the above given coordinates.

First, by virtue of the coordinates $ r, \theta, z_5, \cdots, z_{2N}$,  the kinetic energy and force function can be respectively rewritten as
\begin{equation}\label{kineticenergy}
\mathcal{K}(\mathbf{r}) = \frac{\dot{r}^2}{2}+ \frac{r^2}{2} (g_3 \dot{z}^2_3+\sum_{k=5}^{2N}g_k \dot{z}^2_k + 2\dot{\theta}\sum_{j,k=5}^{2N}\langle\mathcal{E}_j,\mathcal{E}^\bot_k\rangle\dot{z}_j {z}_k+ \dot{\theta}^2),\nonumber
\end{equation}
\begin{equation}\label{force function}
\mathcal{U}(\mathbf{r}) = \frac{\mathcal{U}(z_3 \mathcal{E}_3+\sum_{j = 5}^{2N} {z_j \mathcal{E}_j})}{r}.\nonumber
\end{equation}

It's noteworthy that the variable $\theta$ is not involved in the function $\mathcal{U}(z_3 \mathcal{E}_3+\sum_{j = 5}^{2N} {z_j \mathcal{E}_j})$, this is one of the main advantage of introducing the moving frame.  In particular, the variable $\theta$ is not involved in the Lagrangian $\mathcal{L}$, that is, the variable $\theta$ is an ignorable coordinate.

Since $\mathcal{U}(z_3 \mathcal{E}_3+\sum_{j = 5}^{2N} {z_j \mathcal{E}_j})$ only contains the variables $z_j$ ($j=5, \cdots, 2N$), we will simply
write it as  $U(z)$ henceforth, that is,
\begin{equation}
U(z) = \mathcal{U}(z_3 \mathcal{E}_3+\sum_{j = 5}^{2N} {z_j \mathcal{E}_j})=\mathcal{U}((1 - \sum_{k = 5}^{2N} g_k z^2_k)^{\frac{1}{2}} \widehat{\mathcal{E}}_3+\sum_{j = 5}^{2N} {z_j \mathcal{E}_j}).\nonumber
\end{equation}

Consequently, we can expand $U(z)$ as
\begin{eqnarray}
&U(z) =\mathcal{U}(\hat{\mathcal{E}}_3) + \sum_{k = 5}^{2N} d\mathcal{U}|_{\hat{\mathcal{E}}_3}(\mathcal{E}_k)z_k + d\mathcal{U}|_{\hat{\mathcal{E}}_3}(\hat{\mathcal{E}}_3)\left(\frac{-\sum_{j = 5}^{2N} g_j z^2_j}{2}-\frac{(\sum_{j = 5}^{2N} g_j z^2_j)^2}{8}\right)+\nonumber\\
&  \frac{1}{2}[\sum_{j,k = 5}^{2N} d^2 \mathcal{U}|_{\hat{\mathcal{E}}_3}(\mathcal{E}_j,\mathcal{E}_k) z_j z_k - \sum_{k =5}^{2N} d^2 \mathcal{U}|_{\hat{\mathcal{E}}_3}(\hat{\mathcal{E}}_3,\mathcal{E}_k)z_k\sum_{j = 5}^{2N} g_j z^2_j + d^2 \mathcal{U}|_{\hat{\mathcal{E}}_3}(\hat{\mathcal{E}}_3,\hat{\mathcal{E}}_3)(\frac{\sum_{j = 5}^{2N} g_j z^2_j}{2})^2] \nonumber\\
&+\frac{1}{3!}[\sum_{i,j,k = 5}^{2N} d^3 \mathcal{U}|_{\hat{\mathcal{E}}_3}(\mathcal{E}_i,\mathcal{E}_j,\mathcal{E}_k)z_i z_j z_k+3\sum_{j,k = 5}^{2N} d^3 \mathcal{U}|_{\hat{\mathcal{E}}_3}(\hat{\mathcal{E}}_3,\mathcal{E}_j,\mathcal{E}_k)(\frac{-\sum_{i = 5}^{2N} g_i z^2_i}{2}) z_j z_k] \nonumber\\
&+\frac{1}{4!}[\sum_{h,i,j,k = 5}^{2N} d^4 \mathcal{U}|_{\hat{\mathcal{E}}_3}(\mathcal{E}_h,\mathcal{E}_i,\mathcal{E}_j,\mathcal{E}_k)z_h z_i z_j z_k] +\cdots\nonumber
\end{eqnarray}
where ``$\cdots$" denotes higher order terms of $z_j$ ($j=5, \cdots, 2N$), and $d\mathcal{U}|_{\hat{\mathcal{E}}_3}, d^2\mathcal{U}|_{\hat{\mathcal{E}}_3}, d^3\mathcal{U}|_{\hat{\mathcal{E}}_3},d^4\mathcal{U}|_{\hat{\mathcal{E}}_3} $ denote respectively  the differential,  second order differential, third order differential, fourth order differential of $\mathcal{U}$ at $\hat{\mathcal{E}}_3$.

Then it follows from (\ref{centralconfiguration1}) (\ref{lambda})  (\ref{Hessian11}) (\ref{z3})  that
\begin{eqnarray}\label{expandforcefunction}
&U(z)=  g_3^{\frac{3}{2}}\lambda + \frac{1}{2} \sum_{k = 5}^{2N} g_3 g_k\lambda_k z^2_k  +\frac{1}{6}\sum_{i,j,k = 5}^{2N} d^3 \mathcal{U}|_{\hat{\mathcal{E}}_3}(\mathcal{E}_i,\mathcal{E}_j,\mathcal{E}_k)z_i z_j z_k\nonumber\\
& + \frac{3 g_3^{\frac{3}{2}}\lambda}{8} (\sum_{j = 5}^{2N} g_j z^2_j)^2+\frac{3}{4} g_3^{\frac{3}{2}}(\sum_{j = 5}^{2N} g_j z^2_j)\sum_{k = 5}^{2N} (\frac{g_{k}\lambda_k}{\sqrt{g_{3}}} - g_{k}\lambda)  z^2_k \nonumber\\
 &+\frac{1}{4!}\sum_{h,i,j,k = 5}^{2N} d^4 \mathcal{U}|_{\hat{\mathcal{E}}_3}(\mathcal{E}_h,\mathcal{E}_i,\mathcal{E}_j,\mathcal{E}_k)z_h z_i z_j z_k+\cdots.\nonumber
\end{eqnarray}

And the Lagrangian $\mathcal{L}$ can be  rewritten as
\begin{eqnarray}
\mathcal{L}& =& \frac{\dot{r}^2}{2}+ \frac{r^2}{2} (g_3 \dot{z}^2_3+\sum_{k=5}^{2N}g_k \dot{z}^2_k + 2\dot{\theta}\sum_{j,k=5}^{2N}\langle\mathcal{E}_j,\mathcal{E}^\bot_k\rangle\dot{z}_j {z}_k+ \dot{\theta}^2) + \frac{U(z)}{r}.\nonumber
\end{eqnarray}

By computing $\frac{d}{dt}\frac{\partial \mathcal{L}}{\partial \dot{z}_k}-\frac{\partial \mathcal{L}}{\partial z_k} (k= 5,\cdots, 2N), \frac{d}{dt}\frac{\partial \mathcal{L}}{\partial \dot{r}}-\frac{\partial \mathcal{L}}{\partial r}, \frac{d}{dt}\frac{\partial \mathcal{L}}{\partial \dot{\theta}}-\frac{\partial \mathcal{L}}{\partial \theta}$, one can obtain the equations of motion.
But, for the sake of clearer notations, we would better introduce the following transformation to simplify the Lagrangian $\mathcal{L}$:
\begin{equation}\label{transformation1}
\begin{array}{c}
  r= \sqrt{g_3}x_0, \\
  z_k = \frac{x_k}{\sqrt{g_k}}, ~~~~~~~~~~k=3, 5,\cdots, 2N.
\end{array}
\end{equation}
Set $q_{jk}=\langle\widehat{\mathcal{E}}_j,\widehat{\mathcal{E}}^\bot_k\rangle$, then the square matrix $Q:=(q_{jk})_{(2N-4)\times (2N-4)}$ is an anti-symmetric orthogonal matrix.
Set $a_{ijk}=d^3 \mathcal{U}|_{\widehat{\mathcal{E}}_3}(\widehat{\mathcal{E}}_i,\widehat{\mathcal{E}}_j,\widehat{\mathcal{E}}_k)$, then $a_{ijk}$ is symmetric with respect to the subscripts $i,j,k$. Similarly, set $a_{hijk}=d^4 \mathcal{U}|_{\widehat{\mathcal{E}}_3}(\widehat{\mathcal{E}}_h,\widehat{\mathcal{E}}_i,\widehat{\mathcal{E}}_j,\widehat{\mathcal{E}}_k)$, then $a_{hijk}$ is symmetric with respect to the subscripts $h,i,j,k$.

As a result,  the Lagrangian $\mathcal{L}$ becomes
\begin{eqnarray}
\mathcal{L}& =& \frac{g_3\dot{x}^2_0}{2}+ \frac{g_3x_0^2}{2} ( \dot{x}^2_3+\sum_{k=5}^{2N}\dot{x}^2_k + 2\dot{\theta}\sum_{j,k=5}^{2N}q_{jk}\dot{x}_j {x}_k+ \dot{\theta}^2) + \frac{U(x)}{\sqrt{g_3} x_0},
\end{eqnarray}
where
\begin{equation}\label{x3}
x_3 = \sqrt{1 - \sum_{k = 5}^{2N}  x^2_k},
\end{equation}
\begin{eqnarray}
&U(x) = \mathcal{U}(x_3 \widehat{\mathcal{E}}_3+\sum_{j = 5}^{2N} {x_j \widehat{\mathcal{E}}_j})\nonumber\\
&=  g_3^{\frac{3}{2}}\lambda + \frac{1}{2} \sum_{k = 5}^{2N} g_3 \lambda_k x^2_k  +\frac{1}{6}\sum_{i,j,k = 5}^{2N} a_{ijk} x_i x_j x_k+ \frac{3 g_3^{\frac{3}{2}}\lambda}{8} (\sum_{j = 5}^{2N}  x^2_j)^2\nonumber\\
& +\frac{3}{4} g_3^{\frac{3}{2}}(\sum_{j = 5}^{2N}  x^2_j)\sum_{k = 5}^{2N} (\frac{\lambda_k}{\sqrt{g_{3}}} - \lambda)  x^2_k+\frac{1}{4!}\sum_{h,i,j,k = 5}^{2N}a_{hijk} x_h x_i x_j x_k+\cdots.\nonumber
\end{eqnarray}

Furthermore, suppose $\lambda(\widehat{\mathcal{E}}_3)=\lambda^*$ and $\lambda_k(\widehat{\mathcal{E}}_3)={\lambda^*_k}$ ($k\geq 3$), then it's easy to see that
 \begin{displaymath}
\begin{array}{c}
\lambda=\lambda(\mathcal{E}_3)=g_3^{-\frac{3}{2}}\lambda^*\\
\lambda_k= \lambda_k(\mathcal{E}_3)=g_3^{-1} {\lambda^*_k}.
\end{array}
\end{displaymath}
Consequently, the function $U(x)$  becomes
\begin{eqnarray}
&U(x)= \lambda^* + \frac{1}{2} \sum_{k = 5}^{2N} {\lambda^*_k} x^2_k +\frac{1}{6}\sum_{i,j,k = 5}^{2N} a_{ijk} x_i x_j x_k + \frac{3 \lambda^*}{8} (\sum_{k = 5}^{2N}  x^2_k)^2\nonumber\\
&+\frac{3}{4} (\sum_{j = 5}^{2N}  x^2_j)\sum_{k = 5}^{2N} ({\lambda^*_k}- \lambda^*)  x^2_k+\frac{1}{24}\sum_{h,i,j,k = 5}^{2N} a_{hijk} x_h x_i x_j x_k+ \cdots,\nonumber
\end{eqnarray}
or
\begin{eqnarray}\label{Uxsim}
&U(x)= \lambda^* + \frac{1}{2} \sum_{k = 5}^{2N} {\lambda^*_k} x^2_k +\frac{1}{6}\sum_{i,j,k = 5}^{2N} a_{ijk} x_i x_j x_k\nonumber\\
&+\frac{3}{4} (\sum_{j = 5}^{2N}  x^2_j)\sum_{k = 5}^{2N} ({\lambda^*_k}- \frac{\lambda^*}{2})  x^2_k+\frac{1}{24}\sum_{h,i,j,k = 5}^{2N} a_{hijk} x_h x_i x_j x_k+ \cdots.
\end{eqnarray}
By computing $\frac{d}{dt}\frac{\partial \mathcal{L}}{\partial \dot{\theta}}-\frac{\partial \mathcal{L}}{\partial \theta}, \frac{d}{dt}\frac{\partial \mathcal{L}}{\partial \dot{x}_k}-\frac{\partial \mathcal{L}}{\partial x_k} (k=0, 5,\cdots, 2N)$, it follows that the equations of motion in the coordinates $ \theta, x_0, x_5, \cdots, x_{2N}$ are the following:
\begin{equation}\label{equations of motion theta}
2\dot{x}_0 (\sum_{j,k=5}^{2N}q_{jk}\dot{x}_j {x}_k+ \dot{\theta}) + {x}_0 (\sum_{j,k=5}^{2N}q_{kj}\ddot{x}_k {x}_j+ \ddot{\theta})=0,
\end{equation}
\begin{equation}\label{equations of motion x}
\left\{
             \begin{array}{lr}
              [\frac{x_k \sum_{j=5}^{2N}(\ddot{x}_j x_j + \dot{x}^2_j)}{x^2_3} + \frac{3{x}_k (\sum_{j=5}^{2N}{x}_j \dot{x}_j)^2}{x^4_3}+ \ddot{x}_k + \ddot{\theta} \sum_{j=5}^{2N}q_{kj} x_j+ 2 \dot{\theta} \sum_{j=5}^{2N}q_{kj} \dot{x}_j]&\\
             + 2\frac{\dot{x}_0}{x_0}[\frac{x_k \sum_{j=5}^{2N}\dot{x}_j x_j}{x^2_3} + \dot{x}_k + \dot{\theta} \sum_{j=5}^{2N}q_{kj} x_j] - \frac{1}{g_3^{\frac{3}{2}} x_0^3}\frac{\partial U(x)}{\partial x_k}=0,  ~~~~~~~~~~k=5,\cdots, 2N  &\\
             \ddot{x}_0- {x_0} (\dot{x}^2_3+\sum_{j=5}^{2N}\dot{x}^2_j + 2\dot{\theta}\sum_{j,k=5}^{2N}q_{jk}\dot{x}_j {x}_k+ \dot{\theta}^2)+ \frac{ U(x)}{g_3^{\frac{3}{2}} x_0^2}=0. &
             \end{array}
\right.
\end{equation}

It's noteworthy that the degeneracy of $x_1, x_2, x_3, x_4$ according to intrinsic symmetrical characteristic  of the $N$-body problem  has  been reduced  in the coordinates $\theta, x_0, x_5, \cdots, x_{2N}$.

It's also noteworthy  that, by using the coordinates $\theta, x_0, x_5, \cdots, x_{2N}$, the angular momentum $\mathcal{J}$ can be represented as
\begin{eqnarray}\label{angular momentum}
\mathcal{J} = \sum_{j=1}^{N} {m_j {\mathbf{r}}_j \times {\dot{\mathbf{r}}}_j}
= \sum_{j=1}^{N} {m_j {\mathbf{r}}^\bot_j \cdot {\dot{\mathbf{r}}}_j}
=\langle{\mathbf{r}}^\bot,{\dot{\mathbf{r}}}\rangle
= g_3 x_0^2(\dot{\theta}+\sum_{j,k=5}^{2N}q_{jk}\dot{x}_j {x}_k)=\frac{\partial \mathcal{L}}{\partial \dot{\theta}},
\end{eqnarray}
where ${{\mathbf{r}}^\bot_j \cdot {\dot{\mathbf{r}}}_j}$ denotes the Euclidean scalar product of ${\mathbf{r}}^\bot_j$ and $ {\dot{\mathbf{r}}}_j$ in ${\mathbb{R}}^2$.

Therefore, the equation (\ref{equations of motion theta}) is none other than angular momentum conservation $\dot{\mathcal{J}}=0$.
As a matter of fact, one can further reduce the equations of motion (\ref{equations of motion theta}) and (\ref{equations of motion x}) by the relation (\ref{angular momentum}):
\begin{equation}\label{equations of motion theta1}
\dot{\theta}=\frac{ \mathcal{J}}{g_3 x_0^2}-\sum_{j,k=5}^{2N}q_{jk}\dot{x}_j {x}_k,
\end{equation}
\begin{equation}\label{equations of motion x1}
\left\{
             \begin{array}{lr}
              \frac{x_k \sum_{j=5}^{2N}(\ddot{x}_j x_j + \dot{x}^2_j)}{x^2_3} + \frac{3{x}_k (\sum_{j=5}^{2N}{x}_j \dot{x}_j)^2}{x^4_3}+ \ddot{x}_k - (\frac{2 \mathcal{J}\dot{x}_0}{g_3 x_0^3}+\sum_{j,k=5}^{2N}q_{jk}\ddot{x}_j {x}_k) \sum_{j=5}^{2N}q_{kj} x_j&\\
              + 2 (\frac{ \mathcal{J}}{g_3 x_0^2}-\sum_{j,k=5}^{2N}q_{jk}\dot{x}_j {x}_k) \sum_{j=5}^{2N}q_{kj} \dot{x}_j+ 2\frac{\dot{x}_0}{x_0}[\frac{x_k \sum_{j=5}^{2N}\dot{x}_j x_j}{x^2_3}&\\
              + \dot{x}_k + \frac{ \mathcal{J}}{g_3 x_0^2}-\sum_{j,k=5}^{2N}q_{jk}\dot{x}_j {x}_k\sum_{j=5}^{2N}q_{kj} x_j] - \frac{g_3^{-\frac{3}{2}}}{x_0^3}\frac{\partial U(x)}{\partial x_k}=0,  ~~~~~~~~~~k=5,\cdots, 2N  &\\
             \ddot{x}_0- {x_0} [\dot{x}^2_3+\sum_{j=5}^{2N}\dot{x}^2_j + (\frac{ \mathcal{J}}{g_3 x_0^2})^2-(\sum_{j,k=5}^{2N}q_{jk}\dot{x}_j {x}_k)^2]+ \frac{g_3^{-\frac{3}{2}} U(x)}{x_0^2}=0. &
             \end{array}
\right.
\end{equation}

Note that the system (\ref{equations of motion x1}) of $x=(x_0, x_5, \cdots, x_{2N})$ are closed (and  autonomous),  and once $x$ are solved by the
system (\ref{equations of motion x1}), the variables $\theta$ can also be solved by quadrature. As a result, it's easy to see that a  periodic orbit  of the system (\ref{equations of motion x1}) yields a (relative) periodic orbit of the $N$-body problem. Thus our task is now finding periodic orbits  of the system (\ref{equations of motion x1}). It's necessary to transform the system (\ref{equations of motion x1}) to Hamiltonian form.

\section{Hamiltonian Near Relative Equilibria}
\indent\par

In this section,  based upon the coordinates $ \theta, x_0, x_5, \cdots, x_{2N}$, we will obtain the Hamiltonian of the planar $N$-body problem near relative equilibrium solutions which does not contain the degeneracy according to intrinsic symmetrical characteristic of the $N$-body problem.

Since  the variable $\theta$ is an ignorable coordinate, we can reduce the Lagrangian $\mathcal{L}$ to just a function of the variables $x=(x_5, \cdots, x_{2N})$:
\begin{eqnarray}\label{Lagrangianx}
L(x_0, x)& =&\mathcal{L}-\frac{\partial \mathcal{L}}{\partial \dot{\theta}}\dot{\theta}\nonumber\\
& =& \frac{g_3 \dot{x}_0^2}{2 }+ \frac{g_3 x_0^2}{2} [ \dot{x}^2_3+\sum_{k=5}^{2N}\dot{x}^2_k - (\frac{ \mathcal{J}}{g_3 x_0^2}-\sum_{j,k=5}^{2N}q_{jk}\dot{x}_j {x}_k)^2] + \frac{ U(x)}{\sqrt{g_3} x_0}.
\end{eqnarray}
It's easy to see that the system (\ref{equations of motion x1}) is the Euler-Lagrange equations of the the action functional corresponding to the Lagrangian $L$  (\ref{Lagrangianx}).

It follows from the Legendre Transform that
the corresponding Hamiltonian is
\begin{eqnarray}\label{Hamiltonian}
H(x_0,x,y_0,y) & =& y_0\dot{x}_0+ \sum_{k=5}^{2N}y_k \dot{x}_k -L(x)\nonumber\\
& =&\frac{ y_0^2}{2 g_3}+ \frac{g_3 x_0^2}{2} [\dot{x}^2_3+\sum_{k=5}^{2N}\dot{x}^2_k + \frac{g_3^2 \mathcal{J}^2}{x_0^4}-(\sum_{j,k=5}^{2N}q_{jk}\dot{x}_j {x}_k)^2] - \frac{ U(x)}{\sqrt{g_3} x_0},
\end{eqnarray}
where
\begin{equation}
\begin{array}{lr}
y_0=\frac{\partial {L}}{\partial \dot{x}_0}= g_3 \dot{x}_0,\\
  y_k=\frac{\partial {L}}{\partial \dot{x}_k}= g_3x^2_0 [\frac{x_k \sum_{j=5}^{2N}\dot{x}_j x_j}{x^2_3} + \dot{x}_k +\sum_{j=5}^{2N}q_{kj} {x}_j(\frac{ \mathcal{J}}{g_3x_0^2}-\sum_{i,j=5}^{2N}q_{ij}\dot{x}_i {x}_j)].
\end{array}
\nonumber
\end{equation}

A straight forward computation shows that:
\begin{displaymath}
\begin{array}{c}
  \dot{x}^2_3 = \frac{ (\sum_{k=5}^{2N}\dot{x}_k x_k)^2}{x^2_3} =(\sum_{k=5}^{2N}\frac{ x_k y_k}{g_3 x_0^2} )^2 + O(\parallel (x,y)\parallel^6),\\
  \dot{x}_k= \frac{y_k-\mathcal{J}\sum_{j=5}^{2N}q_{kj} {x}_j}{g_3x_0^2} -\frac{x_k \sum_{j=5}^{2N}\dot{x}_j x_j}{x^2_3} +\sum_{j=5}^{2N}q_{kj} {x}_j\sum_{i,j=5}^{2N}q_{ij}\dot{x}_i {x}_j\\
  =\frac{y_k-\mathcal{J}\sum_{j=5}^{2N}q_{kj} {x}_j}{g_3 x_0^2}-x_k \frac{\sum_{j=5}^{2N} x_j y_j}{ g_3 x_0^2}+\sum_{j=5}^{2N}q_{kj} {x}_j(\frac{\sum_{i,j=5}^{2N}q_{ij}{x}_j y_i-\mathcal{J}\sum_{j=5}^{2N}{x}^2_j}{g_3 x_0^2})+O(\parallel (x,y)\parallel^5)
\end{array}
\end{displaymath}

As a result, we have
\begin{eqnarray}\label{Hamiltonian1}
 H(x_0,x,y_0,y)
& =&\frac{ y_0^2}{2g_3}+ \frac{1}{2g_3 x_0^2} [\mathcal{J}^2+\sum_{k=5}^{2N}(y_k-\mathcal{J}\sum_{j=5}^{2N}q_{kj} {x}_j)^2 +(\sum_{j,k=5}^{2N}q_{kj}{x}_j y_k-\mathcal{J}\sum_{j=5}^{2N}{x}^2_j)^2\nonumber\\
&-& (\sum_{k=5}^{2N}x_k y_k)^2+ O(\parallel (x,y)\parallel^6)] - \frac{ U(x)}{\sqrt{g_3} x_0}.
\end{eqnarray}

We will consider  the existence of periodic
orbits in a neighbourhood of a relative equilibrium $ A(\omega t)\mathcal{E}_3$ of the Newtonian $N$-body problem.  Then a straight forward computation shows that the angular momentum $\mathcal{J}$ of the relative equilibrium $A(\omega t)\mathcal{E}_3$ is just $g_3 \omega$ and $\lambda =  \omega^2$; so, $\mathcal{J}= g_3^{\frac{1}{4}}\sqrt{\lambda^*}$, and the Hamiltonian becomes
\begin{eqnarray}
H(x_0,x,y_0,y)
& =&\frac{ y_0^2}{2g_3}+ \frac{1}{2g_3 x_0^2} [\sqrt{g_3}\lambda^* +\sum_{k=5}^{2N}(y_k-g_3^{\frac{1}{4}}\sqrt{\lambda^*}\sum_{j=5}^{2N}q_{kj} {x}_j)^2\nonumber\\
& +& (\sum_{j,k=5}^{2N}q_{kj}{x}_j y_k-g_3^{\frac{1}{4}}\sqrt{\lambda^*}\sum_{j=5}^{2N}{x}^2_j)^2
- (\sum_{k=5}^{2N}x_k y_k)^2+ O(\parallel (x,y)\parallel^6)] - \frac{ U(x)}{\sqrt{g_3} x_0}.\nonumber
\end{eqnarray}

Without
loss of generality, one can suppose that $\|\mathcal{E}_3\|=1$ or $g_3=1$. As a matter of fact, by the change of variables $t \rightarrow  \frac{t}{\sqrt{g_3}}$, here $t \rightarrow \frac{t}{\sqrt{g_3}}$ means replace $t$ by $\frac{t}{\sqrt{g_3}}$ everywhere, the Hamiltonian becomes
\begin{eqnarray}
 H(x_0,x,y_0,y)
& =&\frac{ y_0^2}{2\sqrt{g_3}}+ \frac{1}{2\sqrt{g_3} x_0^2} [\sqrt{g_3}\lambda^* +\sum_{k=5}^{2N}(y_k-g_3^{\frac{1}{4}}\sqrt{\lambda^*}\sum_{j=5}^{2N}q_{kj} {x}_j)^2\nonumber\\
& +& (\sum_{j,k=5}^{2N}q_{kj}{x}_j y_k-g_3^{\frac{1}{4}}\sqrt{\lambda^*}\sum_{j=5}^{2N}{x}^2_j)^2
- (\sum_{k=5}^{2N}x_k y_k)^2+ O(\parallel (x,y)\parallel^6)] - \frac{ U(x)}{ x_0}.\nonumber
\end{eqnarray}
Furthermore, by the change of variables $y_0 \rightarrow \frac{y_0}{\sqrt[4]{g_3}},y \rightarrow \frac{y}{\sqrt[4]{g_3}}$, the Hamiltonian becomes
\begin{eqnarray}\label{Hamiltonian2}
H(x_0,x,y_0,y)
 &=&\frac{ y_0^2}{2}+ \frac{1}{2 x_0^2} [\lambda^* +\sum_{k=5}^{2N}(y_k-\sqrt{\lambda^*}\sum_{j=5}^{2N}q_{kj} {x}_j)^2
+ (\sum_{j,k=5}^{2N}q_{kj}{x}_j y_k-\sqrt{\lambda^*}\sum_{j=5}^{2N}{x}^2_j)^2\nonumber\\
&-& (\sum_{k=5}^{2N}x_k y_k)^2+ O(\parallel (x,y)\parallel^6)] - \frac{ U(x)}{ x_0}.
\end{eqnarray}
It's noteworthy that the above changes of variables are not symplectic but symplectic with multiplier. 

By the coordinates $ (r, \theta, x_5, \cdots, x_{2N})$, the relative equilibrium $ A(\omega t)\mathcal{E}_3$ is just a solution of equations of motion (\ref{equations of motion theta}) and (\ref{equations of motion x}) such that
\begin{displaymath}
r=1, \theta=\omega t, x_5=0,\cdots, x_{2N}=0.
\end{displaymath}
By the coordinates $ (x_0,x,y_0,y)$, the relative equilibrium $ A(\omega t)\mathcal{E}_3$ is just an equilibrium solution of the Hamiltonian system with  Hamiltonian (\ref{Hamiltonian2}) such that
\begin{displaymath}
x_0=1,  x=0, y_0=0, y=0.
\end{displaymath}
As a matter of notational convenience, the above equilibrium point of the Hamiltonian (\ref{Hamiltonian2}) will be translated to the origin by substituting $x_0$ for $x_0+1$. Then the Hamiltonian (\ref{Hamiltonian2}) becomes:
\begin{eqnarray}\label{Hamiltonian3}
&&H(x_0,x,y_0,y)=\frac{-\omega^2}{2} +\frac{y_0^2}{2}+ \frac{1}{2} [ \omega^2 x_0^2+\sum_{k=5}^{2N}y_k^2-2\omega\sum_{j,k=5}^{2N}q_{kj} {x}_j y_k+\sum_{k=5}^{2N}( \omega^2-\lambda^*_k ){x}^2_k] \nonumber\\
&-&[x_0 (\omega^2 x_0^2+\sum_{k=5}^{2N}y_k^2-2\omega\sum_{j,k=5}^{2N}q_{kj} {x}_j y_k+\sum_{k=5}^{2N}( \omega^2-\frac{\lambda^*_k}{2} ){x}^2_k)+\frac{1}{6}\sum_{i,j,k = 5}^{2N} a_{ijk} x_i x_j x_k ]\nonumber\\
&+& \frac{1}{2} [(\sum_{j,k=5}^{2N}q_{kj}{x}_j y_k-\omega\sum_{k=5}^{2N}{x}^2_k)^2- (\sum_{k=5}^{2N}x_k y_k)^2]+ \frac{x_0}{6}\sum_{i,j,k = 5}^{2N} a_{ijk} x_i x_j x_k\nonumber\\
&+& \frac{3x_0^2}{2} [\omega^2 x_0^2+\sum_{k=5}^{2N}y_k^2-2\omega\sum_{j,k=5}^{2N}q_{kj} {x}_j y_k+\sum_{k=5}^{2N}( \omega^2-\frac{\lambda^*_k}{3} ){x}^2_k]\nonumber\\
&-&   \frac{3}{4} (\sum_{k = 5}^{2N}  x^2_k)\sum_{k = 5}^{2N} (\lambda^*_k - \frac{\omega^2}{2})  x^2_k-\frac{1}{24}\sum_{h,i,j,k = 5}^{2N} a_{hijk} x_h x_i x_j x_k+\cdots
\end{eqnarray}
where $\cdots$ denotes higher order terms and $\omega=\sqrt{\lambda^*}$. Without loss of generality, we will sometimes omit  the constant term $\frac{-\omega^2}{2}$ of the Hamiltonian (\ref{Hamiltonian3}) in the following content.

The problem is now reduced to find periodic
orbits of the Hamiltonian system with  Hamiltonian (\ref{Hamiltonian3}) in a neighbourhood of the origin.

First, it is easy to see that the manifold $\{(x_0,x,y_0,y):x=y=0\}$ is an invariant manifold and on this invariant manifold the problem is reduced to the following problem of  single degree of freedom:
\begin{eqnarray}
 H(x_0,0,y_0,0)&=&\frac{y_0^2}{2}+ \frac{\omega^2 }{2(1+x_0)^{2}} - \frac{\omega^2}{1+x_0}\nonumber\\
 &=&\frac{-\omega^2}{2} +\frac{y_0^2}{2}+ \frac{\omega^2 x_0^2}{2}  -\omega^2 x_0^3
+ \frac{3\omega^2 x_0^4}{2}+\cdots .\nonumber
\end{eqnarray}
Then it's easy to see that the invariant manifold $\{(x_0,x,y_0,y):x=y=0\}$  is constituted by periodic
orbits.
However, this fact is trivial. As a matter of fact, when $x=y=0$, i.e., on the invariant manifold, the primary $N$-body problem is reduced to the two-body problem, the periodic
orbits are just Keplerian elliptic orbits generated by the central configuration $\mathcal{E}_3$. These periodic
orbits on the invariant manifold $\{(x_0,x,y_0,y):x=y=0\}$ are called the trivial family of periodic
orbits near relative equilibria.

Next, let's explore other nontrivial periodic orbits near relative equilibria.

\section{Three-body Problem}
\indent\par
In this section, we  discuss the problem of periodic orbits near five relative equilibrium solutions of the planar three-body problem in detail. We will show that there are abundant periodic
orbits near the five relative equilibria. More specifically, we first rediscover the well known Lyapunov's orbits or Weinstein's orbits of the planar three-body problem. Then we further find that there are abundant Conley-Zender's periodic
orbits for the planar three-body problem.

\subsection{Lagrange relative equilibrium}
\indent\par
First, let's consider the problem  near a Lagrange relative equilibrium.

As in \cite{yu2019stability}, set $\beta=m_1 m_2+ m_3 m_2+ m_1 m_3$. Without loss of generality, suppose
\begin{displaymath}
\begin{array}{c}
  m_1+m_2+m_3=1, \\
  \mathbf{r}= (-\frac{\sqrt{3} m_3}{2 \sqrt{\beta }},\frac{2 m_2+m_3}{2 \sqrt{\beta }},-\frac{\sqrt{3} m_3}{2 \sqrt{\beta }},-\frac{2 m_1+m_3}{2 \sqrt{\beta }},\frac{\sqrt{3} \left(m_1+m_2\right)}{2 \sqrt{\beta }},-\frac{m_1-m_2}{2 \sqrt{\beta }})^\top.\\
\end{array}
\end{displaymath}
Here $ \mathbf{r}$ is a Lagrange central configuration such that $\|\mathbf{r}\|=1$.

Then, as which has been obtained in \cite{yu2019stability}, we have
\begin{equation}\label{eigenvectors}
\begin{array}{c}
  \mathcal{E}_3=\mathbf{r},\\
  \mathcal{E}_4 = \mathbf{r}^\bot,\\
  \omega^2=\lambda=\lambda^*=\beta ^{3/2},
\end{array}\nonumber
\end{equation}
\begin{equation}\label{eigenvalues}
  \lambda_5=\frac{3}{2}  \left(1-\alpha\right) \beta ^{3/2}, ~~~~
  \lambda_6=\frac{3}{2} \left(\alpha+1\right) \beta ^{3/2}, ~~~~\nonumber
\end{equation}
\begin{equation}\label{eigenvectors}
\begin{array}{c}
  \mathcal{E}_5=\sqrt{\frac{3 m_1 m_2}{4 \beta  m_3 \left(2\alpha^2 +\alpha-3 \alpha m_2\right)}} \left(\frac{m_1-m_3 }{\frac{ m_1}{m_3}},\frac{ 3 m_2-2 \alpha-1}{\frac{\sqrt{3} m_1}{m_3}},\frac{ m_2-\alpha-m_1}{\frac{  m_2}{m_3}},\frac{ \alpha+3 m_3-1}{\frac{\sqrt{3}  m_2}{m_3}},\alpha-m_2+m_3,\frac{\alpha+3 m_1-1}{\sqrt{3}}\right)^\top, \\
  \mathcal{E}_6 = \mathcal{E}_5^\bot,
\end{array}\nonumber
\end{equation}
where
\begin{equation}
\alpha=\sqrt{1-3 \beta }.\nonumber
\end{equation}

Let $\Omega$ be the space of masses of the planar three-body problem, then $\Omega$ could  be represent as
\begin{equation}
\Omega=\{(\beta,m_1): \beta\in(0,\frac{1}{3}], m_1\in[\frac{1}{3},1), \beta- m_1(1-m_1)>0, 4\beta\leq 1+2 m_1-3 m_1^2\}.\nonumber
\end{equation}

As a result,  the
quadratic part $H_2$ of the Hamiltonian (\ref{Hamiltonian3}) is
\begin{eqnarray}\label{HamiltonianH2}
H_2=\frac{y_0^2}{2}+ \frac{1}{2} [ \omega^2 x_0^2+y_5^2+y_6^2-2\omega ({x}_5 y_6-{x}_6 y_5)+( \omega^2-\lambda_5 ){x}^2_5+( \omega^2-\lambda_6 ){x}^2_6], \nonumber
\end{eqnarray}
and the $6$ eigenvalues are
\begin{equation}\label{eigenvalues1}
  \pm \omega_0 \mathbf{i}, ~~~~\pm \omega_1 \mathbf{i}, ~~~~\pm \omega_2 \mathbf{i}, ~~~~ ~~~~\nonumber
\end{equation}
where  $\textbf{i}$ denotes the imaginary unit, and
\begin{equation}
\omega_0=\omega= \beta ^{3/4},~~~~\omega _1=\sqrt{\frac{{1-\gamma} }{{2}}} \omega_0,~~~~\omega _2=\sqrt{\frac{{1+\gamma} }{{2}}} \omega_0,\nonumber
\end{equation}
\begin{equation}
\gamma= \sqrt{1-27 \beta }.\nonumber
\end{equation}

\textbf{A. $\beta > \frac{1}{27}$.}
When
 \begin{equation}
\beta > \frac{1}{27},\nonumber
\end{equation}
or more precisely,
 \begin{equation}
m_1 m_2+ m_3 m_2+ m_1 m_3 >\frac{(m_1+m_2+m_3)^2}{27}.\nonumber
\end{equation}
The  eigenvalues $\pm \omega_0 \mathbf{i}$ are only purely imaginary, periodic
orbits in a neighbourhood of the origin are all in the central manifold $\mathcal{W}^{c}_{loc}: x=y=0$, thus there are only the trivial family of periodic
orbits in a neighbourhood of the the origin and these periodic
orbits constitute the central manifolds $\mathcal{W}^{c}_{loc}: x=y=0$.

\textbf{B. $\beta = \frac{1}{27}$.}

When
 \begin{equation}
\beta = \frac{1}{27},\nonumber
\end{equation}
i.e.,
\begin{equation}
m_1 m_2+ m_3 m_2+ m_1 m_3 =\frac{(m_1+m_2+m_3)^2}{27}.\nonumber
\end{equation}
Except the trivial family of periodic
orbits near the origin, we don't know if there are other  periodic
orbits near the origin.

\textbf{C. $\beta < \frac{1}{27}$.} This is the main case we will discussed.

When
 \begin{equation}
\beta < \frac{1}{27},\nonumber
\end{equation}
i.e.,
 \begin{equation}
m_1 m_2+ m_3 m_2+ m_1 m_3 < \frac{(m_1+m_2+m_3)^2}{27}.\nonumber
\end{equation}

Without loss of generality, suppose $m_1\geq m_2\geq m_3$. Then it's easy to see that
 \begin{equation}
m_1>\frac{1}{18} \left(\sqrt{69}+9\right)>0.961478, m_2+m_3<0.038521.\nonumber
\end{equation}

The  $6$ eigenvalues are
\begin{equation}
  \pm \omega_0 \mathbf{i}, ~~~~\pm \omega_1 \mathbf{i}, ~~~~\pm \omega_2 \mathbf{i}, ~~~~ ~~~~\nonumber
\end{equation}
are all purely imaginary.

First, it's easy to see that $\frac{\omega_1}{\omega_0}, \frac{\omega_2}{\omega_0} \in (0,1)$.
By Lyapunov Theorem \ref{Lyapunov}, to the eigenvalues $\pm \omega_0 \mathbf{i}$ there corresponds a
one parameter family of periodic orbits  that lie near the origin and have the approximate period $\frac{2\pi}{\omega_0}$. However, these periodic orbits are just the trivial family of periodic
orbits.

It's also easy to see that $\frac{\omega_1}{\omega_2} \in (0,1)$, $\frac{\omega_0}{\omega_2}\in (1,\sqrt{2})$.
By Lyapunov Theorem \ref{Lyapunov}, to the eigenvalues $\pm \omega_2 \mathbf{i}$ there corresponds a
one parameter family of periodic orbits  that lie near the origin and have the approximate period $\frac{2\pi}{\omega_2}$.

For the eigenvalues $\pm \omega_1 \mathbf{i}$, it's easy to see that $\frac{\omega_1}{\omega_1}, \frac{\omega_0}{\omega_1} \notin \mathbb{N}$ provided that $\beta$ dose not belong to the set
\begin{equation}
 \{\beta: \beta=\frac{1-(1-\frac{2}{n^2})^2}{27} ~or~ \frac{1-(1-\frac{2}{n^2+1})^2}{27}, ~integer~  n\geq 2\}.\nonumber
\end{equation}
Therefore,  by Lyapunov Theorem \ref{Lyapunov}, to the eigenvalues $\pm \omega_1 \mathbf{i}$ there corresponds a
one parameter family of periodic orbits  that lie near the origin and have the approximate period $\frac{2\pi}{\omega_1}$ provided that $\beta$ dose not belong to the above set.

The above results by Lyapunov Theorem have been obtained by Siegel \cite{siegel1971lectures}.

Thus generically there are three one parameter family of periodic orbits  that lie near the origin, in other words, there are three two-dimensional invariant manifold constituted by periodic
orbits. However, we can further prove that periodic
orbits are unexpectedly abundant: generically  the relative measure of the closure of the set of periodic orbits near the origin is close to 1.

To prove this fact, it's necessary to get the  Birkhoff normal form of the Hamiltonian (\ref{Hamiltonian3}).

We confine ourselves to the following space of masses of the planar three-body problem
\begin{equation}\label{space of massesps}
\Omega_{ps}=\{(\beta,m_1)\in \Omega: \beta\in(0,\frac{1}{27})\setminus \{\frac{1}{75},\frac{32}{2187},\frac{16}{675},\frac{1}{36},\frac{64}{1875}\}, m_1\in(\frac{\sqrt{69}+9}{18} ,1)\}
\end{equation}
which has removed  masses corresponding to special resonant cases.
Then, as in \cite{yu2019stability}, we know that
the  Birkhoff normal form of the Hamiltonian (\ref{Hamiltonian3}) is
\begin{equation}\label{Birkhoff normal form}
\begin{aligned}
&\mathcal{H}(\varrho,\varphi)
= \omega _0 \varrho_0 - \omega _1 \varrho_1 + \omega _2 \varrho _2  +\frac{1}{2}[\omega_{00}\varrho_0 ^2+\omega_{11}\varrho_1 ^2+\omega_{22}\varrho_2 ^2\\
&+2\omega_{01}\varrho_0 \varrho_1 +2\omega_{02}\varrho_0  \varrho_2 +2\omega_{12}\varrho_1 \varrho_2 ] +\cdots,
\end{aligned}
\end{equation}
where $\varrho_j$ ($j=0,1,2$) are action variables, and
\begin{equation}
\omega _{00}=-3,\nonumber
\end{equation}
\begin{equation}
\omega _{01}=-\frac{\sqrt{\gamma +1} \left(21 \gamma ^3-40 \gamma ^2+15 \gamma +4\right)}{12 \sqrt{6} \sqrt{\beta } \gamma  (2 \gamma -1)}, \nonumber
\end{equation}
\begin{equation}
\omega _{02}=-\frac{\sqrt{\gamma +1} \left(21 \gamma ^2+19 \gamma -4\right)}{4 \sqrt{2} \gamma  (2 \gamma +1)},\nonumber
\end{equation}
\begin{equation}
\begin{aligned}
&\omega _{12}=\frac{\sqrt{3 \beta } }{4 (18225 \beta ^2-1107 \beta +16) m_1 m_2 m_3} [(360855 \beta ^2-32265 \beta +624) m_1^3+\\
&(-360855 \beta ^2+32265 \beta -624) m_1^2+3 \beta  (120285 \beta ^2-10755 \beta +208) m_1-4 \beta ^2 (432 \beta +43)],
\end{aligned}\nonumber
\end{equation}
\begin{equation}\label{omega11}
\omega _{11}=\frac{(\gamma -1) \left(1211 \gamma ^4-1336 \gamma ^3+279 \gamma ^2+158 \gamma -76\right)}{72 \gamma ^2 \left(10 \gamma ^2-11 \gamma +3\right)}-\frac{3 \beta ^3 \left(31 \gamma ^2+286 \gamma -236\right)}{8 (\gamma -1) \gamma ^2 (5 \gamma -3) m_1 m_2 m_3},\nonumber
\end{equation}
\begin{equation}\label{omega22}
\omega _{22}=-\frac{(\gamma +1) \left(1211 \gamma ^4+1336 \gamma ^3+279 \gamma ^2-158 \gamma -76\right)}{72 \gamma ^2 \left(10 \gamma ^2+11 \gamma +3\right)}-\frac{3 \beta ^3 \left(31 \gamma ^2-286 \gamma -236\right)}{8 \gamma ^2 (\gamma +1) (5 \gamma +3) m_1 m_2 m_3}.\nonumber
\end{equation}

\begin{equation}\label{degerate}
det\left(
\begin{array}{ccc}
 \omega _{00} & \omega _{01} & \omega _{02} \\
 \omega _{01} & \omega _{11} & \omega _{12} \\
 \omega _{02} & \omega _{12} & \omega _{22} \\
\end{array}
\right)=\frac{-27 \beta  }{128 (16-675 \beta )^2 (1-36 \beta )^2 \gamma ^4 m_1^2 m_2^2 m_3^2} f_{deg},
\end{equation}
where
\begin{equation}\label{degeratef}
\begin{aligned}
& f_{deg} = \frac{2(1-36 \beta )^2}{3}  (52542675 \beta ^3+178185258 \beta ^2-9896841 \beta -47632) \beta ^4\\
&-11 (397050199920 \beta ^5-40790893923 \beta ^4+4055047758 \beta ^3-243771759 \beta ^2+6417616 \beta \\
&-59392) \beta ^3 m_1+(5465578392450 \beta ^6+19309935720393 \beta ^5-3995019640449 \beta ^4\\
&+327340481715 \beta ^3-13039336341 \beta ^2+250520816 \beta -1857536) \beta ^2 m_1^2+(2408448\\
&-15298708984020 \beta ^6-29436067209393 \beta ^5+7048034089254 \beta ^4-562788423405 \beta ^3\\
&+20645100208 \beta ^2-359200768 \beta ) \beta  m_1^3+3 [(1821859464150 \beta ^6+4980794507091 \beta ^5\\
&-1182106602432 \beta ^4+94244985459 \beta ^3-3452615664 \beta ^2+59975680 \beta -401408) \\
&m_1^4 (2 \beta +m_1^2-2 m_1+1)].\nonumber
\end{aligned}
\end{equation}
Therefore, the Hamiltonian (\ref{Hamiltonian3}) is nondegenerate if and only if
\begin{equation}
f_{deg}\neq 0.\nonumber
\end{equation}
Thus the set $V_{f_{deg}}$ of points $(\beta,m_1)$ such that the Hamiltonian (\ref{Hamiltonian3}) is degenerate is  a real algebraic variety.
Moreover, the real algebraic variety  $V_{f_{deg}}$ is union of a finite number of zero-dimensional points and one-dimensional ``curves".

To make the direct-viewing understanding of the real algebraic variety  $V_{f_{deg}}$ and the  spaces of masses $\Omega,\Omega_{ps}$,  please see  \cite{yu2019stability}.

As a result, it follows from  Theorem \ref{Conley and Zender}  that
\begin{theorem}\label{Conley and Zenderapply1}
For every choice of positive masses of the planar three-body problem satisfying  $(\beta,m_1)\in \Omega_{ps}\setminus V_{f_{deg}}$, there are abundant periodic orbits near every Lagrange relative equilibrium, and the relative measure of the closure of the set of periodic orbits in  the polydisc $\|\varrho\|<\epsilon$ is
at least $1-O(\epsilon^{\frac{1}{4}})$.
\end{theorem}

Furthermore, it follows from  Remark \ref{KAM1} and the following  Theorem \ref{Diophantine frequency} that the relative measure in the above theorem is
at least $1-O(\exp (-\tilde{c}\epsilon^{\frac{-1}{\upsilon+1}}))$ generically, here $\upsilon>6$ and $\tilde{c}$ are positive constants.

\begin{theorem}\label{Diophantine frequency}\emph{(\cite{yu2019stability})}
The set $\Gamma_r$ of $\beta \in (0,\frac{1}{27})$ corresponding to resonant frequency vectors $\varpi=(\omega _0,-\omega _1,\omega _2)$  is  countable and dense. The set $\Gamma_d$ of $\beta \in (0,\frac{1}{27})$ corresponding to $(c,\upsilon)$-Diophantine frequency vectors $\varpi=(\omega _0,-\omega _1,\omega _2)$ is a  set of full measure for $\upsilon>6$.
\end{theorem}

\subsection{Euler relative equilibrium}
\indent\par
Let's now consider the problem  near a Euler relative equilibrium.

As $\mathbf{r}$ is an Euler (collinear) central configuration such that $\|\mathbf{r}\|=1$ and whose center of masses is at the origin,
without loss of generality, suppose
\begin{displaymath}
\begin{array}{c}
  m_1+m_2+m_3=1, \\
  \mathbf{r}= (\xi_1,0,\xi_3,0,\xi_5,0)^\top,\\
\end{array}
\end{displaymath}
where
\begin{displaymath}
\begin{array}{c}
  \xi_1 = -\kappa  \left(m_2 (\sigma+\frac{1}{2}) +m_3\right), \\
   \xi_3 = \kappa  \left(m_1 (\sigma+\frac{1}{2})+m_3 (\sigma-\frac{1}{2})\right),\\
    \xi_5 =\kappa \left(m_2 (\frac{1}{2}-\sigma)+m_1\right),
\end{array}
\end{displaymath}
and the parameters $\sigma\in (-\frac{1}{2},\frac{1}{2}), \kappa >0$
satisfy the equations of  central configurations:
\begin{equation}\label{centralconfigurationeuler}
\left\{
             \begin{array}{lr}
             \frac{m_2}{(\sigma+\frac{1}{2})^2}+m_3=\lambda \kappa^3  \left(m_2 (\sigma+\frac{1}{2}) +m_3\right)&\\
             -\frac{m_1}{(\sigma+\frac{1}{2})^2}+\frac{m_3}{(\frac{1}{2}-\sigma)^2}=-\lambda \kappa^3 \left(m_1 (\sigma+\frac{1}{2})+m_3 (\sigma-\frac{1}{2})\right)  &\\
            m_1+\frac{m_2}{(\frac{1}{2}-\sigma)^2}=\lambda \kappa^3 \left(m_2 (\frac{1}{2}-\sigma)+m_1\right) &
             \end{array}
\right.
\end{equation}
and
\begin{displaymath}
I(\mathbf{r})=\kappa^2(m_1 m_2 (\sigma+\frac{1}{2})^2+m_2 m_3 (\sigma-\frac{1}{2})^2+m_1 m_3)=1.
\end{displaymath}

Then
\begin{equation}\label{eigenvectors}
\begin{array}{c}
  \mathcal{E}_3=\mathbf{r},\\
  \mathcal{E}_4 = \mathbf{r}^\bot,\\
  \lambda=\frac{m_1 m_2}{\kappa  (\sigma+\frac{1}{2})}+\frac{m_1 m_3}{\kappa }+\frac{m_3 m_2}{\kappa  (\frac{1}{2}-\sigma)},
\end{array}
\end{equation}
and the matrix $\lambda \mathbb{I} + \mathfrak{M}^{-1} \mathfrak{B}$ is
\begin{equation}
    \begin{pmatrix}
        \begin{smallmatrix}
           \lambda +\frac{ m_3 +\frac{m_2}{(\sigma+\frac{1}{2})^3}}{\kappa ^3 /2} & 0 & -\frac{2 m_2}{\kappa ^3 (\sigma+\frac{1}{2})^3} & 0 & -\frac{2 m_3}{\kappa ^3} & 0 \\
 0 & \lambda -\frac{\frac{m_2}{(\sigma+\frac{1}{2})^3}+m_3}{\kappa ^3} & 0 & \frac{m_2}{\kappa ^3 (\sigma+\frac{1}{2})^3} & 0 & \frac{m_3}{\kappa ^3} \\
 -\frac{2 m_1}{\kappa ^3 (\sigma+\frac{1}{2})^3} & 0 & \lambda +\frac{ \frac{m_1}{(\sigma+\frac{1}{2})^3}+\frac{m_3}{(\frac{1}{2}-\sigma)^3}}{\kappa ^3/2} & 0 & -\frac{2 m_3}{\kappa ^3 (\frac{1}{2}-\sigma)^3} & 0 \\
 0 & \frac{m_1}{\kappa ^3 (\sigma+\frac{1}{2})^3} & 0 & \lambda -\frac{\frac{m_3}{(\frac{1}{2}-\sigma)^3}+\frac{m_1}{(\sigma+\frac{1}{2})^3}}{\kappa ^3} & 0 & \frac{m_3}{\kappa ^3 (\frac{1}{2}-\sigma)^3} \\
 -\frac{2 m_1}{\kappa ^3} & 0 & -\frac{2 m_2}{\kappa ^3 (\frac{1}{2}-\sigma)^3} & 0 & \lambda +\frac{m_1+\frac{m_2}{(\frac{1}{2}-\sigma)^3}}{\kappa ^3/2} & 0 \\
 0 & \frac{m_1}{\kappa ^3} & 0 & \frac{m_2}{\kappa ^3 (\frac{1}{2}-\sigma)^3} & 0 & \lambda -\frac{\frac{m_2}{(\frac{1}{2}-\sigma)^3}+m_1}{\kappa ^3} \\
        \end{smallmatrix}
    \end{pmatrix}\nonumber
\end{equation}
As a result of (\ref{Hessian11}), $\lambda_5,\lambda_6$ and $\{ \mathcal{E}_5, \mathcal{E}_6\}$ can be obtained by calculating eigenvalues and eigenvectors of the above matrix.

A straight forward computation shows that
the remaining  eigenvalues and eigenvectors  are
\begin{equation}\label{eigenvalueseuler}
\begin{array}{c}
   \lambda_5=-\lambda +\frac{2 \left(m_1+m_3\right)}{\kappa ^3}+\frac{2 \left(m_1+m_2\right)}{\kappa ^3 (\sigma+\frac{1}{2})^3}+\frac{2 \left(m_2+m_3\right)}{\kappa ^3 (\frac{1}{2}-\sigma)^3}=3 \lambda -2\lambda_6, \\
  \lambda_6=\frac{3\lambda -\lambda_{5}}{2}=\frac{-64 \sigma ^4+160 \sigma ^2+28}{\kappa ^5 \left(4 \sigma ^2-1\right)^3},
\end{array}\nonumber
\end{equation}
\begin{equation}\label{eigenvectorseuler}
\begin{array}{c}
  \mathcal{E}_5=\left(\sqrt{\frac{m_2 m_3}{m_1}} \kappa  (\frac{1}{2}-\sigma),0,-\sqrt{\frac{m_1 m_3}{m_2}} \kappa ,0,\sqrt{\frac{m_1 m_2}{m_3}} \kappa  (\sigma+\frac{1}{2}),0\right)^\top, \\
  \mathcal{E}_6 = \mathcal{E}_5^\bot.
\end{array}\nonumber
\end{equation}

For convenience's sake, since the equations (\ref{centralconfigurationeuler}) can not be solved explicitly,   one can use the parameters $\sigma, \kappa$ as independent variables. Then
\begin{displaymath}
\begin{array}{c}
 \lambda= \frac{-64 \sigma ^4+160 \sigma ^2+28}{\kappa ^5 \left(4 \sigma ^2-1\right)^3} -\frac{8 \left(4 \sigma ^2+7\right)}{\kappa ^3 \left(16 \sigma ^4-40 \sigma ^2-7\right)}, \\
  m_1=-\frac{(2 \sigma +1)^2 \left(\kappa ^3 \lambda  (2 \sigma -1)^3+8\right)}{32 \sigma ^4-80 \sigma ^2-14}, \\
  m_2=-\frac{\left(1-4 \sigma ^2\right)^2 \left(\kappa ^3 \lambda -1\right)}{16 \sigma ^4-40 \sigma ^2-7}, \\
  m_3=\frac{(1-2 \sigma )^2 \left(\kappa  \lambda  (2 \sigma +1)^3-8\right)}{32 \sigma ^4-80 \sigma ^2-14},
\end{array}
\end{displaymath}
and the parameters $\sigma, \kappa$ satisfy one of the following conditions:
\begin{itemize}
                   \item $\sigma =0\land \kappa >2$;
  \item $0<\sigma <\frac{1}{2}\land f_{\kappa}(\sigma) <\kappa ^2<g_{\kappa}(\sigma)$;
  \item $0>\sigma >-\frac{1}{2}\land f_{\kappa}(-\sigma) <\kappa ^2<g_{\kappa}(-\sigma)$.
                 \end{itemize}
Where
\begin{displaymath}
\begin{array}{c}
 f_{\kappa}(\sigma)= -\frac{4 \left(-16 \sigma ^4+40 \sigma ^2+7\right)^2}{\left(4 \sigma ^2-1\right)^3 \left(16 \sigma ^4-8 \sigma ^2+49\right)}, \\
 g_\kappa(\sigma)=-\frac{\left(-16 \sigma ^4+40 \sigma ^2+7\right)^2}{4 \sigma  (2 \sigma -1)^3 \left(16 \sigma ^4+32 \sigma ^3+40 \sigma ^2+24 \sigma +21\right)}.
\end{array}
\end{displaymath}

As a result,  the
quadratic part $H_2$ of the Hamiltonian (\ref{Hamiltonian3}) is
\begin{eqnarray}\label{HamiltonianH2}
H_2=\frac{y_0^2}{2}+ \frac{1}{2} [ \omega^2 x_0^2+y_5^2+y_6^2-2\omega ({x}_5 y_6-{x}_6 y_5)+( \omega^2-\lambda_5 ){x}^2_5+( \omega^2-\lambda_6 ){x}^2_6], \nonumber
\end{eqnarray}
and the $6$ eigenvalues are
\begin{equation}\label{eigenvalues1}
  \pm \omega_0 \mathbf{i}, ~~~~\pm \omega_1 \mathbf{i}, ~~~~\pm \omega_2 , ~~~~ ~~~~\nonumber
\end{equation}
where
\begin{equation}
\begin{array}{c}
 \omega_0=\omega= \sqrt{\lambda}, \\
  \omega _1=\frac{\sqrt{\sqrt{9 \lambda_6^2-10 \lambda_6 \omega ^2+\omega ^4}+\lambda_6+\omega ^2}}{\sqrt{2}}, \\
  \omega _2=\frac{\sqrt{\sqrt{9 \lambda_6^2-10 \lambda_6 \omega ^2+\omega ^4}-\lambda_6-\omega ^2}}{\sqrt{2}}.
\end{array}\nonumber
\end{equation}
It's noteworthy that $\omega _1>\omega _0$. Furthermore, it is straightforward to  show that  $\omega _1,\omega _0$ at most satisfies a resonance
relations of order 3 :$\omega _1=2\omega _0$.

So by Lyapunov Theorem \ref{Lyapunov}, to the eigenvalues $\pm \omega_1 \mathbf{i}$ there corresponds a
one parameter family of periodic orbits  that lie near the origin and have the approximate period $\frac{2\pi}{\omega_1}$.
Thus  there are two one parameter family of periodic orbits  that lie near the origin, in other words, there are two two-dimensional invariant manifold constituted by periodic
orbits.

The above results by Lyapunov Theorem have been obtained by Siegel \cite{siegel1971lectures}.

Note that all of  periodic orbits near the origin lie on a $4$-dimensional central manifold of the Hamiltonian system with  Hamiltonian (\ref{Hamiltonian3}).
It is straightforward to  show that the
quadratic part $H_2$ of the reduced Hamiltonian  on the $4$-dimensional central manifold is positive definite (see the following (\ref{hamil})). So by Theorem \ref{weinstein}, one can again prove that there are two one parameter family of periodic orbits  that lie near the origin.

However, we can further prove that periodic
orbits are unexpectedly abundant (see the following Theorem \ref{Conley and Zenderapply2}):  the relative measure of the closure of the set of periodic orbits near the origin in the four-dimensional central manifold is close to 1.

To prove this fact, by Theorem \ref{Conley and Zendergeneral}, it's necessary to get the  Birkhoff normal form of the Hamiltonian (\ref{Hamiltonian3}) on the central manifold.

First, some tedious computation further yields
\begin{equation}
\begin{aligned}
U(x_3\mathcal{E}_3 +x_5\mathcal{E}_5+x_6\mathcal{E}_6)=\lambda+\frac{\lambda _5 x_5^2+ \lambda _6 x_6^2}{2} +
a_{30} x_5^3+a_{12} x_5 x_6^2+ a_{40} x_5^4 +a_{22} x_5^2 x_6^2+a_{04} x_6^4 +\cdots,
\end{aligned}\nonumber
\end{equation}
where
\begin{equation}
\begin{array}{lr}
  a_{30}=\frac{\frac{m_3^2 \xi_5^3}{(\sigma+\frac{1}{2}) ^4}+\frac{m_1^2 \xi_1^3}{(\sigma-\frac{1}{2})^4}-m_2^2 \xi_3^3}{\kappa ^4 \sqrt{m_1 m_2 m_3}}=-\frac{2}{3}  a_{12},
\end{array}
\nonumber
\end{equation}
\begin{equation}
\begin{array}{lr}
  a_{40}= a_4+\frac{3}{8} (5 \lambda -4 \lambda _6), \\
  a_{04}= \frac{3}{8} (a_4-\lambda +2 \lambda _6),\\
   a_{22}=\frac{3}{4} (-4 a_4+2 \lambda -\lambda _6),\\
   a_4=\frac{\frac{m_3^3 \xi_5^4}{(\sigma+\frac{1}{2}) ^5}-\frac{m_1^3 \xi_1^4}{(\sigma-\frac{1}{2})^5}+m_2^3 \xi_3^4}{\kappa ^5 m_1 m_2 m_3}.
\end{array}\nonumber
\end{equation}
Thus
\begin{equation}\label{Hamiltonianeuler1}
\begin{aligned}
&H(r,x,s,y)=\frac{-\omega^2}{2} +\frac{y_0^2}{2}+ \frac{1}{2} [\omega^2 x_0^2+y_5^2+y_6^2+2\omega  ({x}_6 y_5- {x}_5 y_6)+( \omega^2-\lambda_5 ){x}^2_5+( \omega^2-\lambda_6 ){x}^2_6]\\
&-x_0 [\omega^2 x_0^2+y_5^2+y_6^2+2\omega  ({x}_6 y_5- {x}_5 y_6)+( \omega^2-\frac{\lambda_5}{2} ){x}^2_5+( \omega^2-\frac{\lambda_6}{2} ){x}^2_6]-(
a_{30} x_5^3+a_{12} x_5 x_6^2)\\
&+\frac{1}{2} [({x}_6 y_5- {x}_5 y_6+\omega({x}^2_5+{x}^2_6))^2- (x_5 y_5+x_6 y_6)^2]-[  a_{40} x_5^4 +a_{22} x_5^2 x_6^2+a_{04} x_6^4]\\
&+x_0[
a_{30} x_5^3+a_{12} x_5 x_6^2]+\frac{3x_0^{2}}{2} [\omega^2 x_0^2+y_5^2+y_6^2+2\omega  ({x}_6 y_5- {x}_5 y_6)+\frac{3\omega^2-\lambda _5 }{3} {x}^2_5+\frac{3\omega^2-\lambda _6 }{3}{x}^2_6]+\cdots
\end{aligned}\nonumber
\end{equation}

As a matter of notational convenience, set
\begin{displaymath}
q_0=x_0,q_1=x_5,q_2=x_6,p_0=y_0,p_1=y_5,p_2=y_6.
\end{displaymath}
Then the Hamiltonian $H$ for the three-body problem is
\begin{equation}
H = -\frac{\omega_0^2}{2}+ H_2 + H_3+ H_4+\cdots,
\end{equation}
where
\begin{equation}
\begin{array}{lr}
 H_2 = \frac{1}{2}[p_0^2+ p_1^2+ p_2^2+2\omega _0 \left(p_1 q_2-p_2 q_1\right)+\omega _0^2 q_0^2 +( \omega_0^2-\lambda_5 )  q_1^2+( \omega_0^2-\lambda_5 )  q_2^2 ], \\
  \begin{aligned}
H_3 = -q_0 [\omega_0^2 q_0^2+p_1^2+p_2^2+2\omega_0  (q_2 p_1- q_1 p_2)+( \omega_0^2-\frac{\lambda_5}{2} )q^2_1+( \omega_0^2-\frac{\lambda_6}{2} )q^2_2]-(
a_{30} q_1^3+a_{12} q_1 q_2^2),
\end{aligned} \\
  \begin{aligned}
&H_4 = \frac{1}{2} [(q_2 p_1- q_1 p_2+\omega_0(q^2_1+q^2_2))^2- (q_1 p_1+q_2 p_2)^2]-[  a_{40} q_1^4 +a_{22} q_1^2 q_2^2+a_{04} q_2^4]\nonumber\\
&+q_0[
a_{30} q_1^3+a_{12} q_1 q_2^2]+\frac{3q_0^{2}}{2} [\omega_0^2 q_0^2+p_1^2+p_2^2+2\omega_0  (q_2 p_1- q_1 p_2)+\frac{3\omega_0^2-\lambda _5 }{3} q^2_1+\frac{3\omega_0^2-\lambda _6 }{3}q^2_2].
\end{aligned}
\end{array}\nonumber
\end{equation}

We now look for a change of variables from $(p,q)$ to
$(\mathfrak{p},\mathfrak{q})$ such that $H_2$
 takes the form
\begin{equation}
 \frac{\omega_0(\mathfrak{p}_0^2+\mathfrak{q}_0^2)}{2}+ \frac{\omega_1(\mathfrak{p}_1^2+\mathfrak{q}_1^2)}{2}+\omega_2 \mathfrak{p}_2 \mathfrak{q}_2.\nonumber
\end{equation}

Let $\mathbb{J}$ denote the usual symplectic matrix $\left(
                                                       \begin{array}{cc}
                                                           & -\mathbb{I} \\
                                                         \mathbb{I} &   \\
                                                       \end{array}
                                                     \right)
$.
A straight forward computation shows that  the eigenvalues of the matrix $\mathbb{J}\frac{\partial^2 H_2}{\partial^2 (p,q)}$ are $\pm \omega_0 \mathbf{i},\pm \omega_1 \mathbf{i},\pm\omega_2$.

Note that  we can restrict our
attention to the variables $p_1,p_2,q_1,q_2$.
For the  eigenvalues $\omega_1 \mathbf{i},\omega_2,-\omega_2$, the corresponding
eigenvector are
\begin{equation}
             \begin{array}{lr}
             \left(\frac{\omega_0  \left(2 \lambda_6 +\omega_1^2-3 \lambda\right)}{-2 \lambda_6 +\omega_1^2+3 \lambda},\frac{\mathbf{i} \omega_1 \left(-2 \lambda_6 +\omega_1^2+\lambda\right)}{-2 \lambda_6 +\omega_1^2+3 \lambda},-\frac{2 \mathbf{i} \omega_1 \omega_0 }{-2 \lambda_6 +\omega_1^2+3 \lambda},1\right)^\bot, &  \\
              \left(\frac{\omega_0  \left(-2 \lambda_6 +\omega_2^2+3 \lambda\right)}{2 \lambda_6 +\omega_2^2-3 \lambda},\frac{2 \omega_2 \lambda}{2 \lambda_6 +\omega_2^2-3 \lambda}+\omega_2,\frac{2 \omega_2 \omega_0 }{2 \lambda_6 +\omega_2^2-3 \lambda},1\right)^\bot, &  \\
               \left(\frac{\omega_0  \left(-2 \lambda_6 +\omega_2^2+3 \lambda\right)}{2 \lambda_6 +\omega_2^2-3 \lambda},\omega_2 \left(-\frac{2 \lambda}{2 \lambda_6 +\omega_2^2-3 \lambda}-1\right),-\frac{2 \omega_2 \omega_0 }{2 \lambda_6 +\omega_2^2-3 \lambda},1\right)^\bot.
             \end{array}\nonumber
\end{equation}

It follows that
we can  introduce the following symplectic transformation to reduce the Hamiltonian $H$:
\begin{equation}\label{symplectic transformation}
\left\{
             \begin{array}{lr}
             p_0  =  \sqrt{\omega_0}\mathfrak{p}_0 &  \\
             q_0  =  \frac{\mathfrak{q}_0}{\sqrt{\omega_0}} &  \\
             p_1  = \frac{\omega_0   \left(2 \lambda_6 +\omega_1^2-3 \lambda\right)\mathfrak{q}_1}{\sqrt{r_1} \left(-2 \lambda_6 +\omega_1^2+3 \lambda\right)}+\frac{ \omega_0  \left(-2 \lambda_6 +\omega_2^2+3 \lambda\right)\mathfrak{p}_2}{\sqrt{r_2} \left(2 \lambda_6 +\omega_2^2-3 \lambda\right)}+\frac{\omega_0   \left(-2 \lambda_6 +\omega_2^2+3 \lambda\right)\mathfrak{q}_2}{\sqrt{r_2} \left(2 \lambda_6 +\omega_2^2-3 \lambda\right)}&  \\
              p_2  =  \frac{\omega_1  \left(-2 \lambda_6 +\omega_1^2+\lambda\right)\mathfrak{p}_1}{\sqrt{r_1} \left(-2 \lambda_6 +\omega_1^2+3 \lambda\right)}+\frac{\omega_2 \left(-\frac{2 \lambda}{2 \lambda_6 +\omega_2^2-3 \lambda}-1\right)\mathfrak{p}_2}{\sqrt{r_2}}+\frac{ \left(\frac{2 \omega_2 \lambda}{2 \lambda_6 +\omega_2^2-3 \lambda}+\omega_2\right)\mathfrak{q}_2}{\sqrt{r_2}}&  \\
              q_1  =  -\frac{2 \omega_1  \omega_0 \mathfrak{p}_1}{\sqrt{r_1} \left(-2 \lambda_6 +\omega_1^2+3 \lambda\right)}-\frac{2 \omega_2  \omega_0 \mathfrak{p}_2}{\sqrt{r_2} \left(2 \lambda_6 +\omega_2^2-3 \lambda\right)}+\frac{2 \omega_2 \omega_0  \mathfrak{q}_2}{\sqrt{r_2} \left(2 \lambda_6 +\omega_2^2-3 \lambda\right)}&  \\
              q_2  = \frac{\mathfrak{q}_1}{\sqrt{r_1}}+\frac{\mathfrak{p}_2}{\sqrt{r_2}}+\frac{\mathfrak{q}_2}{\sqrt{r_2}}&
             \end{array}
\right.
\end{equation}
here
\begin{equation}
\begin{array}{lr}
r_1=\frac{\omega_1 \left(-\lambda_6 +2 \omega_1^2-\lambda\right)}{-2 \lambda_6 +\omega_1^2+3 \lambda}, \\
r_2= -\frac{2 \omega_2 \left(\lambda_6 +2 \omega_2^2+\lambda\right)}{2 \lambda_6 +\omega_2^2-3 \lambda}.
\end{array}
\nonumber
\end{equation}

Then the Hamiltonian $H$ becomes
\begin{equation}\label{hamil}
H(\mathfrak{p},\mathfrak{q}) =-\frac{\omega_0^2}{2}+\frac{\omega_0 (\mathfrak{p}_0^2+\mathfrak{q}_0^2)}{2}+ \frac{\omega_1(\mathfrak{p}_1^2+\mathfrak{q}_1^2)}{2}+\omega_2\mathfrak{p}_2 \mathfrak{q}_2+ H_3(\mathfrak{p},\mathfrak{q})+H_4(\mathfrak{p},\mathfrak{q}) +\cdots.
\end{equation}

But we'd better  introduce the following  complex symplectic transformation to reduce the Hamiltonian $H$:

\begin{equation}\label{fu symplectic transformation}
\left\{
             \begin{array}{lr}
             \mathfrak{p} _0=\frac{\zeta_0}{\sqrt{2}}+\frac{\textbf{i} \eta_0}{\sqrt{2}} &  \\
             \mathfrak{q} _0=\frac{\eta_0}{\sqrt{2}}+\frac{\textbf{i} \zeta_0}{\sqrt{2}} &  \\
             \mathfrak{p} _1=\frac{\zeta_1}{\sqrt{2}}+\frac{\textbf{i} \eta_1}{\sqrt{2}}&  \\
              \mathfrak{q} _1=\frac{\eta_1}{\sqrt{2}}+\frac{\textbf{i} \zeta_1}{\sqrt{2}}&  \\
              \mathfrak{p} _2=\zeta_2&  \\
              \mathfrak{q} _2=\eta_2&
             \end{array}
\right.
\end{equation}
Then the Hamiltonian $H$ becomes
\begin{equation}
H(\zeta,\eta)  = -\frac{\omega_0^2}{2}+\textbf{i}\omega _0 \zeta _0 \eta _0 + \textbf{i}\omega _1 \zeta _1 \eta _1 + \omega _2 \zeta _2 \eta _2 + H_3(\zeta,\eta)+H_4(\zeta,\eta) +\cdots.
\nonumber
\end{equation}

We perform the change of variables $(\zeta,\eta)\mapsto (u,v)$
with a generating function
\begin{equation}
 u _0 \eta _0 +u _1 \eta _1 +u _2 \eta _2 + S_3(u,\eta)+S_4(u,\eta) +\cdots,
\nonumber
\end{equation}
 such that in the new variables $(u,v)$ the Hamiltonian function reduces to a Birkhoff normal form of degree 4 up to terms of degree higher than 4:
\begin{equation}
\begin{aligned}
&H(\zeta,\eta) = H(\zeta(u,v),\eta(u,v)) = \mathcal{H}(u,v) \\
&=\textbf{i} \omega _0 u_0 v_0 + \textbf{i}\omega _1 u_1 v _1 + \omega _2 u _2 v_2 -\frac{1}{2}[\omega_{00}(u_0 v_0)^2+\omega_{11}(u_1 v_1)^2+\omega_{22}(u_2 v_2)^2\\
&+2\omega_{01}(u_0 v_0 u_1 v_1)+2\omega_{02}(u_0 v_0 u_2 v_2)+2\omega_{12}(u_1 v_1 u_2 v_2)] +\cdots,
\end{aligned}\nonumber
\end{equation}
 where $S_3$ and $S_4$
are forms of degree 3 and 4 in $u,\eta$,
and
\begin{equation}
\zeta=u+\frac{\partial S_3}{\partial \eta}+\frac{\partial S_4}{\partial \eta}+\cdots,~~~~ ~~~~ ~~~~v=\eta+\frac{\partial S_3}{\partial u}+\frac{\partial S_4}{\partial u}+\cdots.\nonumber
\end{equation}

First, it is straightforward to  show that up to resonance
relation of order 4, $\omega _1,\omega _0$ at most satisfies a resonance
relation of order 3 ($\omega _1=2\omega _0$) for a curve in the space of masses. Fortunately, this resonance
relation does not appear in  the process of obtaining the  Birkhoff normal form. So for every choice of masses, there is no resonance for the Hamiltonian function $H$.

We make use of the relation
\begin{equation}\label{connectrelation}
H(u+\frac{\partial S_3}{\partial \eta}+\frac{\partial S_4}{\partial \eta}+\cdots,\eta) = \mathcal{H}(u,\eta+\frac{\partial S_3}{\partial u}+\frac{\partial S_4}{\partial u}+\cdots)
\end{equation}
to find the  Birkhoff normal form of degree 4.

Equating  the forms of  order 3 in $u,\eta$ of (\ref{connectrelation}) we obtain
\begin{equation}\label{connectrelation3}
\textbf{i} \omega _0 (\frac{\partial S_3}{\partial \eta_0} \eta _0-\frac{\partial S_3}{\partial u_0} u _0) + \textbf{i}\omega _1 (\frac{\partial S_3}{\partial \eta_1} \eta _1-\frac{\partial S_3}{\partial u_1} u _1) + \omega _2 (\frac{\partial S_3}{\partial \eta_2} \eta _2-\frac{\partial S_3}{\partial u_2} u _2) + H_3(u,\eta)= 0.\nonumber
\end{equation}

It follows that $S_3$ can be  determined. Then by equating  the forms of  order 4 in $u,\eta$ of (\ref{connectrelation}) we obtain
\begin{equation}\label{connectrelation4}
\begin{aligned}
&\textbf{i} \omega _0 (\frac{\partial S_4}{\partial \eta_0} \eta _0-\frac{\partial S_4}{\partial u_0} u _0) + \textbf{i}\omega _1 (\frac{\partial S_4}{\partial \eta_1} \eta _1-\frac{\partial S_4}{\partial u_1} u _1) + \omega _2 (\frac{\partial S_4}{\partial \eta_2} \eta _2-\frac{\partial S_4}{\partial u_2} u _2) + H_{3\rightarrow 4}+H_4(u,\eta)\\
&+\frac{1}{2}[\omega_{00}(u_0 v_0)^2+\omega_{11}(u_1 v_1)^2-\omega_{22}(u_2 v_2)^2+2\omega_{01}(u_0 v_0 u_1 v_1)-2\textbf{i}\omega_{02}(u_0 v_0 u_2 v_2)\\
&-2\textbf{i}\omega_{12}(u_1 v_1 u_2 v_2)]=0 ,
\end{aligned}\nonumber
\end{equation}
where $H_{3\rightarrow 4}$ is the forms of  order 4 of $H_3(u+\frac{\partial S_3}{\partial \eta},\eta)$.

It follows that $S_4$ and the  Birkhoff normal form of degree 4 can be  determined.

However, we remark that it is only necessary to determine $\omega_{00},\omega_{01},\omega_{11}$ for the  Birkhoff normal form on the center manifold.

For the Birkhoff normal form of degree 4 on the center manifold:
\begin{equation}
\begin{aligned}
&\mathcal{H}(u,v)
=\textbf{i} \omega _0 u_0 v_0 + \textbf{i}\omega _1 u_1 v _1  -\frac{1}{2}[\omega_{00}(u_0 v_0)^2+\omega_{11}(u_1 v_1)^2+2\omega_{01}(u_0 v_0 u_1 v_1)] +\cdots,
\end{aligned}\nonumber
\end{equation}
we can switch to action-angle variables, the above Birkhoff normal form becomes
\begin{equation}\label{Birkhoff normal form}
\begin{aligned}
&\mathcal{H}(\varrho,\varphi)
= \omega _0 \varrho_0 + \omega _1 \varrho_1  +\frac{1}{2}[\omega_{00}\varrho_0 ^2+\omega_{11}\varrho_1 ^2+2\omega_{01}\varrho_0 \varrho_1 ] +\cdots,
\end{aligned}
\end{equation}
where $\varrho_j=\textbf{i}u_j v_j$ ($j=0,1$) are action variables, and
\begin{equation}
\omega _{00}=-3,\nonumber
\end{equation}
\begin{equation}
\omega _{01}=-\frac{(2 \tau -1) \sqrt{\frac{2 \tau ^2+7 \tau -4}{\tau }} \left(7 \tau ^4+42 \tau ^3+41 \tau ^2-76 \tau -28\right)}{(\tau -1) \left(2 \tau ^2+7 \tau +2\right) \left(8 \tau ^2+19 \tau -16\right)}, \nonumber
\end{equation}
\begin{equation}
\begin{aligned}
&\omega _{11}=\frac{-16 \tau ^8-125 \tau ^7+58 \tau ^6+1436 \tau ^5+454 \tau ^4-2739 \tau ^3+2640 \tau ^2-1140 \tau +128}{2 (\tau -1)^2 \tau  (\tau +4) (2 \tau -1) \left(8 \tau ^2+19 \tau -16\right)}\\
&+\frac{27 \tau  \left(\frac{9 a_{30}^2 \tau  \left(28 \tau ^6+131 \tau ^5+434 \tau ^4+544 \tau ^3-472 \tau ^2+140 \tau -120\right)}{(4 \tau +1) \left(4 \tau ^2+7 \tau -6\right)}-a_4 \lambda  (\tau +4) \left(3 \tau ^4+8 \tau ^3+24 \tau ^2+8\right)\right)}{8 \lambda ^2 (\tau +4)^2 (2 \tau -1) \left(\tau ^2-1\right)^2},
\end{aligned}\nonumber
\end{equation}
where
\begin{equation}
\tau =\frac{1}{4} \left(5-\frac{9 \lambda _6}{\lambda }\right)+\frac{3}{4} \sqrt{9 \left(\frac{\lambda _6}{\lambda }\right){}^2-\frac{10 \lambda _6}{\lambda }+1}. \nonumber
\end{equation}

Then
\begin{equation}
\begin{aligned}
&\omega _{00}\omega _{11}-\omega _{01}^2=-\frac{(\tau +4) (2 \tau -1)^3 \left(7 \tau ^4+42 \tau ^3+41 \tau ^2-76 \tau -28\right)^2}{(\tau -1)^2 \tau  \left(2 \tau ^2+7 \tau +2\right)^2 \left(8 \tau ^2+19 \tau -16\right)^2}\\
&+\frac{3 }{8 \lambda ^2 \tau  (\tau +4)^2 (2 \tau -1) (4 \tau +1) \left(\tau ^2-1\right)^2 \left(4 \tau ^2+7 \tau -6\right) \left(8 \tau ^2+19 \tau -16\right)}\\
&[-243 a_{30}^2  \tau ^3 (224 \tau ^8+1580 \tau ^7+5513 \tau ^6+10502 \tau ^5-384 \tau ^4-16552 \tau ^3+9252 \tau ^2-4520 \tau \\
&+1920)+27 a_4 \lambda  \tau ^2(384 \tau ^{10}+4240 \tau ^9+19016 \tau ^8+45319 \tau ^7+49694 \tau ^6-16688 \tau ^5\\
&-54352 \tau ^4+20392 \tau ^3-17776 \tau ^2+5824 \tau +3072) +4 \lambda ^2 (\tau +1)^2 (256 \tau ^{12}+3536 \tau ^{11}\\
&+12848 \tau ^{10}-15853 \tau ^9-161192 \tau ^8-157864 \tau ^7+277966 \tau ^6+136889 \tau ^5-377438 \tau ^4\\
&+243876 \tau ^3-35208 \tau ^2-17888 \tau +3072)].
\end{aligned}\nonumber
\end{equation}

A straight forward computation shows that for every choice of masses, $\omega _{00}\omega _{11}-\omega _{01}^2 \neq 0$.  Therefore, the reduced Hamiltonian  on the $4$-dimensional central manifold  is nondegenerate.
Aa a result, it follows from  Theorem \ref{Conley and Zendergeneral}  that
\begin{theorem}\label{Conley and Zenderapply2}
For every choice of positive masses of the planar three-body problem, there are abundant periodic orbits near every Euler relative equilibrium, and the relative measure of the closure of the set of periodic orbits in  the polydisc $\|\varrho\|<\epsilon$ on the $4$-dimensional center manifold is
at least $1-O(\epsilon^{\frac{1}{4}})$.
\end{theorem}

Furthermore, one can show that the set  of masses  such that the frequency vector $\varpi=(\omega _0,\omega _1)$ is resonant of order $m>4$  is  a real algebraic variety in the spaces of masses. In fact, it suffices to observe that $\omega _1/\omega _0$ is an algebraic function of masses. Then it follows that generically frequency vectors $\varpi=(\omega _0,\omega _1)$ are nonresonant up to order $m>4$. According to Remark \ref{KAM1}, it follows that, for any $m>4$, the relative measure in the above theorem is
at least $1-O(\epsilon^{\frac{m-3}{4}})$ generically.

\section{Periodic Orbits Near Euler-Moulton Relative Equilibria of the $N$-body Problem}
\indent\par
It's natural to conjecture that generically there are abundant periodic orbits near a relative equilibrium solution of the general $N$-body problem, although it's not easy to prove this. But as we would intuitively expect, for  example, suppose that the
quadratic part $H_2$ of the Hamiltonian for a relative equilibrium solution generated by a nondegenerate central configuration $\mathcal{E}_3$ has $d$ imaginary eigenvalues. Then the elements such as $\lambda, \lambda_k, \mathcal{E}_k, \omega_{jk}$ etc would smoothly depend on the masses, and  the
quadratic part $H_2$ corresponding to  central configurations near $\mathcal{E}_3$ would generically have $d$ unequal imaginary eigenvalues which are nonresonant up to order $4$; and then, the Hamiltonian corresponding to  central configurations near $\mathcal{E}_3$ would generically be nondegenerate, as a result, it should prove  that generically there are abundant periodic orbits near the relative equilibrium on the $2d$-dimensional center manifold by Theorem \ref{Conley and Zendergeneral}. Unfortunately, we cannot rigorously prove this  intuitive view.

However, in this section we will prove rigorously that there are abundant periodic orbits near every Euler-Moulton relative equilibrium of the $N$-body Problem on the $2N-2$-dimensional center manifold.
\begin{theorem}\label{Conley and Zenderapply3}
For almost every choice of positive masses of the planar $N$-body problem, there are abundant periodic orbits near the corresponding Euler-Moulton relative equilibrium, and the relative measure of the closure of the set of periodic orbits in  the polydisc $\|\varrho\|<\epsilon$ on the $2N-2$-dimensional center manifold is
at least $1-O(\epsilon^{\frac{1}{4}})$.
\end{theorem}

By Theorem \ref{Conley and Zendergeneral}, the key point is judging that the reduced Hamiltonian of the planar $N$-body problem for every Euler-Moulton relative equilibrium on the center manifold is nondegenerate. Because  Euler-Moulton central configurations, i.e., collinear central configurations are nondegenerate, the elements such as $\lambda,  \lambda_k, \mathcal{E}_k, \omega_{jk}$ etc smoothly depend on the masses on the whole, indeed, algebraically depend on the masses on the whole. As a result, it is easy to believe that the corresponding reduced Hamiltonian  on the center manifold is nondegenerate except a proper algebraic subset of the mass space. However, one can not simply claim Theorem \ref{Conley and Zenderapply3} is correct, since we can not simply exclude constant value functions, or more precisely, the elements such as the frequencies $\omega_k$ might be constant and resonant for all the masses of the $N$-body Problem, then we  cannot even obtain the Birkhoff normal form of degree 4; even if the frequencies $\omega_k$ are not resonant and we could obtain the Birkhoff normal form of degree 4, the determinant  $det(\omega_{jk})$ on the center manifold might be invariably zero, that is, the reduced Hamiltonian  on the center manifold is degenerate for all the masses of the $N$-body Problem. This it's necessary to prove rigorously that the elements such as $\omega_k,  det(\omega_{jk})$, as algebraic functions of  the masses, are not constant. We will prove this by the perturbed  method and the inductive method.

We remark that  in the following it's shown that the
quadratic part $H_2$ of the reduced Hamiltonian  on the $2N-2$-dimensional central manifold is always positive definite (see the following (\ref{firststep})). So by Theorem \ref{weinstein}, one can  prove that there are $N-1$ one parameter family of periodic orbits  near every  Euler-Moulton relative equilibrium.

\subsection{Collinear Central Configurations}
\indent\par
As a preliminary  to the following sections, let's recall the results of collinear central configurations.

Suppose $\mathcal{E}_3=\mathbf{r} = (\xi_1, 0, \xi_2, 0, \cdots, \xi_{N}, 0)^\top\in (\mathbb{R}\times 0)^{N} \subset \mathbb{R}^{2N}$ is a collinear central configuration, then the matrix $B_{jk} = \frac{m_j m_k}{r^3_{jk}}\left(
                                    \begin{array}{cc}
                                      -2 & 0 \\
                                      0 & 1 \\
                                    \end{array}
                                  \right)
$, so $\mathfrak{M}^{-1} \mathfrak{B}$ becomes:
\begin{center}
$ \left(
       \begin{array}{cccc}
         A_{11} &\frac{ m_2}{r^3_{12}} D &\cdots & \frac{ m_N}{r^3_{1N}} D \\
         \frac{ m_1}{r^3_{12}} D &A_{11} &\cdots & \frac{ m_N}{r^3_{2N}} D \\
         \vdots &\vdots & \ddots & \vdots \\
         \frac{ m_1}{r^3_{1N}} D & \frac{ m_2}{r^3_{2N}} D& \cdots & A_{NN}\\
       \end{array}
     \right)
$
\end{center}
where $D=\left(
                                    \begin{array}{cc}
                                      -2 & 0 \\
                                      0 & 1 \\
                                    \end{array}
                                  \right)$
and the diagonal blocks are given by:
\begin{displaymath}
A_{kk} = -\sum_{1\leq j\leq N, j\neq k} \frac{ m_j}{r^3_{jk}} D.
\end{displaymath}

\cite{yu2019stability} has pointed out that
$\{\mathcal{E}_5, \mathcal{E}_6,  \cdots, \mathcal{E}_{2N-1}, \mathcal{E}_{2N}\}$ can be considered as
\begin{displaymath}
\{\mathcal{E}_5,\mathcal{E}_6=\mathcal{E}_5^\bot,  \cdots,\mathcal{E}_{2N-1},\mathcal{E}_{2N}=\mathcal{E}_{2N-1}^\bot\},
\end{displaymath} then $Q$ becomes block diagonal with  block $J=\left(
           \begin{array}{cc}
             0 & -1 \\
             1 & 0 \\
           \end{array}
         \right)$:
         \begin{equation}\label{qjk}
Q = \left(
       \begin{array}{ccc}

         J &   &  \\
           & \ddots &   \\
          &   & J\\
       \end{array}
     \right)
\end{equation}
We also have
\begin{equation}\label{lineareigenvalues}
\lambda_{2k}=\frac{3\sqrt{g_3}\lambda -\lambda_{2k-1}}{2}<0, ~~~~~~~\lambda_{2k-1} >3\sqrt{g_3}\lambda, ~~~~~~~ for ~~~ 3\leq k \leq N.
\end{equation}

For simplicity, we confine ourselves to the configuration space $\mathbb{R}^{N}$ and  suppose $\mathbf{e}_2 = (\xi_1, \xi_2,\cdots, \xi_{N})^\top\in  \mathbb{R}^{N}$ is a collinear central configuration corresponding to $\mathcal{E}_3$ in the following.

The equations (\ref{centralconfiguration}) of central configurations become:
\begin{equation}\label{collinearcentralconfiguration}
\sum_{j=1,j\neq k}^N \frac{m_j}{|\mathbf{r}_j-\mathbf{r}_k|^3}(\mathbf{r}_j-\mathbf{r}_k)=-\lambda \mathbf{r}_k,1\leq k\leq N,
\end{equation}
here $\mathbf{r}_k=\xi_k \in  \mathbb{R}$.

Without loss of generality, we fix $\lambda=1$, $m_1=1$ and assume $\xi_1<\xi_2<\cdots<\xi_N$.
Then the equations (\ref{collinearcentralconfiguration}) are $N$ algebraic equations of $N$ Unknowns  $\mathbf{r}_k \in  \mathbb{R}$ ($k=1,2,\cdots,N$) with $N-1$ parameters $m_k $ ($k=2,\cdots,N$). The equations (\ref{collinearcentralconfiguration}) are independent at the central configuration $\mathbf{e}_2$ if and only if the matrix
\begin{center}
$ \mathcal{D}= \left(
       \begin{array}{cccc}
         1+\sum_{1\leq j\leq N, j\neq 1} \frac{ 2m_j}{r^3_{j1}} &-\frac{ 2m_2}{r^3_{12}}  &\cdots & -\frac{ 2m_N}{r^3_{1N}}  \\
         -\frac{ 2m_1}{r^3_{12}}  &1+\sum_{1\leq j\leq N, j\neq 2} \frac{ 2m_j}{r^3_{j2}} &\cdots & -\frac{ 2m_N}{r^3_{2N}}  \\
         \vdots &\vdots & \ddots & \vdots \\
         -\frac{ 2m_1}{r^3_{1N}}  & -\frac{ 2m_2}{r^3_{2N}} & \cdots & 1+\sum_{1\leq j\leq N, j\neq N} \frac{ 2m_j}{r^3_{jN}}\\
       \end{array}
     \right)
$
\end{center}
is nonsingular at $\mathbf{e}_2$; and  the central configuration  $\mathbf{e}_2$ is degenerate if and only if the matrix $ \mathcal{D}$ is degenerate at $\mathbf{e}_2$.

When masses parameters $m_k>0 $ ($k=2,\cdots,N$), it's easy to see that all the eigenvalues  of $\mathcal{D}$  are positive, this is also true even for one of  masses parameters $m_k $ ($k=2,\cdots,N$) being zero (and for some of  masses parameters $m_k $ ($k=2,\cdots,N$) being zero under some additional conditions). Therefore collinear central configuration are all nondegenerate, furthermore, the elements $\xi_k$ ($k=1,\cdots,N$) (or the vector $\mathbf{e}_2$) are algebraic functions of masses parameters.

Since the matrix $\mathcal{D}$ is symmetric linear mapping with respect to the scalar product $\langle,\rangle$, there are $N$ orthogonal eigenvectors $\mathbf{e}_1,\mathbf{e}_2,\cdots,\mathbf{e}_N$ of $\mathcal{D}$  with respect to the scalar product $\langle,\rangle$.
The eigenvectors $\mathbf{e}_1,\mathbf{e}_2,\cdots,\mathbf{e}_N$ and the corresponding eigenvalues $\iota_1,\iota_2,\cdots,\iota_N$ are also algebraic functions of masses parameters. One can fix $\mathbf{e}_1 = (1, 1,\cdots, 1)^\top$ and $\iota_1= 1= \lambda$. The eigenvalue $\iota_2= 3 \lambda=3$.

It's noteworthy that $\iota_k(\mathbf{e}_2)/\iota_1(\mathbf{e}_2)$ or $\iota_k(\mathbf{e}_2)/\lambda(\mathcal{E}_3)$ ($k=1,\cdots,N$) do not depend on the scale $\sqrt{g_3}$ of the central configuration $\mathbf{e}_2$. The similar is true for $\frac{\lambda_k (\mathcal{E}_3)}{\sqrt{g_3}\lambda(\mathcal{E}_3)}$. Indeed,
it's obvious that for $k=3,\cdots,N$, we have
\begin{displaymath}
\begin{array}{c}
  \lambda_{2k-1}(\mathcal{E}_3)= \sqrt{g_3}\iota_k (\mathbf{e}_2),\\
  \lambda_{2k}(\mathcal{E}_3)=\frac{3 -\iota_k(\mathbf{e}_2)/\lambda(\mathcal{E}_3)}{2}\sqrt{g_3}\lambda(\mathcal{E}_3).
\end{array}
\end{displaymath}
In particular, we have
\begin{equation}\label{lambdalambdak}
\left\{
             \begin{array}{lr}
             \lambda^*=g_3^{\frac{3}{2}}\lambda=g_3^{\frac{3}{2}}=\|\mathbf{e}_2\|^3, &\\
             \lambda_{2k-1}^*= \iota_k \lambda^*, &k=3,\cdots,N\\
             \lambda_{2k}^*=\frac{3 -\iota_k}{2}\lambda^*, &k=3,\cdots,N
             \end{array}
\right.
\end{equation}

\subsection{Induction  on Collinear Central Configurations}
\indent\par
In this subsection we will obtain some information of the elements  $\iota_k,\mathbf{e}_k$ ($k=1,\cdots,N$) by the perturbed method and the inductive method.

\textbf{A. The two-body case.} First, let's consider the two-body problem: $N=2$.

By the equations of central configurations
\begin{equation}\label{collinearcentralconfigurationtwo-body}
\left\{
             \begin{array}{lr}
             \frac{ m_2(\xi_2-\xi_1)}{r^3_{12}}=-\xi_1, &\\
             \frac{ m_1(\xi_1-\xi_2)}{r^3_{12}}=-\xi_2, &
             \end{array}
\right.
\end{equation}
we have
\begin{displaymath}
\mathbf{e}_2 = (\xi_1, \xi_2)^\top~~ ~~~~ ~~~~ ~~~~ ~~r_{12}=\xi_2-\xi_1=\sqrt[3]{m_1+m_2},
\end{displaymath}
where
\begin{displaymath}
\xi_1=-\frac{m_2}{({m_1+m_2})^{\frac{2}{3}}} ~~ ~~~~ ~~~~ ~~~~ ~~\xi_2=\frac{m_1}{({m_1+m_2})^{\frac{2}{3}}}.
\end{displaymath}

Suppose $m_2=\varepsilon_1$ is small, then
\begin{equation}\label{two-body}
\begin{array}{c}
  \mathbf{e}_2(\varepsilon_1) = (0, 1)^\top+O(\varepsilon_1), \\
  r_{12}(\varepsilon_1)=\xi_2-\xi_1=\sqrt[3]{1+\varepsilon_1}=1+c_1\varepsilon_1+O(\varepsilon_1^2),
\end{array}
\end{equation}
where $c_1=\frac{1}{3}$.

\textbf{B. The  three-body case.} Let's firstly consider the restricted three-body problem: $N=3, m_3=0$.

The equations (\ref{collinearcentralconfiguration}) of central configurations reduce to the equations (\ref{collinearcentralconfigurationtwo-body}) of central configurations for the two-body problem and the equation
\begin{equation}\label{collinearcentralconfigurationrestrictedthree-body}
 \frac{m_1}{(r_{12}+r_{23})^3}(\xi_1-\xi_3)+\frac{m_2}{r_{23}^3}(\xi_2-\xi_3) = - \xi_3.\nonumber
\end{equation}
Then $r_{12}$ is same as that in (\ref{two-body}) and we further have
\begin{displaymath}
\frac{\varepsilon_1}{r_{23}^3}=1+\frac{2}{r_{12}^3}+ O(r_{23}),
\end{displaymath}
or
\begin{displaymath}
r_{23}(\varepsilon_1)=c_2 \varepsilon_1^ {\frac{1}{3}}+ c_2^5 \varepsilon_1^{\frac{2}{3}}+O(\varepsilon_1),
\end{displaymath}
where $c_2=\frac{1}{\sqrt[3]{3\iota_1}}=\frac{1}{\sqrt[3]{3}}$.

By
\begin{center}
$ \mathcal{D}= \left(
       \begin{array}{ccc}
         1+\frac{ 2m_2}{r^3_{12}} &-\frac{ 2m_2}{r^3_{12}}  & 0  \\
         -\frac{ 2m_1}{r^3_{12}}  &1+\frac{ 2m_1}{r^3_{12}} & 0  \\
         -\frac{ 2m_1}{r^3_{13}}  & -\frac{ 2m_2}{r^3_{23}} &  1+\frac{ 2m_1}{r^3_{13}}  +\frac{ 2m_2}{r^3_{23}}\\
       \end{array}
     \right)
$
\end{center}
it follows that
\begin{displaymath}
\begin{array}{c}
  \iota_3(\varepsilon_1)= 1+\frac{ 2m_1}{r^3_{13}}  +\frac{ 2m_2}{r^3_{23}}=3\iota_2-12c_2 \varepsilon_1^ {\frac{1}{3}}+ O(\varepsilon_1^{\frac{2}{3}})=9-12c_2 \varepsilon_1^ {\frac{1}{3}}+ O(\varepsilon_1^{\frac{2}{3}}), \\
  \mathbf{e}_2(\varepsilon_1) =(\xi_1, \xi_2, \xi_3)^\top= (0, 1,1)^\top+(0, 0,c_2 \varepsilon_1^ {\frac{1}{3}})^\top+ O(\varepsilon_1^{\frac{2}{3}}),\\
  \mathbf{e}_3 (\varepsilon_1)= (0, 0,1)^\top.
\end{array}
\end{displaymath}

Let's further consider the  three-body problem: $N=3, m_3>0$.
Suppose $m_3=\varepsilon_2 = O (\varepsilon_1^ {100})$ is small, then it follows from the implicit function theorem that
\begin{equation}\label{three-body1}
\begin{array}{c}
  \iota_3(\varepsilon_1,\varepsilon_2)= \iota_3(\varepsilon_1)+O(\varepsilon_2)=3^2\iota_1-12c_2 \varepsilon_1^ {\frac{1}{3}}+ O(\varepsilon_1^{\frac{2}{3}}), \\
  \mathbf{e}_2 (\varepsilon_1,\varepsilon_2)= \mathbf{e}_2(\varepsilon_1)+O(\varepsilon_2)=(0, 1,1)^\top+(0, 0,c_2 \varepsilon_1^ {\frac{1}{3}})^\top+ O(\varepsilon_1^{\frac{2}{3}}),\\
  \mathbf{e}_3 (\varepsilon_1,\varepsilon_2)=\mathbf{e}_3(\varepsilon_1)+O(\varepsilon_2)= (0, 0,1)^\top +  O(\varepsilon_2).
\end{array}
\end{equation}
In particular, we have
\begin{equation}\label{three-body2}
\begin{array}{c}
  r_{12}(\varepsilon_1,\varepsilon_2)=r_{12}(\varepsilon_1)+O(\varepsilon_2)=1+c_1\varepsilon_1+O(\varepsilon_1^2),\\
  r_{23}(\varepsilon_1,\varepsilon_2)=r_{23}(\varepsilon_1)+O(\varepsilon_2)=c_2 \varepsilon_1^ {\frac{1}{3}}+ c_2^5 \varepsilon_1^{\frac{2}{3}}+O(\varepsilon_1).
\end{array}
\end{equation}

\textbf{C. The  four-body case.} Let's firstly consider the restricted four-body problem: $N=4, m_4=0$.

The equations (\ref{collinearcentralconfiguration}) of central configurations reduce to  the equations of central configurations for the three-body problem and the equation
\begin{equation}\label{collinearcentralconfigurationrestrictedfour-body}
\frac{m_1}{(r_{12}+r_{23}+r_{34})^2}+\frac{m_2}{(r_{23}+r_{34})^2}+\frac{m_3}{r_{34}^2} =  \xi_4.\nonumber
\end{equation}
Then $r_{12}, r_{23}$ are  same as that in (\ref{three-body2}) and we further have
\begin{displaymath}
\frac{\varepsilon_2}{r_{34}^3}=1+\frac{2}{(r_{12}+r_{23})^3}+\frac{2\varepsilon_1}{r_{23}^3}+ O(r_{34}),
\end{displaymath}
or
\begin{displaymath}
r_{34}(\varepsilon_1,\varepsilon_2)=c_3 \varepsilon_2^ {\frac{1}{3}}+ O((\varepsilon_1\varepsilon_2)^{\frac{1}{3}}),
\end{displaymath}
where $c_3=\frac{1}{\sqrt[3]{3^2\iota_1}}=\frac{1}{\sqrt[3]{3^2}}$.

By
\begin{center}
$ \mathcal{D}= \left(
       \begin{array}{cccc}
         1+\frac{ 2m_2}{r^3_{12}}+\frac{ 2m_3}{r^3_{13}} &-\frac{ 2m_2}{r^3_{12}}  & -\frac{ 2m_3}{r^3_{13}} &0   \\
         -\frac{ 2m_1}{r^3_{12}}  &1+\frac{ 2m_1}{r^3_{12}}+\frac{ 2m_3}{r^3_{23}} &-\frac{ 2m_3}{r^3_{23}} & 0  \\
         -\frac{ 2m_1}{r^3_{13}}  & -\frac{ 2m_2}{r^3_{23}} &  1+\frac{ 2m_1}{r^3_{13}}  +\frac{ 2m_2}{r^3_{23}}&0\\
         -\frac{ 2m_1}{r^3_{14}}  & -\frac{ 2m_2}{r^3_{24}} & -\frac{ 2m_3}{r^3_{34}}& 1+\frac{ 2m_1}{r^3_{14}}  +\frac{ 2m_2}{r^3_{24}}+\frac{ 2m_3}{r^3_{34}}\\
       \end{array}
     \right)
$
\end{center}
it follows that
\begin{displaymath}
\begin{array}{c}
  \iota_4(\varepsilon_1,\varepsilon_2)= 1+\frac{ 2m_1}{r^3_{14}}  +\frac{ 2m_2}{r^3_{24}}+\frac{ 2m_3}{r^3_{34}}=3 \iota_3(\varepsilon_1,\varepsilon_2)+O(\varepsilon_2^{\frac{1}{3}})=27-36c_2 \varepsilon_1^ {\frac{1}{3}}+ O(\varepsilon_1^{\frac{2}{3}}), \\
  \mathbf{e}_2(\varepsilon_1,\varepsilon_2) = (\xi_1, \xi_2, \xi_3,\xi_4)^\top=(0, 1,1,1)^\top+(0, 0,c_2 \varepsilon_1^ {\frac{1}{3}},c_2 \varepsilon_1^ {\frac{1}{3}})^\top+ O(\varepsilon_1^{\frac{2}{3}}),\\
  \mathbf{e}_4 (\varepsilon_1,\varepsilon_2)= (0, 0,0,1)^\top.
\end{array}
\end{displaymath}

Let's further consider the  three-body problem: $N=3, m_3>0$.

Suppose $m_4=\varepsilon_3 = O (\varepsilon_2^ {100})$ is small, then it follows from the implicit function theorem that
\begin{equation}\label{four-body1}
\begin{array}{c}
\iota_3(\varepsilon_1,\varepsilon_2,\varepsilon_3)=\iota_3(\varepsilon_1,\varepsilon_2)+O(\varepsilon_3)= 3^2\iota_1-12c_2 \varepsilon_1^ {\frac{1}{3}}+ O(\varepsilon_1^{\frac{2}{3}}), \\
  \iota_4(\varepsilon_1,\varepsilon_2,\varepsilon_3)=\iota_4(\varepsilon_1,\varepsilon_2)+O(\varepsilon_3)=27-36c_2 \varepsilon_1^ {\frac{1}{3}}+ O(\varepsilon_1^{\frac{2}{3}}), \\
 \mathbf{e}_2 (\varepsilon_1,\varepsilon_2,\varepsilon_3)= (\xi_1, \xi_2, \xi_3,\xi_4)^\top = \mathbf{e}_2 (\varepsilon_1,\varepsilon_2)+O(\varepsilon_3)=(0, 1,1+c_2 \varepsilon_1^ {\frac{1}{3}},1+c_2 \varepsilon_1^ {\frac{1}{3}})^\top+ O(\varepsilon_1^{\frac{2}{3}}),\\
 \mathbf{e}_4 (\varepsilon_1,\varepsilon_2,\varepsilon_3)=\mathbf{e}_4 (\varepsilon_1,\varepsilon_2)+O(\varepsilon_3)= (0, 0,0,1)^\top+O(\varepsilon_3),
\end{array}\nonumber
\end{equation}
and
\begin{equation}\label{four-body2}
\begin{array}{c}
 r_{12}(\varepsilon_1,\varepsilon_2,\varepsilon_3)= r_{12}(\varepsilon_1,\varepsilon_2)+O(\varepsilon_3)=r_{12}(\varepsilon_1)+O(\varepsilon_2)=1+c_1\varepsilon_1+O(\varepsilon_1^2),\\
r_{23}(\varepsilon_1,\varepsilon_2,\varepsilon_3)=  r_{23}(\varepsilon_1,\varepsilon_2)+O(\varepsilon_3)=r_{23}(\varepsilon_1)+O(\varepsilon_2)=c_2 \varepsilon_1^ {\frac{1}{3}}+ c_2^5 \varepsilon_1^{\frac{2}{3}}+O(\varepsilon_1),\\
r_{34}(\varepsilon_1,\varepsilon_2,\varepsilon_3)=  r_{34}(\varepsilon_1,\varepsilon_2)+O(\varepsilon_3)=c_3 \varepsilon_2^ {\frac{1}{3}}+ O((\varepsilon_1\varepsilon_2)^{\frac{1}{3}}).
\end{array}\nonumber
\end{equation}

\textbf{D. The  $N$-body case.} In short, an easy induction gives
\begin{equation}\label{n-body1}
\begin{array}{c}
\iota_{k}(\varepsilon_1,\cdots,\varepsilon_{n})=\iota_{k}(\varepsilon_1,\cdots,\varepsilon_{n-1})+O(\varepsilon_{n}),~~~~~n+1\geq k \geq 3;\\
  \iota_{n+2}(\varepsilon_1,\cdots,\varepsilon_{n})= 3\iota_{n+1}(\varepsilon_1,\cdots,\varepsilon_{n})+O(\varepsilon_{n}^ {\frac{1}{3}}), ~~~~~n\geq 1;\\
  \mathbf{e}_k (\varepsilon_1,\cdots,\varepsilon_{n})= \mathbf{e}_k (\varepsilon_1,\cdots,\varepsilon_{n-1})+O(\varepsilon_n),~~~~~n+1\geq k\geq 2;\\
  \mathbf{e}_{n+2} (\varepsilon_1,\cdots,\varepsilon_{n})= (\underbrace{0, \cdots,0}_{n+1},1)^\top ,~~~~~n\geq 1;\\
  r_{n+1,n+2}(\varepsilon_1,\cdots,\varepsilon_{n})=c_{n+1} \varepsilon_{n}^ {\frac{1}{3}}+O((\varepsilon_{n-1}\varepsilon_n)^{\frac{1}{3}}),~~~~~n\geq 2;\\
  c_{n+1}=\frac{1}{\sqrt[3]{3^n\iota_1}}=\frac{1}{\sqrt[3]{3^n}},~~~~~n\geq 1.
\end{array}
\end{equation}
To sum up, we have
\begin{equation}\label{n-body1}
\begin{array}{c}
\iota_{1}=1,\\
\iota_{2}=3,\\
  \iota_{n}(\varepsilon_1,\cdots,\varepsilon_{N-1})=3^{n-1}(1-\frac{4}{3}c_2 \varepsilon_1^ {\frac{1}{3}})+ O(\varepsilon_1^{\frac{2}{3}}), ~~~~~N\geq n>2;
\end{array}
\end{equation}
and
\begin{equation}\label{n-body2}
\begin{array}{c}
  \mathbf{e}_{N} (\varepsilon_1,\cdots,\varepsilon_{N-1})= (\underbrace{0, \cdots,0}_{N-1},1)^\top+O(\varepsilon_{N-1}),\\
  \mathbf{e}_{2} (\varepsilon_1,\cdots,\varepsilon_{N-1})=(\xi_1,\cdots,\xi_{N})^\top,
\end{array}
\end{equation}
where
\begin{equation}\label{n-body3}
 \xi_{n+1}-\xi_n=r_{n,n+1},~~~~~N-1\geq n\geq1,
\end{equation}
and
\begin{equation}\label{n-body}
\begin{array}{c}
  r_{12}(\varepsilon_1,\cdots,\varepsilon_{N-1})=1+c_1\varepsilon_1+O(\varepsilon_1^2),\\
  r_{23}(\varepsilon_1,\cdots,\varepsilon_{N-1})=  c_2 \varepsilon_1^ {\frac{1}{3}}+ c_2^5 \varepsilon_1^{\frac{2}{3}}+O(\varepsilon_1),\\
  r_{n,n+1}(\varepsilon_1,\cdots,\varepsilon_{N-1})=c_n \varepsilon_{n-1}^ {\frac{1}{3}}+O((\varepsilon_{n-2}\varepsilon_{n-1})^{\frac{1}{3}}), ~~~~~N-1\geq n\geq 3.
\end{array}\nonumber
\end{equation}

\subsection{The Resonance of  Frequencies}
\indent\par
In this subsection, we will prove that the frequencies of the reduced Hamiltonian  on the center manifold are nonresonant up to order $4$ generically, and thus, we could obtain the Birkhoff normal form of degree 4 for all the masses of the $N$-body Problem except at most  a proper algebraic subset of the mass space.

Recall that the Hamiltonian is
\begin{eqnarray}\label{Hamiltonian3Euler}
&&H(x_0,x,y_0,y)=\frac{-\omega^2}{2} +\frac{y_0^2}{2}+ \frac{1}{2} [ \omega^2 x_0^2+\sum_{k=5}^{2N}y_k^2-2\omega\sum_{j,k=5}^{2N}q_{kj} {x}_j y_k+\sum_{k=5}^{2N}( \omega^2-\lambda^*_k ){x}^2_k] \nonumber\\
&-&[x_0 (\omega^2 x_0^2+\sum_{k=5}^{2N}y_k^2-2\omega\sum_{j,k=5}^{2N}q_{kj} {x}_j y_k+\sum_{k=5}^{2N}( \omega^2-\frac{\lambda^*_k}{2} ){x}^2_k)+\frac{1}{6}\sum_{i,j,k = 5}^{2N} a_{ijk} x_i x_j x_k ]\nonumber\\
&+& \frac{1}{2} [(\sum_{j,k=5}^{2N}q_{kj}{x}_j y_k-\omega\sum_{k=5}^{2N}{x}^2_k)^2- (\sum_{k=5}^{2N}x_k y_k)^2]+ \frac{x_0}{6}\sum_{i,j,k = 5}^{2N} a_{ijk} x_i x_j x_k\nonumber\\
&+& \frac{3x_0^2}{2} [\omega^2 x_0^2+\sum_{k=5}^{2N}y_k^2-2\omega\sum_{j,k=5}^{2N}q_{kj} {x}_j y_k+\sum_{k=5}^{2N}( \omega^2-\frac{\lambda^*_k}{3} ){x}^2_k]\nonumber\\
&-&  [ \frac{3 \omega^2}{8} (\sum_{k = 5}^{2N}  x^2_k)^2+\frac{3}{4} (\sum_{j = 5}^{2N}  x^2_j)\sum_{k = 5}^{2N} (\lambda^*_k - \omega^2)  x^2_k]-\frac{1}{24}\sum_{h,i,j,k = 5}^{2N} a_{hijk} x_h x_i x_j x_k+\cdots,
\end{eqnarray}
here $\omega=\sqrt{\lambda^*} =  g_3^{\frac{3}{4}}$ and $\lambda^*,\lambda^*_k$ satisfy the relations (\ref{lambdalambdak}).

By (\ref{qjk}), the
quadratic part $H_2$ of the Hamiltonian (\ref{Hamiltonian3Euler}) is
\begin{equation}\label{Hamiltonian2H2Euler}
\begin{aligned}
&H_2=\frac{1}{2}[y_0^2+\omega^2 x_0^{2}]\\
&+ \sum_{k = 3}^{N}\frac{1}{2} [y_{2k-1}^2+y_{2k}^2+2\omega ({x}_{2k} y_{2k-1}-{x}_{2k-1} y_{2k})+\omega^2(1 -\iota_k){x}^2_{2k-1}+\omega^2( \frac{\iota_k-1 }{2}){x}^2_{2k}],
\end{aligned}
\end{equation}
it's easy to see that the $2N-2$ eigenvalues are
\begin{equation}\label{eigenvalues1}
  \pm \omega_0 \mathbf{i}, ~~~~\pm \omega_{2k-1} \mathbf{i}, ~~~~\pm \omega_{2k} , ~~~~ k=3,\cdots,N~~~~\nonumber
\end{equation}
where
\begin{equation}
\begin{array}{c}
 \omega_0=\omega= \sqrt{\lambda^*}, \\
  \omega _{2k-1}=\frac{\omega}{2} \sqrt{\sqrt{9 \iota_k^2-34 \iota_k+25}-\iota_k+5} , ~~~~~~for~~ k \in \{3,\cdots,N\},\\
  \omega _{2k}=\frac{\omega}{2} \sqrt{\sqrt{9 \iota_k^2-34 \iota_k+25}+\iota_k-5},~~~~~~for~~ k \in \{3,\cdots,N\}.
\end{array}\nonumber
\end{equation}
We will prove that $\omega _0,\omega _{2k-1}$ ($3 \leq k\leq N$) are nonresonant up to order $4$ generically. Obviously, we need only consider  the resonance of the frequencies $\mu_k$ ($1 \leq k\leq N-1$), here
\begin{equation}\label{resonance}
\begin{array}{c}
                                                                                                          \mu _1=1, \\
                                                                                                          \mu _{k}=\frac{1}{2} \sqrt{\sqrt{9 \iota_{k+1}^2-34 \iota_{k+1}+25}-\iota_{k+1}+5}, ~~~~~~~~2 \leq k\leq N-1.
                                                                                                        \end{array}
  \nonumber
\end{equation}

By virtue of (\ref{n-body1}), we have
\begin{equation}\label{muk}
\mu _{k}=a_k+b_k \varepsilon_1^ {\frac{1}{3}}+ O(\varepsilon_1^{\frac{2}{3}}), ~~~~~~~~2 \leq k\leq N-1,
\end{equation}
where
\begin{displaymath}
\begin{array}{c}
  a_k= \frac{1}{2} \sqrt{\sqrt{3^{2 k+2}-34\cdot 3^k+25}-3^k+5},\\
  b_k= \frac{ 3^{k-1} \left(\frac{17-3^{k+2}}{\sqrt{9^{k+1}-34\cdot 3^k+25}}+1\right) c_2 }{\sqrt{\sqrt{9^{k+1}-34\cdot 3^k+25}-3^k+5}}.
\end{array}
\end{displaymath}

It's easy to see that $a_k$ is monotonically increasing with respect to $k$. As a matter of fact, a straight forward computation shows that $a_{k+1}/a_k$ is monotonically decreasing with respect to $k\geq 2$ and for $k\geq 2$ we have
\begin{equation}\label{akak1}
\sqrt{3}< a_{k+1}/a_k \leq \sqrt{\frac{\sqrt{1417}-11}{4 \sqrt{7}-2}}\approx 1.76186.
\end{equation}

Set $a_1=\mu _1=1$ and note that $a_2=\sqrt{2 \sqrt{7}-1}\approx 2.07159$. Let's now prove that $a _{k}$ ($1 \leq k\leq N-1$) are nonresonant up to order $4$.

To prove that $a _{k}$ ($1 \leq k\leq N-1$) are nonresonant of order $3$, it suffices to show that both of the following statements are impossible:
\begin{description}
  \item[1)] there exist two elements $a_{k_1},a_{k_2}$ ($1\leq k_1<k_2\leq N-1$) such that $a_{k_2}=2a_{k_1}$;
  \item[2)] there exist three elements $a_{k_1},a_{k_2},a_{k_3}$ ($1\leq k_1<k_2<k_3 \leq N-1$) such that $a_{k_3}=a_{k_2}+a_{k_1}$.
\end{description}
By virtue of (\ref{akak1}) and $a_{k+1}/a_k\approx 2.07159$,  it can easily be seen  that the case \textbf{1)} is impossible. For the case \textbf{1)}, since we have the inequality
\begin{displaymath}
a_{k_2}+a_{k_1}< (\frac{1}{3^{\frac{k_3-k_2}{2}}}+\frac{1}{3^{\frac{k_3-k_1}{2}}})a_{k_3}\leq (\frac{1}{3^{\frac{1}{2}}}+\frac{1}{3^{\frac{2}{2}}})a_{k_3}<a_{k_3},
\end{displaymath}
we know that the case \textbf{2)} is also impossible. Therefore $a _{k}$ ($1 \leq k\leq N-1$) are nonresonant of order $3$.

To prove that $a _{k}$ ($1 \leq k\leq N-1$) are nonresonant of order $4$, it suffices to show that all of the following statements are impossible:
\begin{description}
  \item[3)] there exist two elements $a_{k_1},a_{k_2}$ ($1\leq k_1<k_2\leq N-1$) such that $a_{k_2}=3a_{k_1}$;
  \item[4)] there exist three elements $a_{k_1},a_{k_2},a_{k_3}$ ($1\leq k_1<k_2<k_3 \leq N-1$) such that $a_{k_3}= a_{k_2}+ 2a_{k_1}$ or $a_{k_3}=2a_{k_2}\pm a_{k_1}$;
\item[5)] there exist four elements $a_{k_1},a_{k_2},a_{k_3},a_{k_4}$ ($1\leq k_1<k_2<k_3 <k_4 \leq N-1$) such that $a_{k_4}= a_{k_3}+ a_{k_2}+ a_{k_1}$.
\end{description}

First, it can easily be checked  that the case \textbf{3)} is impossible.

Next, for the case \textbf{4)}, if $k_2\leq k_3-2$, then  we have the inequalities
\begin{displaymath}
\begin{array}{c}
  a_{k_2}+ 2a_{k_1}<  (\frac{1}{3^{\frac{1}{2}}}+\frac{2}{3^{\frac{4}{2}}})a_{k_3}<a_{k_3}, \\
 2a_{k_2}\pm a_{k_1}\leq 2a_{k_2}+ a_{k_1}<  (\frac{2}{3^{\frac{2}{2}}}+\frac{1}{3^{\frac{3}{2}}})a_{k_3}<a_{k_3};
\end{array}
\end{displaymath}
if $k_2= k_3-1, k_1= k_2-1$, then it can easily be checked  that we have the inequalities
\begin{displaymath}
 \begin{array}{c}
   a_{k_2}+ 2a_{k_1} >  a_{k_3}, \\
   2a_{k_2}+ a_{k_1}> a_{k_3},\\
2a_{k_2}- a_{k_1}< a_{k_3};
 \end{array}
\end{displaymath}
 if $k_2= k_3-1, k_1\leq k_2-2$, then we have the inequality
\begin{displaymath}
 a_{k_2}+ 2a_{k_1}<  (\frac{1}{3^{\frac{2}{2}}}+\frac{2}{3^{\frac{3}{2}}})a_{k_3}<a_{k_3};
\end{displaymath}
if $k_2= k_3-1, k_1= k_2-2$, then it can easily be checked  that we have the inequalities
\begin{displaymath}
 \begin{array}{c}
   2a_{k_2}+ a_{k_1}> a_{k_3},\\
2a_{k_2}- a_{k_1}< a_{k_3};
 \end{array}
\end{displaymath}
if $k_2= k_3-1, k_1\leq k_2-3$, then we have the inequality
\begin{displaymath}
 2a_{k_2}\pm a_{k_1} \geq 2a_{k_2}- a_{k_1}> (\frac{2}{1.762}-\frac{1}{3^{\frac{k_3-k_1}{2}}})a_{k_3}\geq (\frac{2}{1.762}-\frac{1}{3^{\frac{4}{2}}})a_{k_3}>a_{k_3};
\end{displaymath}
to summarize, the case \textbf{4)} is  impossible.

For the case \textbf{5)}, if $k_3\leq k_4-2$ or $k_2\leq k_3-2$ or $k_1\leq k_2-3$, then we have the inequality
\begin{displaymath}
a_{k_4}> a_{k_3}+ a_{k_2}+ a_{k_1};
\end{displaymath}
if $k_3= k_4-1, k_2= k_3-1,k_1= k_2-1$, then it can easily be checked  that we have the inequality
\begin{displaymath}
a_{k_4}< a_{k_3}+ a_{k_2}+ a_{k_1};
\end{displaymath}
 if $k_3= k_4-1, k_2= k_3-1,k_1= k_2-2$, then it can easily be checked  that we have the inequality
\begin{displaymath}
a_{k_4}\neq a_{k_3}+ a_{k_2}+ a_{k_1};
\end{displaymath}
 the case \textbf{5)} is  impossible.

Thus we arrive at the conclusion that $a _{k}$ ($1 \leq k\leq N-1$) are nonresonant up to order $4$. As a result, the frequencies $\omega _0,\omega _{2k-1}$ ($3 \leq k\leq N$) are nonresonant up to order $4$ for sufficiently small $\varepsilon_1$. According to the fact that $\omega _0,\omega _{2k-1}$ ($3 \leq k\leq N$) algebraically depend on the masses on the whole, it follows that $\omega _0,\omega _{2k-1}$ ($3 \leq k\leq N$) are nonresonant up to order $4$ generically.
Therefore, we could obtain the Birkhoff normal form of degree 4 for all the masses of the $N$-body Problem except at most  a proper algebraic subset of the mass space.

\subsection{The Degeneracy of The  Reduced Hamiltonian}
\indent\par
In this subsection, we will prove that the determinant of $\omega_{jk}$ on the center manifold is not zero for appropriate masses by the perturbed method and the inductive method. Then it follows that the reduced Hamiltonian  on the center manifold is nondegenerate for all the masses of the $N$-body problem except at most  a proper algebraic subset of the mass space.

Suppose that the reduced Hamiltonian  on the center manifold is nondegenerate for  the masses $m_1,\cdots,m_{N-1}$ of the $N-1$-body problem, by mathematical methods of induction, we will prove that the reduced Hamiltonian  on the center manifold is nondegenerate for  the masses $m_1,\cdots,m_N,$ of the $N$-body problem such that $m_{N}=\varepsilon$ is sufficiently small.

Note that if we consider the restricted $N$-body problem, that is, when $\varepsilon = 0$, then the definitions of $q_{jk}=\langle\widehat{\mathcal{E}}_j,\widehat{\mathcal{E}}^\bot_k\rangle$,  $a_{ijk}=d^3 \mathcal{U}|_{\widehat{\mathcal{E}}_3}(\widehat{\mathcal{E}}_i,\widehat{\mathcal{E}}_j,\widehat{\mathcal{E}}_k)$ and $a_{hijk}=d^4 \mathcal{U}|_{\widehat{\mathcal{E}}_3}(\widehat{\mathcal{E}}_h,\widehat{\mathcal{E}}_i,\widehat{\mathcal{E}}_j,\widehat{\mathcal{E}}_k)$ may be singular. This is because that
$g_{2N-1},g_{2N} \sim \varepsilon $, and $\|\mathcal{E}_{2N-1}\|,\|\mathcal{E}_{2N}\| \sim \sqrt{\varepsilon} \approx 0$, as a result,
\begin{displaymath}
  \widehat{\mathcal{E}}_{2N-1}= (0, 0,  \cdots, \frac{1}{\sqrt{\varepsilon}}, 0,)^\top+O(\sqrt{\varepsilon}).
\end{displaymath}

For convenience, set
\begin{displaymath}
\mathbf{e}_{k}(\varepsilon)= \|\mathbf{e}_{k}(\varepsilon)\| (e_{k,1}(\varepsilon),  e_{k,2}(\varepsilon),  \cdots, e_{k,N}(\varepsilon))^\top, ~~~~~~~~   k\in \{2,4,\cdots,N\},
\end{displaymath}
here $\mathbf{e}_{k}(\varepsilon)$ ($k=2,4,\cdots,N$) can be treated as the vectors  $\mathbf{e}_{k} (m_2,\cdots,m_{N-1},\varepsilon)$ in the \textbf{subsection 6.2}.

It's clear that
\begin{displaymath}
\begin{array}{c}
  \widehat{\mathcal{E}}_{2k-1}(\varepsilon)= (e_{k,1}(\varepsilon), 0, e_{k,2}(\varepsilon), 0, \cdots, e_{k,N}(\varepsilon), 0,)^\top, ~~~~~~for~~   k\in \{2,\cdots,N\}, \\
  \widehat{\mathcal{E}}_{2k}(\varepsilon)=\widehat{\mathcal{E}}_{2k-1}^ \bot=(0,e_{k,1}(\varepsilon), 0, e_{k,2}(\varepsilon), 0, \cdots,0, e_{k,N}(\varepsilon))^\top, ~~~~~~for~~   k\in \{2,\cdots,N\},
\end{array}
\end{displaymath}
where the central configuration $\widehat{\mathcal{E}}_3(\varepsilon) = (e_{2,1}(\varepsilon), 0, e_{2,2}(\varepsilon), 0, \cdots, e_{2,N}(\varepsilon), 0,)^\top$ satisfies the condition
\begin{displaymath}
e_{2,1}(\varepsilon)< e_{2,2}(\varepsilon)< \cdots <e_{2,N}(\varepsilon).
\end{displaymath}
Furthermore, it's easy to see that
\begin{displaymath}
\begin{array}{c}
  e_{k,j}(\varepsilon)= e_{k,j}(0)+O(\varepsilon), ~~~~~~for~~  j\in \{1,\cdots,N\}~~and~~ k\in \{2,\cdots,N-1\}, \\
  e_{N,j}(\varepsilon)= O(\sqrt{\varepsilon}), ~~~~~~for~~  j\in \{1,\cdots,N-1\}, \\
  e_{N,N}(\varepsilon)=\frac{1}{\sqrt{\varepsilon}}+O(\sqrt{\varepsilon}).
\end{array}
\end{displaymath}
Set $e_{k,ij}(\varepsilon)=e_{k,j}(\varepsilon)-e_{k,i}(\varepsilon)$, then
\begin{displaymath}
\begin{array}{c}
   e_{k,ij}(\varepsilon)=e_{k,j}(\varepsilon)-e_{k,i}(\varepsilon)= e_{k,ij}(0)+O(\varepsilon), ~~~~~~for~~   k\in \{2,\cdots,N-1\}, \\
   e_{N,ij}(\varepsilon)=O(\sqrt{\varepsilon}), ~~~~~~for~~  i,j\in \{1,\cdots,N-1\},  \\
   e_{N,jN}(\varepsilon)=\frac{1}{\sqrt{\varepsilon}}+O(\sqrt{\varepsilon}), ~~~~~~for~~   j\in \{1,\cdots,N-1\}.
\end{array}
\end{displaymath}
It's noteworthy that $|e_{2,ij}(0)|>0$ for $i\neq j$; the vector $\mathcal{E}_{3}(0)$ (or $\mathbf{e}_{2}(0)$) is  the central configuration of the restricted $N$-body problem; the vectors $\widehat{\mathcal{E}}_{k}(0)$ ($3\leq k\leq 2N-2$) and
\begin{displaymath}
{\mathcal{E}}_{2N-1}(0)=(0, 0,  \cdots,0, 0, 1, 0,)^\top,{\mathcal{E}}_{2N}(0)={\mathcal{E}}^\bot_{2N-1}(0)
\end{displaymath} are just
eigenvectors of  the restricted $N$-body problem; especially the vectors
\begin{displaymath}
(e_{k,1}(0), e_{k,2}(0), \cdots, e_{k,N-1}(0))^\top, ~~~~~~~~~~~~k\in \{2,\cdots,N-1\}
\end{displaymath}
are just
eigenvectors of  the collinear $N-1$-body problem.

Thanks to (\ref{lambdalambdak}) and (\ref{n-body1}),  the corresponding eigenvalues $\lambda_k^*(\varepsilon)$ and $\lambda^*(\varepsilon)$ satisfy
\begin{displaymath}
\begin{array}{c}
  \lambda^*(\varepsilon) = \|\mathbf{e}_2(\varepsilon)\|^3=\lambda^*(0)+O(\varepsilon),\\
  \lambda_k^*(\varepsilon)=\lambda_k^*(0)+O(\varepsilon), ~~~~~~for~~   k\in \{5,4,\cdots,2N\}.
\end{array}
\end{displaymath}
where the central configuration $\mathbf{e}_2(\varepsilon)=(\xi_1, \xi_2,\cdots, \xi_{N})^\top$ is the unique solution of  the equations  of central configurations:
\begin{equation}
\sum_{j=1,j\neq k}^N \frac{m_j}{|\mathbf{r}_j-\mathbf{r}_k|^3}(\mathbf{r}_j-\mathbf{r}_k)=- \mathbf{r}_k,1\leq k\leq N,\nonumber
\end{equation}
such that $\xi_1<\xi_2<\cdots<\xi_N$ and $\lambda(\mathbf{e}_2)=1$; and the  eigenvalues $\lambda_k^*(0)$ and $\lambda^*(0)$ correspond the restricted $N$-body problem, especially, they are the same as the corresponding values of the $N-1$-body problem.

Recall that
\begin{equation}
\begin{aligned}
&U(x) = \mathcal{U}(x_3 \widehat{\mathcal{E}}_3+\sum_{k = 5}^{2N} {x_k \widehat{\mathcal{E}}_k}) = \sum_{1\leq i<j\leq N} {\frac{m_i m_j }{r_{ij}}}\\
&= \lambda^* + \sum_{k = 5}^{2N} \frac{{\lambda^*_k}}{2} x^2_k +\sum_{i,j,k = 5}^{2N} \frac{a_{ijk}}{6} x_i x_j x_k +\frac{3}{4} \sum_{k = 5}^{2N}  x^2_k\sum_{k = 5}^{2N} ({\lambda^*_k}- \frac{\lambda^*}{2})  x^2_k+\sum_{h,i,j,k = 5}^{2N} \frac{a_{hijk}}{24} x_h x_i x_j x_k+ \cdots,
\end{aligned}\nonumber
\end{equation}
here
\begin{equation}\label{x3}
x_3 = \sqrt{1 - \sum_{k = 5}^{2N}  x^2_k},\nonumber
\end{equation}
\begin{equation}\label{rij}
 r_{ij}= \sqrt{( e_{2,ij}x_3+\sum_{k = 3}^{N} e_{k,ij}x_{2k-1})^2+(\sum_{k = 3}^{N} e_{k,ij}x_{2k})^2},\nonumber
\end{equation}
and
\begin{displaymath}
\begin{array}{c}
  \lambda^* = \lambda^*(\varepsilon) , ~~~~~~~~~~~~\lambda_k^*=\lambda_k^*(\varepsilon),\\
  a_{ijk} = a_{ijk}(\varepsilon) , ~~~~~~~~~~~~a_{hijk}=a_{hijk}(\varepsilon).
\end{array}
\end{displaymath}

Due to \begin{displaymath}
\begin{aligned}
&{\frac{m_i m_j }{r_{ij}}}={\frac{m_i m_j }{|e_{2,ij}|}}\{1-\sum_{k = 3}^{N} \frac{ e_{k,ij} x_{2k-1}}{e_{2,ij}}+\frac{1}{2}\sum_{k = 5}^{2N}x_{k}^2+\frac{(\sum_{k = 3}^{N} e_{k,ij}x_{2k-1})^2}{e_{2,ij}^2}-\frac{(\sum_{k = 3}^{N} e_{k,ij}x_{2k})^2}{2e_{2,ij}^2}\\
&-\sum_{k = 5}^{2N}x_{k}^2 \sum_{k = 3}^{N} \frac{ e_{k,ij} x_{2k-1}}{e_{2,ij}}-(\sum_{k = 3}^{N} \frac{ e_{k,ij} x_{2k-1}}{e_{2,ij}})^3 + \frac{3(\sum_{k = 3}^{N} e_{k,ij}x_{2k-1})(\sum_{k = 3}^{N} e_{k,ij}x_{2k})^2}{2 e_{2,ij}^3} \\
&+(\frac{3}{8})(\sum_{k = 5}^{2N}x_{k}^2)^2+ \frac{3}{2}(\sum_{k = 3}^{N} \frac{ e_{k,ij} x_{2k-1}}{e_{2,ij}})^2(\sum_{k = 5}^{2N}x_{k}^2)-\frac{3}{4}(\sum_{k = 3}^{N} \frac{ e_{k,ij} x_{2k}}{e_{2,ij}})^2(\sum_{k = 5}^{2N}x_{k}^2) +(\sum_{k = 3}^{N} \frac{ e_{k,ij} x_{2k-1}}{e_{2,ij}})^4\\
&+\frac{3}{8}(\sum_{k = 3}^{N} \frac{ e_{k,ij} x_{2k}}{e_{2,ij}})^4
-3(\sum_{k = 3}^{N} \frac{ e_{k,ij} x_{2k-1}}{e_{2,ij}})^2(\sum_{k = 3}^{N} \frac{ e_{k,ij} x_{2k}}{e_{2,ij}})^2 \}+O(\parallel x\parallel^5),
\end{aligned}
\end{displaymath}
it follows that
\begin{equation}
\begin{aligned}
&U(x) = \mathring{U}+ \frac{1}{2} (\mathring{\lambda}_{2N-1} x^2_{2N-1} +\mathring{\lambda}_{2N} x^2_{2N} )- \frac{\sigma_1}{\sqrt{\varepsilon}} x_{2N-1}^3+ \frac{3\sigma_1}{2\sqrt{\varepsilon}} x_{2N-1}x_{2N}^2\\
&+\frac{3}{4}\{(x^2_{2N-1} + x^2_{2N}) \sum_{k = 5}^{2N} (\mathring{\lambda}_k - \frac{\mathring{\lambda}}{2})  x^2_k+[(\mathring{\lambda}_{2N-1} - \frac{\mathring{\lambda}}{2})  x^2_{2N-1}+(\mathring{\lambda}_{2N} - \frac{\mathring{\lambda}}{2})  x^2_{2N}]\sum_{k = 5}^{2N-2}  x^2_k\}\\
&+ \frac{\sigma_2}{{\varepsilon}} x_{2N-1}^4+ \frac{\frac{3}{8}\sigma_2}{{\varepsilon}} x_{2N}^4- \frac{3\sigma_2}{{\varepsilon}} x_{2N-1}^2 x_{2N}^2+\mathcal{R}_3+\mathcal{R}_4+O_{\sqrt{\varepsilon}}+O(\parallel x\parallel^5),
\end{aligned}\nonumber
\end{equation}
where $O_{\sqrt{\varepsilon}}$ denotes the terms that whose coefficients are bounded with respect to $\sqrt{\varepsilon}$; and
\begin{equation}
\begin{aligned}
&\mathring{U} =  \mathring{\lambda} + \sum_{k = 5}^{2N}  \frac{\mathring{\lambda}_k}{2} x^2_k +\sum_{i,j,k = 5}^{2N-2}  \frac{\mathring{a}_{ijk}}{6} x_i x_j x_k +\frac{3}{4}\sum_{k = 5}^{2N-2}  x^2_k \sum_{k = 5}^{2N-2} (\mathring{\lambda}_k - \frac{\mathring{\lambda}}{2})  x^2_k+\sum_{h,i,j,k = 5}^{2N-2} \frac{\mathring{a}_{hijk}}{24} x_h x_i x_j x_k
\end{aligned}\nonumber
\end{equation}
is the same as the terms of the function $U(x)$   of the $N-1$-body Problem up to degree 4;
\begin{equation}
\begin{aligned}
&\mathcal{R}_3 =\sum_{j< N} {\frac{3 m_j\sum_{k = 3}^{N-1} e_{k,jN}(0) (\frac{1}{2 }x_{2k-1} x_{2N}^2-x_{2k-1} x^2_{2N-1}+ x_{2k} x_{2N-1}x_{2N})
 }{ e_{2,jN}^4(0)}},
\end{aligned}\nonumber
\end{equation}
\begin{equation}
\begin{aligned}
&\mathcal{R}_4 =\sum_{j< N} {\frac{ m_j }{ e_{2,jN}^5(0)}}[\frac{\sum_{k = 3}^{N-1} e_{k,jN}(0)\left( 4 x_{2N-1}^3 x_{2k-1}+\frac{3}{2}x_{2N}^3 x_{2k}- 6x_{2N-1} x_{2N}^2 x_{2k-1}-6x_{2N-1}^2 x_{2N} x_{2k}\right)}{\sqrt{\varepsilon}} \\
&+\frac{\left(\sqrt{\varepsilon} e_{N,jN}(\varepsilon)\right)^4-1}{\varepsilon} (x_{2N-1}^4+\frac{3}{8}x_{2N}^4-3x_{2N-1}^2 x_{2N}^2)+6 x_{2N-1}^2\left(\sum_{k = 3}^{N-1} e_{k,jN}(0) x_{2k-1}\right)^2\\
&-12x_{2N-1} x_{2N}\sum_{k = 3}^{N-1} e_{k,jN}(0) x_{2k-1}\sum_{k = 3}^{N-1} e_{k,jN}(0) x_{2k}+\frac{9}{4} x_{2N}^2\left(\sum_{k = 3}^{N-1} e_{k,jN}(0) x_{2k}\right)^2 ];
\end{aligned}\nonumber
\end{equation}
\begin{displaymath}
\begin{array}{c}
  \mathring{\lambda} = \lambda^*(0) , ~~~~~~~~~~~~\mathring{\lambda}_k = \lambda^*_{k}(0), \\
  \mathring{a}_{ijk}=a_{ijk}(0), ~~~~~~~~~~~~\mathring{a}_{hijk}=a_{hijk}(0);
\end{array}
\end{displaymath}
\begin{displaymath}
\sigma_1= \sum_{j=1}^{N-1}\frac{m_j}{e_{2,jN}(0)^4}>0,~~~~~~~~~~~~\sigma_2= \sum_{j=1}^{N-1} \frac{m_j}{e_{2,jN}(0)^5}>0;
\end{displaymath}
it's noteworthy that $\frac{\left(\sqrt{\varepsilon} e_{N,jN}(\varepsilon)\right)^4-1}{\varepsilon}$ is bounded with respect to $\varepsilon$.

Then the Hamiltonian $H$ becomes
\begin{equation}\label{HamiltonianEulerreduce}
H(x_0,x,y_0,y)= \frac{-{\omega}^2}{2}+ H_2(x_0,x,y_0,y)+H_3(x_0,x,y_0,y)+H_4(x_0,x,y_0,y)+\cdots,
\end{equation}
where
\begin{equation}\label{HamiltonianEulerreduceH2}
\begin{aligned}
&H_2=\frac{y_0^2}{2}+ \frac{1}{2} [ {\omega}^2 x_0^2+\sum_{k=5}^{2N}y_k^2+\sum_{k = 3}^{N}2{\omega} ({x}_{2k} y_{2k-1}-{x}_{2k-1} y_{2k})+\sum_{k=5}^{2N}( {\omega}^2-\lambda^*_k ){x}^2_k]=\\
& \mathring{H}_2+  \frac{ y_{2N-1}^2+y_{2N}^2+2 \mathring{\omega} ({x}_{2N} y_{2N-1}-{x}_{2N-1} y_{2N})+{\mathring{\omega}}^2( 1-\mathring{\iota}_N ) {x}^2_{2N-1}+{\mathring{\omega}}^2(\frac{\mathring{\iota}_N-1}{2} ){x}^2_{2N}}{2}  +O_{\varepsilon},\nonumber\\
\end{aligned}
\end{equation}
and
\begin{equation}\label{HamiltonianEulerreduceH34}
\begin{array}{lr}
\begin{aligned}
&H_3= \mathring{H}_3-\{x_0 [y_{2N-1}^2+y_{2N}^2+2\mathring{\omega} ({x}_{2N} y_{2N-1}-{x}_{2N-1} y_{2N})\\
&+{\mathring{\omega}}^2( 1-\frac{\mathring{\iota}_N}{2} ){x}^2_{2N-1}+{\mathring{\omega}}^2( \frac{\mathring{\iota}_N-1}{4} ){x}^2_{2N}]
- \frac{\sigma_1}{\sqrt{\varepsilon}} x_{2N-1}^3+ \frac{3\sigma_1}{2\sqrt{\varepsilon}} x_{2N-1}x_{2N}^2+\mathcal{R}_3\} +O_{\sqrt{\varepsilon}},\nonumber
\end{aligned}\\
\begin{aligned}
&H_4= \mathring{H}_4+\{x_0[- \frac{\sigma_1}{\sqrt{\varepsilon}} x_{2N-1}^3+ \frac{3\sigma_1}{2\sqrt{\varepsilon}} x_{2N-1}x_{2N}^2+\mathcal{R}_3]\\
&+\frac{3x_0^2}{2 }[ y_{2N-1}^2+y_{2N}^2+2 \mathring{\omega} ({x}_{2N} y_{2N-1}-{x}_{2N-1} y_{2N})+{\mathring{\omega}}^2( 1-\frac{\mathring{\iota}_N}{3} ){x}^2_{2N-1}+{\mathring{\omega}}^2(\frac{\mathring{\iota}_N+3}{6}){x}^2_{2N}]\\
&+\frac{1}{2 }[ (\sum_{j,k=5}^{2N}q_{kj}{x}_j y_k-\mathring{\omega}\sum_{j=5}^{2N}{x}^2_j)^2-(\sum_{j,k=5}^{2N-2}q_{kj}{x}_j y_k-\mathring{\omega}\sum_{j=5}^{2N-2}{x}^2_j)^2- (\sum_{k=5}^{2N}x_k y_k)^2+(\sum_{k=5}^{2N-2}x_k y_k)^2]\\
&-[\frac{3}{4}\left((x^2_{2N-1} + x^2_{2N}) \sum_{k = 5}^{2N} (\mathring{\lambda}_k - \frac{\mathring{\lambda}}{2})  x^2_k+\left((\mathring{\lambda}_{2N-1} - \frac{\mathring{\lambda}}{2})  x^2_{2N-1}+(\mathring{\lambda}_{2N} - \frac{\mathring{\lambda}}{2})  x^2_{2N}\right)\sum_{k = 5}^{2N-2}  x^2_k\right)\\
&+ \frac{\sigma_2}{{\varepsilon}} x_{2N-1}^4+ \frac{\frac{3}{8}\sigma_2}{{\varepsilon}} x_{2N}^4- \frac{3\sigma_2}{{\varepsilon}} x_{2N-1}^2 x_{2N}^2+\mathcal{R}_4]\} +O_{\sqrt{\varepsilon}};
\end{aligned}
\end{array}
\end{equation}
\begin{equation}\label{Hamiltonian3EulerrestrictH234}
\begin{array}{lr}
\mathring{H}_2=\frac{y_0^2}{2}+ \frac{1}{2} [ {\mathring{\omega}}^2 x_0^2+\sum_{k=5}^{2N-2}y_k^2+\sum_{k = 3}^{N-1}2{\mathring{\omega}} ({x}_{2k} y_{2k-1}-{x}_{2k-1} y_{2k})+\sum_{k=5}^{2N-2}( {\mathring{\omega}}^2-\mathring{\lambda}_k ){x}^2_k],\\
\mathring{H}_3=-[x_0 ({\mathring{\omega}}^2 x_0^2+\sum_{k=5}^{2N-2}y_k^2-2{\mathring{\omega}}\sum_{j,k=5}^{2N-2}q_{kj} {x}_j y_k+\sum_{k=5}^{2N-2}( {\mathring{\omega}}^2-\frac{\mathring{\lambda}_k}{2} ){x}^2_k)+\sum_{i,j,k = 5}^{2N-2} \frac{\mathring{a}_{ijk}}{6} x_i x_j x_k ],\\
\begin{aligned}
&\mathring{H}_4=\frac{1}{2} [(\sum_{j,k=5}^{2N-2}q_{kj}{x}_j y_k-{\mathring{\omega}}\sum_{k=5}^{2N-2}{x}^2_k)^2- (\sum_{k=5}^{2N-2}x_k y_k)^2]+ \frac{x_0}{6}\sum_{i,j,k = 5}^{2N-2} \mathring{a}_{ijk} x_i x_j x_k\nonumber\\
&+ \frac{3x_0^2}{2} [{\mathring{\omega}}^2 x_0^2+\sum_{k=5}^{2N-2}y_k^2-2{\mathring{\omega}}\sum_{j,k=5}^{2N-2}q_{kj} {x}_j y_k+\sum_{k=5}^{2N-2}( {\mathring{\omega}}^2-\frac{\mathring{\lambda}_k}{3} ){x}^2_k]\nonumber\\
&-  [ \frac{3 {\mathring{\omega}}^2}{8} (\sum_{k = 5}^{2N-2}  x^2_k)^2+\frac{3}{4} (\sum_{j = 5}^{2N-2}  x^2_j)\sum_{k = 5}^{2N-2} (\mathring{\lambda}_k - {\mathring{\omega}}^2)  x^2_k]-\frac{1}{24}\sum_{h,i,j,k = 5}^{2N-2} \mathring{a}_{hijk} x_h x_i x_j x_k;
\end{aligned}
\end{array}
\end{equation}
and
\begin{eqnarray}\label{Hamiltonian3Eulerrestrict}
\mathring{H}&=&\frac{-{\mathring{\omega}}^2}{2} +\mathring{H}_2+\mathring{H}_3+\mathring{H}_4\nonumber
\end{eqnarray}
is the same as the terms of the Hamiltonian   of the $N-1$-body Problem up to degree 4;
\begin{displaymath}
\omega = \sqrt{\lambda^*(\varepsilon)} , ~~~~~~~~~~~~\mathring{\omega} = \sqrt{\lambda^*(0)} , ~~~~~~~~~~~~\mathring{\iota}_N=\iota_N(0).
\end{displaymath}

According to the result in  the \textbf{subsection 6.3}, it's easy to see that the frequencies of the problem   are nonresonant up to order $4$ for sufficiently small $\varepsilon$. Therefore, the Hamiltonian $H$, $\mathring{H}$ and $H_{sim}$ could be reduced to a Birkhoff normal form of degree 4. Here the Hamiltonian  $H_{sim}$ is a simplification of the Hamiltonian $H$:
\begin{equation}\label{Hamiltonian3Eulersim}
H_{sim}(x_0,x,y_0,y) =H_{sim2}+H_{sim2}+H_{sim4}
\end{equation}
where
\begin{equation}\label{Hamiltonian3EulersimH234}
\begin{array}{lr}
H_{sim2}=\mathring{H}_2+  \frac{ y_{2N-1}^2+y_{2N}^2+2 \mathring{\omega} ({x}_{2N} y_{2N-1}-{x}_{2N-1} y_{2N})+{\mathring{\omega}}^2( 1-\mathring{\iota}_N ) {x}^2_{2N-1}+{\mathring{\omega}}^2(\frac{\mathring{\iota}_N-1}{2} ){x}^2_{2N}}{2},\\
H_{sim3}=\mathring{H}_3+ \frac{\sigma_1}{\sqrt{\varepsilon}} x_{2N-1}^3- \frac{3\sigma_1}{2\sqrt{\varepsilon}} x_{2N-1}x_{2N}^2,\\
H_{sim4}=\mathring{H}_4-(\frac{\sigma_2}{{\varepsilon}} x_{2N-1}^4+ \frac{\frac{3}{8}\sigma_2}{{\varepsilon}} x_{2N}^4- \frac{3\sigma_2}{{\varepsilon}} x_{2N-1}^2 x_{2N}^2).
\end{array}\nonumber
\end{equation}
As a matter of fact, we will show that the reduced Hamiltonian of $H$ on the center manifold is nondegenerate by comparing   the determinants $det(\omega_{jk})$ of the  Birkhoff normal forms for the two Hamiltonian $H$ and $H_{sim}$.

Before comparing  their determinants, let's recall the process of obtaining Birkhoff normal form of degree 4.

As a matter of notational convenience, set
\begin{displaymath}
\begin{array}{c}
  q_0=x_0,p_0=y_0, \\
  q_k=x_{k+4},p_k=y_{k+4},  ~~~~~~for~~ k \in \{1,\cdots,2N-4\}.
\end{array}
\end{displaymath}

The first step of  obtaining Birkhoff normal form is to simplify the
quadratic part $H_2$ of the Hamiltonian. It can easily be seen that  the simplification of the
quadratic part $H_2$  could be achieved by the following transformation
\begin{equation}\label{symplectic transformationeuler}
\left\{
             \begin{array}{lr}
             p_0  =  \sqrt{\omega_0}\mathfrak{p}_0 &  \\
             q_0  =  \frac{\mathfrak{q}_0}{\sqrt{\omega_0}} &  \\
             p_{2k-1}  = \frac{\omega_0   \left(2 \lambda^*_{2k+4} +\omega_{2k-1}^2-3 \lambda^*\right)\mathfrak{q}_{2k-1}}{\sqrt{r_{2k-1}} \left(-2 \lambda^*_{2k+4} +\omega_{2k-1}^2+3 \lambda^*\right)}+\frac{\omega_0   \left(-2 \lambda^*_{2k+4} +\omega_{2k}^2+3 \lambda^*\right)(\mathfrak{p}_{2k}+\mathfrak{q}_{2k})}{\sqrt{r_{2k}} \left(2 \lambda^*_{2k+4} +\omega_{2k}^2-3 \lambda^*\right)},&k =1,\cdots,N-2  \\
              p_{2k}  =  \frac{\omega_{2k-1}  \left(-2 \lambda^*_{2k+4} +\omega_{2k-1}^2+\lambda^*\right)\mathfrak{p}_{2k-1}}{\sqrt{r_{2k-1}} \left(-2 \lambda^*_{2k+4} +\omega_{2k-1}^2+3 \lambda^*\right)}+\frac{ \left(\frac{2 \omega_{2k} \lambda^*}{2 \lambda^*_{2k+4} +\omega_{2k}^2-3 \lambda^*}+\omega_{2k}\right)(-\mathfrak{p}_{2k}+\mathfrak{q}_{2k})}{\sqrt{r_{2k}}},&k =1,\cdots,N-2 \\
              q_{2k-1}  =  -\frac{2 \omega_{2k-1}  \omega_0 \mathfrak{p}_{2k-1}}{\sqrt{r_{2k-1}} \left(-2 \lambda^*_{2k+4} +\omega_{2k-1}^2+3 \lambda^*\right)}+\frac{2 \omega_{2k} \omega_0  (-\mathfrak{p}_{2k}+\mathfrak{q}_{2k})}{\sqrt{r_{2k}} \left(2 \lambda^*_{2k+4} +\omega_{2k}^2-3 \lambda^*\right)},& k =1,\cdots,N-2 \\
              q_{2k}  = \frac{\mathfrak{q}_{2k-1}}{\sqrt{r_{2k-1}}}+\frac{\mathfrak{p}_{2k}+\mathfrak{q}_{2k}}{\sqrt{r_{2k}}},&k =1,\cdots,N-2
             \end{array}
\right.
\end{equation}
here
\begin{equation}
\begin{array}{lr}
r_{2k-1}=\frac{\omega_{2k-1} \left(-\lambda^*_{2k+4} +2 \omega_{2k-1}^2-\lambda^*\right)}{-2 \lambda^*_{2k+4} +\omega_{2k-1}^2+3 \lambda^*}, ~~~~~~~~ k =1,\cdots,N-2,\\
r_{2k}= -\frac{2 \omega_{2k} \left(\lambda^*_{2k+4} +2 \omega_{2k}^2+\lambda^*\right)}{2 \lambda^*_{2k+4} +\omega_{2k}^2-3 \lambda^*}, ~~~~~~~~ k =1,\cdots,N-2;
\end{array}
\nonumber
\end{equation}
and
\begin{equation}
\begin{array}{c}
 \omega_0=\omega= \sqrt{\lambda^*}, \\
  \omega _{2k-1}=\frac{\omega}{2} \sqrt{\sqrt{9 \iota_{k+2}^2-34 \iota_{k+2}+25}-\iota_{k+2}+5} , ~~~~~~for~~ k \in \{1,\cdots,N-2\},\\
  \omega _{2k}=\frac{\omega}{2} \sqrt{\sqrt{9 \iota_{k+2}^2-34 \iota_{k+2}+25}+\iota_{k+2}-5},~~~~~~for~~ k \in \{1,\cdots,N-2\}.
\end{array}\nonumber
\end{equation}

Then the Hamiltonian $H$
becomes
\begin{equation}\label{firststep}
\begin{aligned}
&H(\mathfrak{p},\mathfrak{q}) =-\frac{\omega_0^2}{2}+\frac{\omega_0 (\mathfrak{p}_0^2+\mathfrak{q}_0^2)}{2}+ \sum_{k = 1}^{N-2}[\frac{\omega_{2k-1}(\mathfrak{p}_{2k-1}^2+\mathfrak{q}_{2k-1}^2)}{2}+\omega_{2k}\mathfrak{p}_{2k} \mathfrak{q}_{2k}]\\
&+ H_3(\mathfrak{p},\mathfrak{q})+H_4(\mathfrak{p},\mathfrak{q}) +\cdots,
\end{aligned}
\end{equation}
or by a further  transformation
\begin{equation}\label{fu symplectic transformationeuler}
\left\{
             \begin{array}{lr}
             \mathfrak{p} _0=\frac{\zeta_0}{\sqrt{2}}+\frac{\textbf{i} \eta_0}{\sqrt{2}} &  \\
             \mathfrak{q} _0=\frac{\eta_0}{\sqrt{2}}+\frac{\textbf{i} \zeta_0}{\sqrt{2}} &  \\
             \mathfrak{p} _{2k-1}=\frac{\zeta_{2k-1}}{\sqrt{2}}+\frac{\textbf{i} \eta_{2k-1}}{\sqrt{2}},&k =1,\cdots,N-2  \\
              \mathfrak{q} _{2k-1}=\frac{\eta_{2k-1}}{\sqrt{2}}+\frac{\textbf{i} \zeta_{2k-1}}{\sqrt{2}},&k =1,\cdots,N-2  \\
              \mathfrak{p} _{2k}=\zeta_{2k},&k =1,\cdots,N-2  \\
              \mathfrak{q} _{2k}=\eta_{2k},&k =1,\cdots,N-2
             \end{array}
\right.
\end{equation}
the Hamiltonian $H$
becomes
\begin{equation}
H(\zeta,\eta)  = -\frac{\omega_0^2}{2}+\textbf{i}\omega _0 \zeta _0 \eta _0 + \sum_{k = 1}^{N-2}[\textbf{i}\omega _{2k-1} \zeta _{2k-1} \eta _{2k-1} + \omega _{2k} \zeta _{2k} \eta _{2k}] + H_3(\zeta,\eta)+H_4(\zeta,\eta) +\cdots.
\nonumber
\end{equation}

By the way, it's obvious that the reduced Hamiltonian on  the center manifold is positive definite by (\ref{firststep}). Therefore, one can obtain  $N-1$ one parameter family of periodic orbits  that lie near Euler relative equilibrium by Theorem \ref{weinstein}.

The second step of  obtaining Birkhoff normal form is to simplify the
cubic part $H_3$ of the Hamiltonian. This simplification  could be achieved by a change of variables $(\zeta,\eta)\mapsto (u,v)$
with a generating function
\begin{equation}
 \sum_{k = 0}^{2N-4}u _k \eta _k + S_3(u,\eta),
\nonumber
\end{equation}
and $S_3$ can be  determined by the equation
\begin{equation}\label{connectrelation3euler}
\begin{aligned}
&\textbf{i} \omega _0 (\frac{\partial S_3}{\partial \eta_0} \eta _0-\frac{\partial S_3}{\partial u_0} u _0) + \sum_{k = 1}^{N-2}\textbf{i}\omega _{2k-1} (\frac{\partial S_3}{\partial \eta_{2k-1}} \eta _{2k-1}-\frac{\partial S_3}{\partial u_{2k-1}} u _{2k-1})\\
 &+ \sum_{k = 1}^{N-2}\omega _{2k} (\frac{\partial S_3}{\partial \eta_{2k}} \eta _{2k}-\frac{\partial S_3}{\partial u_{2k}} u _{2k}) + H_3(u,\eta)= 0.
\end{aligned}
\end{equation}

Similarly, the last step of  obtaining Birkhoff normal form is to simplify the
quartic part $H_4$ of the Hamiltonian, and this simplification  could be achieved by a change of variables
with a generating function
\begin{equation}
 \sum_{k = 0}^{2N-4}u _k \eta _k + S_3(u,\eta)+ S_4(u,\eta),
\nonumber
\end{equation}
$S_4$ can be  determined by the equation
\begin{equation}\label{connectrelation4euler}
\begin{aligned}
&\textbf{i} \omega _0 (\frac{\partial S_4}{\partial \eta_0} \eta _0-\frac{\partial S_4}{\partial u_0} u _0) + \sum_{k = 1}^{N-2}\textbf{i}\omega _{2k-1} (\frac{\partial S_4}{\partial \eta_{2k-1}} \eta _{2k-1}-\frac{\partial S_4}{\partial u_{2k-1}} u _{2k-1})\\
&+ \sum_{k = 1}^{N-2}\omega _{2k} (\frac{\partial S_4}{\partial \eta_{2k}} \eta _{2k}-\frac{\partial S_4}{\partial u_{2k}} u _{2k}) + H_{3\rightarrow 4}+H_4(u,\eta)+\frac{1}{2}[\omega_{00}(u_0 v_0)^2\\
&+\sum_{k = 1}^{N-2}\omega_{{2k-1},{2k-1}}(u_{2k-1} v_{2k-1})^2-\sum_{k = 1}^{N-2}\omega_{2k,2k}(u_{2k} v_{2k})^2+2\sum_{k = 1}^{N-2}\omega_{0,{2k-1}}u_0 v_0 u_{2k-1} v_{2k-1}\\
&+2\sum_{1\leq j<k \leq N-2}\omega_{{2j-1},{2k-1}}u_{2j-1} v_{2j-1} u_{2k-1} v_{2k-1}-2\sum_{k = 1}^{N-2}\textbf{i}\omega_{0,{2k}}(u_0 v_0 u_{2k} v_{2k})\\
&-2\sum_{j,k = 1}^{N-2}\textbf{i}\omega_{{2j-1},{2k}}(u_{2j-1} v_{2j-1} u_{2k} v_{2k})-2\sum_{1\leq j<k \leq N-2}\omega_{{2j},{2k}}u_{2j} v_{2j} u_{2k} v_{2k}]=0 ,
\end{aligned}
\end{equation}
where $H_{3\rightarrow 4}$ is the forms of  order 4 of $H_3(u+\frac{\partial S_3}{\partial \eta},\eta)$.

Then the  Birkhoff normal form of degree 4 for the  Hamiltonian of $H$ is
\begin{equation}\label{Birkhoffeuler}
\begin{aligned}
&\mathcal{H}(u,v)
=\textbf{i} \omega _0 u_0 v_0 + \sum_{k = 1}^{N-2}[\textbf{i}\omega _{2k-1} u _{2k-1} v _{2k-1} + \omega _{2k} u _{2k} v _{2k}]  -\frac{1}{2}[\omega_{00}(u_0 v_0)^2\\
&+\sum_{k = 1}^{N-2}\omega_{{2k-1},{2k-1}}(u_{2k-1} v_{2k-1})^2-\sum_{k = 1}^{N-2}\omega_{2k,2k}(u_{2k} v_{2k})^2+2\sum_{k = 1}^{N-2}\omega_{0,{2k-1}}u_0 v_0 u_{2k-1} v_{2k-1}\\
&+2\sum_{1\leq j<k \leq N-2}\omega_{{2j-1},{2k-1}}u_{2j-1} v_{2j-1} u_{2k-1} v_{2k-1}-2\sum_{k = 1}^{N-2}\textbf{i}\omega_{0,{2k}}(u_0 v_0 u_{2k} v_{2k})\\
&-2\sum_{j,k = 1}^{N-2}\textbf{i}\omega_{{2j-1},{2k}}(u_{2j-1} v_{2j-1} u_{2k} v_{2k})-2\sum_{1\leq j<k \leq N-2}\omega_{{2j},{2k}}u_{2j} v_{2j} u_{2k} v_{2k}] ,
\end{aligned}\nonumber
\end{equation}
and the determinant  of the  Birkhoff normal form for the reduced Hamiltonian of $H$ on the center manifold is the determinant of the matrix
\begin{displaymath}
\Omega=\begin{pmatrix}
           \omega_{00} & \omega_{01} & \omega_{03} & \cdots & \omega_{0,2N-5}  \\
 \omega_{01} & \omega_{11} &\omega_{13} & \cdots & \omega_{1,2N-5} \\
\omega_{03} & \omega_{13} & \omega_{33} &\cdots& \omega_{3,2N-5} \\
 \vdots & \vdots &\vdots & \ddots & \vdots  \\
 \omega_{0,2N-5} & \omega_{1,2N-5} & \omega_{3,2N-5} &\cdots & \omega_{2N-5,2N-5}
    \end{pmatrix}
\end{displaymath}

Our task is now to estimate the values of elements in the above matrix $\Omega$. We can complete this task by comparing  the elements  of  the  Hamiltonian $H$ with  that of  the  Hamiltonian $H_{sim}$.

First, thanks to the equation (\ref{connectrelation4euler}), the element $\omega_{jk}$ is determined only by the coefficients of the terms $u_j \eta_j u_{k} \eta_{k}$ in $H_{3\rightarrow 4}$ and $H_{4}$. So it is essential that  to investigate $H_{3}$, $S_{3}$ and $H_{4}$.

Let the transformation
\begin{equation}\label{symplecticcompositetransformation}
\left\{
             \begin{array}{lr}
             p _{k}=\phi_{k}(\zeta,\eta,\varepsilon),&k =0,1,\cdots,2N-4  \\
              q _{k}=\psi_{k}(\zeta,\eta,\varepsilon),&k =0,1,\cdots,2N-4
             \end{array}
\right.
\end{equation}
denote the composite transformation of the transformations (\ref{symplectic transformationeuler}) and (\ref{fu symplectic transformationeuler}), here $\varepsilon$ in the transformation just indicates the fact that the transformation depends on  $\varepsilon$. It can easily be verified that
\begin{equation}
\left\{
             \begin{array}{lr}
             \phi_{k}(\zeta,\eta,\varepsilon)= \phi_{k}(\zeta,\eta,0)+O(\varepsilon),&k =0,1,\cdots,2N-4  \\
             \psi_{k}(\zeta,\eta,\varepsilon)=\psi_{k}(\zeta,\eta,0)+O(\varepsilon),&k =0,1,\cdots,2N-4
             \end{array}\nonumber
\right.
\end{equation}
and
the transformation
\begin{equation}\label{symplecticcompositetransformation0}
\left\{
             \begin{array}{lr}
             p _{k}=\phi_{k}(\zeta,\eta,0),&k =0,1,\cdots,2N-4  \\
              q _{k}=\psi_{k}(\zeta,\eta,0),&k =0,1,\cdots,2N-4
             \end{array}
\right.
\end{equation}
is exactly the transformation  simplifying  the
quadratic part  $H_{sim2}$ of the Hamiltonian $H_{sim}$.

As a matter of notational convenience, in the following we simply write $\phi_{k}$ and $\psi_{k}$ to represent $\phi_{k}(\zeta,\eta,\varepsilon)$ and $\psi_{k}(\zeta,\eta,\varepsilon)$ respectively, similarly,  $\mathring{\phi}_{k}$ and $\mathring{\psi}_{k}$  represent $\phi_{k}(\zeta,\eta,0)$ and $\psi_{k}(\zeta,\eta,0)$ respectively.

A straight forward computation shows  that the parts of $H$ and $H_{sim}$ have the following relations
\begin{equation}\label{HamiltonianEulerreducerelationsH234}
H_2(\phi,\psi)  = \textbf{i}\omega _0 \zeta _0 \eta _0 + \sum_{k = 1}^{N-2}[\textbf{i}\omega _{2k-1} \zeta _{2k-1} \eta _{2k-1} + \omega _{2k} \zeta _{2k} \eta _{2k}] = H_{sim2}(\mathring{\phi},\mathring{\psi})  +O_{\varepsilon},\nonumber
\end{equation}
\begin{equation}
\begin{aligned}
&H_3(\phi,\psi)= H_{sim3}(\mathring{\phi},\mathring{\psi})-\{\mathring{\psi}_0 [\mathring{\phi}_{2N-5}^2+\mathring{\phi}_{2N-4}^2+2\mathring{\omega} (\mathring{\psi}_{2N-4} \mathring{\phi}_{2N-5}-\mathring{\psi}_{2N-5} \mathring{\phi}_{2N-4})\\
&+{\mathring{\omega}}^2( 1-\frac{\mathring{\iota}_N}{2} )\mathring{\psi}_{2N-5}^2+{\mathring{\omega}}^2( \frac{\mathring{\iota}_N-1}{4} )\mathring{\psi}_{2N-4}^2]
+\mathcal{R}_3(\mathring{\phi},\mathring{\psi})\} +O_{\sqrt{\varepsilon}}\\
&= H_{sim3}(\mathring{\phi},\mathring{\psi})+ O_1(\zeta _{2N-5},\zeta _{2N-4}, \eta _{2N-5},\eta _{2N-4})+O_{\sqrt{\varepsilon}},\nonumber
\end{aligned}
\end{equation}
\begin{equation}
\begin{aligned}
&H_4(\phi,\psi)= H_{sim4}(\mathring{\phi},\mathring{\psi})+ O_1(\zeta _{2N-5},\zeta _{2N-4}, \eta _{2N-5},\eta _{2N-4}) +O_{\sqrt{\varepsilon}}+\\
&\sum_{j< N} {\frac{ \sum_{k = 1}^{N-3} e_{k+2,jN}(0)\left( 4 \mathring{\psi}_{2N-5}^3 \mathring{\psi}_{2k-1}+\frac{3}{2}\mathring{\psi}_{2N-4}^3 \mathring{\psi}_{2k}- 6\mathring{\psi}_{2N-5} \mathring{\psi}_{2N-4}^2 \mathring{\psi}_{2k-1}-6\mathring{\psi}_{2N-5}^2 \mathring{\psi}_{2N-4} \mathring{\psi}_{2k}\right) }{ \sqrt{\varepsilon} e_{2,jN}^5(0)/m_j}},
\end{aligned}\nonumber
\end{equation}
where $O_1(\zeta _{2N-5},\zeta _{2N-4}, \eta _{2N-5},\eta _{2N-4})$ denotes the terms which contain at least two of $\zeta _{2N-5}$, $\zeta _{2N-4}$, $\eta _{2N-5}, \eta _{2N-4}$ as a factor and whose coefficients are bounded with respect to $\varepsilon$.

According to the equation (\ref{connectrelation3euler}), the function $S_3$ for the Hamiltonian $H$ has the form
\begin{equation}
\begin{aligned}
&S_3(u,\eta)
=S_{sim3}(u,\eta)+ O_1(u _{2N-5},u _{2N-4}, \eta _{2N-5},\eta _{2N-4}) +O_{\sqrt{\varepsilon}},
\end{aligned}\nonumber
\end{equation}
where the function $S_{sim3}$ corresponds to the Hamiltonian $H_{sim}$.

Consequently, the coefficients $c_{jk}$ of the terms $u_j \eta_j u_{k} \eta_{k}$ in $H_{3\rightarrow 4}$ (or  $H_3(u+\frac{\partial S_3}{\partial \eta},\eta)$) have the forms
\begin{equation}
\begin{aligned}
&c_{jk}
=c_{simjk} +O({\sqrt{\varepsilon}}), &j,k=0,1,\cdots,2N-6\\
&c_{jk}
=c_{simjk}+ O(1) ,&j=0,1,\cdots,2N-6,~~~k=2N-5,2N-4\\
&c_{jk}
=c_{simjk}+ O(\frac{1}{\sqrt{\varepsilon}}) ,&j,k=2N-5,2N-4
\end{aligned}\nonumber
\end{equation}
where $c_{simjk}$ are the coefficients of the terms $u_j\eta_j u_{k}\eta_{k}$ in  $H_{sim3}(u+\frac{\partial S_{sim3}}{\partial \eta},\eta)$. Similarly,
\begin{equation}
\begin{aligned}
&d_{jk}
=d_{simjk} +O({\sqrt{\varepsilon}}), &j,k=0,1,\cdots,2N-6\\
&d_{jk}
=d_{simjk}+ O(1) ,&j=0,1,\cdots,2N-6;~~~k=2N-5,2N-4\\
&d_{jk}
=d_{simjk}+ O(\frac{1}{\sqrt{\varepsilon}}) ,&j,k=2N-5,2N-4
\end{aligned}\nonumber
\end{equation}
where $d_{jk}$ and $d_{simjk}$ are the coefficients of the terms $u_j\eta_j u_{k}\eta_{k}$ in $H_4(u,\eta)$ and $H_{sim4}(u,\eta)$ respectively.

We claim that
\begin{equation}
\begin{aligned}
&c_{simjk}+d_{simjk}
=const=O(1), &j,k=0,1,\cdots,2N-6\\
&c_{simjk}+d_{simjk}
= 0 ,&j=0,1,\cdots,2N-6;~~~k=2N-5,2N-4\\
&c_{simjk}+d_{simjk}
\sim\frac{1}{\varepsilon} ,&j,k=2N-5,2N-4
\end{aligned}\nonumber
\end{equation}
and the determinants $det(\Omega)$ and $ det(\Omega_{sim})$ of the  Birkhoff normal forms for the Hamiltonian $H$ and $H_{sim}$ satisfy
\begin{equation}
\begin{aligned}
det(\Omega_{sim}) \sim \frac{1}{\varepsilon}, ~~~~~~~~~~~~~~~det(\Omega) = det(\Omega_{sim}) +O(\frac{1}{\sqrt{\varepsilon}}).
\end{aligned}\nonumber
\end{equation}

As a matter of fact, the claim will be prove if we can obtain the Birkhoff normal form of the  Hamiltonian  $H_{sim}$.

According to (\ref{Hamiltonian3Eulersim}), it's easy to see that we need only pay  attention to
the Hamiltonian
\begin{equation}
\begin{aligned}
&H_{\varepsilon} =  \frac{ 1}{2}[y_{2N-1}^2+y_{2N}^2+2 \mathring{\omega} ({x}_{2N} y_{2N-1}-{x}_{2N-1} y_{2N})+{\mathring{\omega}}^2( 1-\mathring{\iota}_N ) {x}^2_{2N-1}+{\mathring{\omega}}^2(\frac{\mathring{\iota}_N-1}{2} ){x}^2_{2N}]\\
&+ \frac{\sigma_1}{\sqrt{\varepsilon}} x_{2N-1}^3- \frac{3\sigma_1}{2\sqrt{\varepsilon}} x_{2N-1}x_{2N}^2-(\frac{\sigma_2}{{\varepsilon}} x_{2N-1}^4+ \frac{\frac{3}{8}\sigma_2}{{\varepsilon}} x_{2N}^4- \frac{3\sigma_2}{{\varepsilon}} x_{2N-1}^2 x_{2N}^2).
\end{aligned}\nonumber
\end{equation}
Because the other part $\mathring{H}$  of the Hamiltonian $H_{sim}$, by the inductive hypothesis, could be directly transformed into the Birkhoff normal form
\begin{equation}\label{Birkhoffeuler00}
\begin{aligned}
&\mathring{H}
=\textbf{i} \mathring{\omega} _0 u_0 v_0 + \sum_{k = 1}^{N-3}[\textbf{i}\mathring{\omega} _{2k-1} u _{2k-1} v _{2k-1} + \mathring{\omega} _{2k} u _{2k} v _{2k}]  -\frac{1}{2}[\mathring{\omega}_{00}(u_0 v_0)^2\\
&+\sum_{k = 1}^{N-3}\mathring{\omega}_{{2k-1},{2k-1}}(u_{2k-1} v_{2k-1})^2-\sum_{k = 1}^{N-3}\mathring{\omega}_{2k,2k}(u_{2k} v_{2k})^2+2\sum_{k = 1}^{N-3}\mathring{\omega}_{0,{2k-1}}u_0 v_0 u_{2k-1} v_{2k-1}\\
&+2\sum_{1\leq j<k \leq N-3}\mathring{\omega}_{{2j-1},{2k-1}}u_{2j-1} v_{2j-1} u_{2k-1} v_{2k-1}-2\sum_{k = 1}^{N-3}\textbf{i}\mathring{\omega}_{0,{2k}}(u_0 v_0 u_{2k} v_{2k})\\
&-2\sum_{j,k = 1}^{N-3}\textbf{i}\mathring{\omega}_{{2j-1},{2k}}(u_{2j-1} v_{2j-1} u_{2k} v_{2k})-2\sum_{1\leq j<k \leq N-3}\mathring{\omega}_{{2j},{2k}}u_{2j} v_{2j} u_{2k} v_{2k}] .
\end{aligned}\nonumber
\end{equation}

An argument similar to the one used in  the previous sections  shows that the Birkhoff normal form of  the  Hamiltonian  $H_{\varepsilon}$ is
\begin{equation}
\begin{aligned}
&\textbf{i}\mathring{\omega} _{2N-5} u _{2N-5} v _{2N-5} + \mathring{\omega} _{2N-4} u _{2N-4} v _{2N-4}  -\frac{1}{2}[\mathring{\omega}_{{2N-5},{2N-5}}(u_{2N-5} v_{2N-5})^2\\
&-\mathring{\omega}_{2N-4,2N-4}(u_{2N-4} v_{2N-4})^2
-2\textbf{i}\mathring{\omega}_{{2N-5},{2N-4}}(u_{2N-5} v_{2N-5} u_{2N-4} v_{2N-4})],
\end{aligned}\nonumber
\end{equation}
where
\begin{equation}
\begin{aligned}
&\mathring{\omega}_{{2N-5},{2N-5}}=-\frac{6144 f_{num} \mathring{\iota}_N  \left(9 \mathring{\iota}_N ^3-61 \mathring{\iota}_N ^2+127 \mathring{\iota}_N -75\right) \sigma_1^2}{f_{den} \mathring{\omega}^4 \varepsilon },
\end{aligned}\nonumber
\end{equation}
\begin{equation}
\begin{aligned}
&\mathring{\omega}_{{2N-5},{2N-4}}=\frac{9 \sqrt{2} \sigma_1^2 }{\varepsilon \mathring{\iota}_N  \sqrt{(\mathring{\iota}_N -3) \mathring{\iota}_N } \left(9 \mathring{\iota}_N ^2-34 \mathring{\iota}_N +25\right) \left(27 \mathring{\iota}_N ^2-95 \mathring{\iota}_N +50\right) \mathring{\omega}^4  }\\
&\left[162 \mathring{\iota}_N ^4-1029 \mathring{\iota}_N ^3+2457 \mathring{\iota}_N ^2-2695 \mathring{\iota}_N +825-\mathring{c} \mathring{\iota}_N  \left(81 \mathring{\iota}_N ^4-555 \mathring{\iota}_N ^3+1397 \mathring{\iota}_N ^2-1545 \mathring{\iota}_N +550\right)\right];
\end{aligned}\nonumber
\end{equation}
\begin{equation}
\begin{aligned}
&f_{num}=3 \mathring{c} \mathring{\iota}_N  [9 (9 \Delta -106) \mathring{\iota}_N ^4-(165 \Delta +1774) \mathring{\iota}_N ^3+(1177 \Delta +17576) \mathring{\iota}_N ^2-5 (827 \Delta +7161) \mathring{\iota}_N \\
&+4450 (\Delta +5)+243 \mathring{\iota}_N ^5]+[-81 (7 \Delta -80) \mathring{\iota}_N ^4+3 (363 \Delta +4430) \mathring{\iota}_N ^3-(9097 \Delta +131244) \mathring{\iota}_N ^2\\
&+(30667 \Delta +276755) \mathring{\iota}_N -36300 (\Delta +5)-1701 \mathring{\iota}_N ^5],
\end{aligned}\nonumber
\end{equation}
\begin{equation}
\begin{aligned}
&f_{den}=(\Delta -\mathring{\iota}_N +5)^2 (\Delta +\mathring{\iota}_N -5) [5 (\Delta +3)-3 \mathring{\iota}_N ] [(3 \Delta -34) \mathring{\iota}_N +5 (\Delta +5)+9 \mathring{\iota}_N ^2]^3 \\
&[(3 \Delta +34) \mathring{\iota}_N +5 (\Delta -5)-9 \mathring{\iota}_N ^2],
\end{aligned}\nonumber
\end{equation}
\begin{displaymath}
\mathring{c}=\frac{\sigma_2 \mathring{\omega}^2}{\sigma_1^2}, ~~~~~~~~~~~~~~~\Delta=\sqrt{(\mathring{\iota}_N -1) (9 \mathring{\iota}_N -25)};
\end{displaymath}
here we omit the long expression of $\mathring{\omega}_{{2N-4},{2N-4}}$ which is similar to that of $\mathring{\omega}_{{2N-5},{2N-5}}$.

Let $det(\mathring{\Omega})$ be the determinant  of the  Birkhoff normal form for the reduced Hamiltonian of $\mathring{H}$ on the center manifold, then
\begin{displaymath}
det(\Omega_{sim})=det(\mathring{\Omega}) \mathring{\omega}_{{2N-5},{2N-5}}.
\end{displaymath}
To prove $det(\Omega)\neq 0$, it suffices   to prove $f_{num}\neq 0$.

First, according to the definitions of $\sigma_2, \sigma_2$ and $ \mathring{\omega}$, it follows that $\mathring{c}$ does not depend on the scale  of the central configuration $\mathbf{e}_2(0)$ of the restricted $N$-body problem, here recall that the central configuration $\mathbf{e}_2(0)=(\xi_1, \xi_2,\cdots, \xi_{N})^\top$ is the unique solution of  the equations  of central configurations of the restricted $N$-body problem:
\begin{equation}
\sum_{j=1,j\neq k}^N \frac{m_j}{|\mathbf{r}_j-\mathbf{r}_k|^3}(\mathbf{r}_j-\mathbf{r}_k)=- \mathbf{r}_k,1\leq k\leq N,\nonumber
\end{equation}
such that $\xi_1<\xi_2<\cdots<\xi_N$ and $\lambda(\mathbf{e}_2)=1$. Hence $\mathring{c}=\frac{\sigma_2 \mathring{\lambda}}{\sigma_1^2}$ can be rewritten  as
\begin{displaymath}
\mathring{c}=\sum_{j=1}^{N-1} \frac{m_j}{(\xi_N-\xi_j)^5} /(\sum_{j=1}^{N-1}\frac{m_j}{(\xi_N-\xi_j)^4})^2.
\end{displaymath}
Thanks to (\ref{n-body3}), it follows that
\begin{displaymath}
\mathring{c}=[ \frac{1}{c_{N-1}^5 m_{N-1}^{\frac{2}{3}}}(1+o(1))] /[\frac{1}{c_{N-1}^8 m_{N-1}^{\frac{2}{3}}}(1+o(1))]=c_{N-1}^3 (1+o(1)),
\end{displaymath}
where $o(1)$ denotes infinitesimal.
Due to (\ref{n-body1}), we have
\begin{displaymath}
\mathring{c}=\frac{1}{3^{N-2}}+o(1).
\end{displaymath}

Recall that $\mathring{\iota}_N=3^{N-1}+o(1)$, as a result, we have
\begin{equation}
\begin{aligned}
&f_{num}=2 (\iota -3) [243 \iota ^4-324 \iota ^3-2310 \iota ^2+6540 \iota -3125\\
&+\sqrt{9 \iota ^2-34 \iota +25}  \left(81 \iota ^3+45 \iota ^2+883 \iota -625\right)]+o(1),
\end{aligned}\nonumber
\end{equation}
where $\iota=3^{N-1}$.

A straight forward computation shows that the function
\begin{equation}
243 \iota ^4-324 \iota ^3-2310 \iota ^2+6540 \iota -3125+\sqrt{9 \iota ^2-34 \iota +25}  \left(81 \iota ^3+45 \iota ^2+883 \iota -625\right)\nonumber
\end{equation}
is never zero for any $N\geq 3$.

Thus we can  summarize what we have proved as  the following proposition.
\begin{proposition}
For sufficiently small masses $m_k$ ($k=2,\cdots,N$),  the reduced Hamiltonian  on the center manifold near every Euler-Moulton relative equilibrium is nondegenerate for  the corresponding planar $N$-body problem.
\end{proposition}

By this proposition, it's clear that
 Theorem \ref{Conley and Zenderapply3} holds.

\indent\par


\newpage

\bibliographystyle{plain}



\end{document}